\documentclass[reqno,letter]{amsart}
\title[The meromorphic $R$--matrix of the Yangian]{The meromorphic $R$--matrix of the Yangian}

\author[S. Gautam]{Sachin Gautam}
\address{Department of Mathematics, The Ohio State University, 231 W 18th St.
Columbus, OH (USA)}
\email{gautam.42@osu.edu}

\author[V. Toledano Laredo]{Valerio Toledano Laredo}
\address{Department of Mathematics, Northeastern University, 360 Huntington Avenue, Boston, MA (USA)}
\email{V.ToledanoLaredo@neu.edu}

\author[C. Wendlandt]{Curtis Wendlandt}
\address{Department of Mathematics, The Ohio State University, 231 W 18th St.
Columbus, OH (USA)}
\email{Wendlandt.4@osu.edu}
\dedicatory{To Kolya Reshetikhin, on his 60th birthday.}

\usepackage{amsmath,amssymb}
\usepackage{mathrsfs}
\usepackage[usenames,dvipsnames]{color}
\usepackage{enumitem}
\usepackage{hyperref}
\hypersetup{
    colorlinks=true,       			% false: boxed links; true: colored links
    linkcolor=blue, %Cerulean,          	% color of internal links
    citecolor=blue,        			% color of links to bibliography
    filecolor=blue,      				% color of file links
    urlcolor=blue           			% color of external links
}

\usepackage[all]{xy}
\xyoption{arc}
\newtheorem*{thm}{Theorem}
\newtheorem*{prop}{Proposition}
\newtheorem*{lem}{Lemma}
\newtheorem*{cor}{Corollary}

\newenvironment{pf}{\paragraph{{\sc Proof}}}{\qed\par\medskip}
\theoremstyle{definition}

\newtheorem*{rem}{Remark}

\numberwithin{equation}{section}

\numberwithin{figure}{section}

\newcommand {\remark}{\noindent {\bf Remark. }}

\newcommand{\Y}{\mathcal{Y}}

\newcommand{\spec}{\sigma}

\newcommand{\ForRes}[2]{\oint #1\, d#2}

\newcommand{\id}{\mathbf{1}}

\newcommand{\veps}{\varepsilon}
\newcommand{\lp}{\left(}
\newcommand{\rp}{\right)}

\newcommand{\lX}{\mathfrak{X}}

\newcommand{\g}{\mathfrak{g}}
\newcommand{\h}{\mathfrak{h}}

\newcommand{\gn}{\mathfrak{n}}
\newcommand{\Sym}{\mathfrak{S}}

\newcommand{\vup}[1]{v_+}%{\left|\uparrow_{#1}\right\rangle}
\newcommand{\vdown}[1]{v_-}%{\left|\downarrow_{#1}\right\rangle}

\newcommand{\bfA}{\mathbf{A}}

\newcommand{\A}{\mathcal{A}}
\newcommand{\cC}{\mathcal{C}}
\newcommand{\D}{\mathcal{D}}

\newcommand{\HH}{\mathcal{H}}

\newcommand{\J}{\mathcal{J}}
\newcommand{\K}{\mathcal{K}}
\newcommand{\LL}{\mathcal{L}}

\newcommand{\RR}{\mathcal{R}}

\newcommand{\V}{\mathcal{V}}

\newcommand{\calZ}{\mathcal Z}
\newcommand{\calP}{\mathcal P}
\newcommand{\C}{\mathbb{C}}
\newcommand{\nC}{\mathbb{C}^{\times}}
\newcommand{\N}{\mathbb{Z}_{\geq 0}}
\newcommand{\IP}{\mathbb{P}}
\newcommand{\Q}{\mathbb{Q}}

\newcommand{\Z}{\mathbb{Z}}

\renewcommand{\Re}{\operatorname{Re}}

\newcommand{\opp}{{\scriptstyle{\operatorname{op}}}}

\newcommand{\Rloopintro}
{\operatorname{Rep}_{\scriptscriptstyle{\operatorname{fd}}}(\qloop)}

\newcommand {\wh}[1]{\widehat{#1}}
\newcommand {\ol}[1]{\overline{#1}}

\newcommand {\Aut}{\operatorname{Aut}}
\newcommand {\Ker}{\operatorname{Ker}}
\newcommand{\End}{\operatorname{End}}
\newcommand{\Hom}{\operatorname{Hom}}
\newcommand{\ad}{\operatorname{ad}}
\newcommand{\Ad}{\operatorname{Ad}}

\newcommand{\qKZ}{q\operatorname{KZ}}
\newcommand {\aand}{\qquad\text{and}\qquad}

\newcommand {\Rep}{\operatorname{Rep}}
\newcommand {\ffd}{\operatorname{fd}}
\newcommand {\fd}{finite--dimensional }
\newcommand {\lhs}{left--hand side }
\newcommand {\rhs}{right--hand side }
\newcommand {\wrt}{with respect to }

\newcommand{\Ryang}{\Rep_{\ffd}(\Yhg)}
\newcommand{\Rync}{\Rep_{\ffd}^{\scriptscriptstyle{\operatorname{NC}}}(\Yhg)}

\newcommand{\ds}{\displaystyle}

\newcommand {\Comment}[1]{}
\newcommand {\Omit}[1]{}

\newcommand {\gr}{\operatorname{gr}}

\newcommand {\qloop}{U_q(L\g)}

\newcommand {\Yhg}{Y_\hbar(\g)}

\newcommand{\Yhsl}[1]{Y_{\hbar}(\sl_{#1})}

\newcommand {\isom}{\stackrel{\sim}{\rightarrow}}

\newcommand{\cbin}[2]{\left(\begin{array}{c} #1\\ #2\end{array}\right)}

\renewcommand {\sl}{\mathfrak{sl}}

\newcommand {\eg}{{\it e.g., }}

\newcommand {\Id}{\operatorname{Id}}
\newcommand {\ie}{{\it i.e.}, }

\newcommand {\KD}{Drinfeld--Kohno }

\newcommand{\bfI}{{\mathbf I}}

\newcommand {\sfr}{\mathsf{r}}

\usepackage{upgreek}

\newcommand{\kmDelta}{\Delta}
\newcommand{\kmtensor}[1]{\underset{\scriptscriptstyle{#1}}{\otimes}}
\newcommand{\kmdelta}[1]{\Delta_{#1}}

\newcommand{\ddelta}[1]{\underset{\scriptscriptstyle{{\operatorname{D}},#1}}{\Delta}}

\newcommand{\dtensor}[1]{\underset{\scriptscriptstyle{\operatorname{D},#1}}{\otimes}}

\newcommand{\fop}{\sfr}%{\mathsf{r}^-}
\newcommand{\fopop}{\sfr^{21}}%{\mathsf{r}^+}
\newcommand{\Fop}[1]{\mathcal{R}^-_{#1}}
\newcommand{\UF}{\mathcal{R}^-}

\newcommand {\updown}{\uparrow/\downarrow}
\newcommand {\downup}{\downarrow/\uparrow}

\newcommand{\Rup}{\mathcal{R}^{\uparrow}}
\newcommand{\Rdown}{\mathcal{R}^{\downarrow}}
\newcommand{\UR}{\mathcal{R}^{\varepsilon}}
\newcommand{\URupdown}{\mathcal{R}^{\updown}}

\newcommand{\Rrat}{\mathsf{R}}

\newcommand{\Rratd}{\mathsf{R}^0}

\newcommand{\Rupd}{\mathcal{R}^{0,\uparrow}}
\newcommand{\Rdownd}{\mathcal{R}^{0,\downarrow}}
\newcommand{\URd}{\mathcal{R}^{0,\varepsilon}}
\newcommand{\URdupdown}{\mathcal{R}^{0,\updown}}

\newcommand{\logFop}{\omega}

\newcommand{\Top}[1]{\operatorname{T}(#1)}

\newcommand{\ellg}{\ell}
\newcommand{\filt}[2]{\mathcal{F}_{#1}(#2)}

\newcommand{\hext}[2]{{#1}[\negthinspace[#2]\negthinspace]}
\newcommand{\Linfinity}[2]{{#1}[#2;{#2}^{-1}]\negthinspace]}
\newcommand{\Ltwo}[3]{#1[#2;{#2}^{-1}]\negthinspace][#3;{#3}^{-1}]\negthinspace]}

\newcommand{\Ynaughthg}{Y^0_{\hbar}(\mathfrak{g})}
\newcommand{\Ypmhg}{Y^{\pm}_{\hbar}(\mathfrak{g})}

\newcommand{\Yplushg}{Y^+_{\hbar}(\mathfrak{g})}
\newcommand{\Yminushg}{Y^-_{\hbar}(\mathfrak{g})}

\newcommand{\Primitive}[2]{\operatorname{Prim}^{#2}(#1)}
\newcommand{\primitive}[1]{#1\otimes 1 + 1\otimes #1}

\newcommand{\adT}{\mathcal{T}}

\newcommand{\Triv}{{\mathbf 1}}

\newcommand{\scR}{\mathscr{R}}

\newcommand {\sfGamma}{\mathsf{\Gamma}}

\newcommand {\curtiscomment}[1]{}%{\footnote{\textcolor{magenta}{C:\,#1}}}
\newcommand {\sfR}{\mathsf{R}}
\newcommand {\calR}{\RR}

\newcommand {\sop}{^{\scriptscriptstyle{\operatorname{op}}}}
\newcommand {\YBE}{Yang--Baxter equation }
\newcommand {\KMA}{Kac--Moody algebra }
\newcommand {\KT}{Khoroshkin--Tolstoy }
\newcommand {\QYBE}{quantum Yang--Baxter equation }
\newcommand {\RRD}{\RR}%{\RR^{(D)}}
\newcommand {\DYg}{D\Yhg}
\newcommand {\ud}{{\uparrow/\downarrow}}
\newcommand {\RRou}[1]{\RR^{0,\uparrow}_{#1}}%{\RR^{0,\ud}_{#1}}
\newcommand {\RRod}[1]{\RR^{0,\downarrow}_{#1}}%{\RR^{0,\ud}_{#1}}
\newcommand{\Rq}{\scR}
\newcommand {\flip}{(1\,2)}

\newcommand {\oT}{\operatorname{T}}
\newcommand {\bin}[2]{\lp\begin{matrix}#1\\#2\end{matrix}\rp}

\newcommand {\sfQ}{{\mathsf Q}}

\begin{document}
\maketitle

\begin{abstract}
Let $\g$ be a complex semisimple Lie algebra and $\Yhg$ its Yangian.
Drinfeld proved that the universal $R$--matrix $\RR(s)$ of $\Yhg$ gives rise to rational
solutions of the $\operatorname{QYBE}$ on irreducible, finite--dimensional representations
of $\Yhg$. This result was recently extended by Maulik--Okounkov to symmetric
Kac--Moody algebras, and representations arising from geometry.
We show that rationality ceases to hold on arbitrary finite-dimensional
representations, if one requires such solutions to be natural
and  compatible with tensor products. Equivalently, the tensor category of
finite--dimensional representations of $\Yhg$ does not admit rational commutativity constraints.
We construct instead two meromorphic commutativity constraints, which are
related by a unitarity condition. Each possesses an asymptotic expansion in $s$
which has the same formal properties as $\RR(s)$, and therefore coincides
with it by uniqueness. In particular, we give a
constructive proof of the existence of $\RR(s)$.
Our construction relies on the Gauss decomposition
$\RR^+(s)\cdot \RR^0(s)\cdot \RR^-(s)$ of $\RR(s)$. The
divergent abelian term $\RR^0$ was resummed on \fd
representations by the first two authors in \cite{sachin-valerio-III}.
In the present paper, we construct $\RR^{\pm}(s)$, prove that they
are rational on finite--dimensional representations, and that they
intertwine the standard coproduct of $\Yhg$ and the deformed Drinfeld coproduct
introduced in \cite{sachin-valerio-III}.
\end{abstract}

\setcounter{tocdepth}{1}
\tableofcontents
\newpage

\section{Introduction}\label{sec: intro}
%==============

\subsection{} % intro and universal R
%--------------

Let $\g$ be a complex, semisimple Lie algebra with an invariant
inner product $(\cdot,\cdot)$, and $\Yhg$ the corresponding Yangian,
which is a Hopf algebra deforming the current algebra $U(\g[z])$
introduced by Drinfeld \cite{drinfeld-qybe}. We assume that $\hbar
\in\nC$ is fixed throughout. Drinfeld proved that $\Yhg$ possesses
a unique universal $R$--matrix. Specifically, let $\Delta:\Yhg\to\Yhg
\otimes \Yhg$ be the coproduct of $\Yhg$, and $\tau_s:\Yhg\to\Yhg$
%,$s\in\C$, 
the one--parameter group of automorphisms which quantizes
the shift automorphism $z\mapsto z+s$ of $U(\g[z])$. Note that, if
$s$ is considered %regarded
as a variable, $\tau_s$ may be regarded as
%is
a homomorphism $\Yhg\to
\Yhg[s]$. Then, the following holds.

\begin{thm}\cite[Thm. 3]{drinfeld-qybe}
\label{thm:intro1}\hfill
\begin{enumerate}[font=\upshape]
%
% existence and uniqueness
\item
There is a unique formal series 
\[\RRD(s)
=1 + \sum_{k=1}^{\infty} \RRD_k s^{-k}\in\hext{\Yhg^{\otimes 2}}{s^{-1}}\]
such that the following holds\footnote
{Our conventions differ slightly from those of \cite{drinfeld-qybe}, where
the intertwining equation \eqref{eq:delta delta opp} is written as $\id\otimes\tau_s\circ\Delta\sop(a)
=\RR(s)^{-1}\cdot\id\otimes\tau_s\circ\Delta(a)\cdot\RR(s)$. Thus, our
$\RRD(s)$ is Drinfeld's $\RR(-s)^{-1}$.} in $\Linfinity{\Yhg^{\otimes 2}}{s}$
\begin{equation}\label{eq:delta delta opp}
\tau_s\otimes\id\circ \Delta\sop(a)
= 
\RRD(s)\cdot
\tau_s\otimes\id\circ \Delta(a)
\cdot
\RRD(s)^{-1}
\end{equation}
for any $a\in\Yhg$, and
\begin{align}
\Delta\otimes\id(\RRD(s)) &= \RRD_{13}(s)\cdot\RRD_{23}(s)
\label{eq:cable 1}\\
\id\otimes\Delta(\RRD(s)) &= \RRD_{13}(s)\cdot\RRD_{12}(s)
\label{eq:cable 2}
\end{align}
%
% additional properties
\item $\RRD$ satisfies the following identities
\[\begin{array}{ll}
\text{$\bullet$ 1--jet:}		& \RRD(s) = 1+ \hbar s^{-1} \Omega_{\g} + O(s^{-2})\\%\\[.1ex]
\text{$\bullet$ Unitarity:}	&\RRD(s)^{-1} = \RRD_{21}(-s)\\%\\[.01ex]
\text{$\bullet$ Translation:}&\tau_a\otimes \tau_b(\RRD(s)) = \RRD(s+a-b)%\\
\end{array}\]
where $\Omega_{\g}\in\g\otimes\g$ is the Casimir tensor corresponding
to $(\cdot,\cdot)$.
% QYBE
\item $\RRD$ is a solution of the \QYBE (QYBE)\footnote{The QYBE may be viewed as an equation in $\hext{\Linfinity
{\Yhg^{\otimes 3}}{s_1}}{s_2^{-1}}$ by expanding it as if $|s_2|\gg |s_1|$,
that is setting $(s_1+s_2)^{-1}=\sum_{k\geq 0}(-1)^{k}s_1^ks_2^{-k-1}$,
or as an equation in $\hext{\Linfinity{\Yhg^{\otimes 3}}{s_2}}{s_1^{-1}}$,
by setting $(s_1+s_2)^{-1}=\sum_{k\geq 0}(-1)^{k}s_2^ks_1^{-k-1}$. The
precise statement of (3) above is that \eqref{eq:QYBE} holds in either of
these cases.}
\begin{equation}\label{eq:QYBE}
\RRD_{12}(s_1)\RRD_{13}(s_1+s_2)\RRD_{23}(s_2)
= 
\RRD_{23}(s_2)\RRD_{13}(s_1+s_2)\RRD_{12}(s_1)
\end{equation}
\end{enumerate}
\end{thm}

\subsection{}\label{ss:rational soln} % Drinfeld factorisation for irreducible representations rational QYBE
%--------------

One of the main goals of this paper is to clarify the analytic nature
of the formal power series $\RRD(s)$, and of the solutions of the
QYBE obtained from it. Let $V_1,V_2$ be two \fd representations
of $\Yhg$, and $\RRD_{V_1,V_2}(s)\in\hext{\End(V_1\otimes V_2)}
{s^{-1}}$ the corresponding evaluation of $\RRD(s)$. Drinfeld proved
that $\RRD_{V_1,V_2}(s)$ has a zero radius of convergence in general
\cite[Examples  1,2]{drinfeld-qybe}, but nevertheless gives rise to a
rational solution of the QYBE as follows.
\begin{thm}\cite[Thm. 4]{drinfeld-qybe}\label{thm:intro2}
Assume that $V_1$ and $V_2$ are irreducible with highest weight
vectors $v_1,v_2$, and let $\rho_{V_1,V_2}(s)\in 1+s^{-1}\hext{\C}{s^{-1}}$
be the matrix element of $\RRD_{V_1,V_2}(s)$ given by
\[\RRD_{V_1,V_2}(s)\,v_1\otimes v_2=\rho_{V_1,V_2}(s)\,v_1\otimes v_2\]
Then,
\begin{enumerate}
\item $\sfR_{V_1,V_2}(s)=\RRD_{V_1,V_2}(s)\cdot \rho_{V_1,V_2}(s)^{-1}$
is a rational function of $s$.
% rational QYBE
\item If $V_1=V=V_2$, \eqref{eq:QYBE} together with the factorisation
\begin{equation}\label{eq:factorisation of D}
\RRD_{V_1,V_2}(s) = \sfR_{V_1,V_2}(s)\cdot \rho_{V_1,V_2}(s)
\end{equation}
imply that $\sfR_{V,V}(s)$ is a rational solution of the QYBE.
\end{enumerate}
\end{thm}
More recently, a geometric construction of $R$--matrices % solutions of the QYBE
corresponding to the (extended) Yangian of a symmetric \KMA was given
by Maulik--Okounkov \cite{maulik-okounkov-qgqc}, which provides in
particular an alternative construction of rational solutions of the QYBE
on the equivariant cohomology of Nakajima quiver varieties.

\subsection{}
%--------------

One of the byproducts of this paper is to extend the factorisation
\eqref{eq:factorisation of D}
%\[\RRD_{V_1,V_2}(s) = \sfR_{V_1,V_2}(s)\cdot \rho_{V_1,V_2}(s)\]
to an arbitrary pair of (not necessarily irreducible) \fd representations.
In this case, the divergent factor $\rho_{V_1,V_2}(s)$ takes values
in $\End(V_1\otimes V_2)[\![s^{-1}]\!]$, and intertwines the action of
$\Yhg$ given by $\Delta_s=\tau_s\otimes\id\circ\Delta$, whereas
the rational factor $\sfR_{V_1,V_2}(s)$ intertwines $\Delta_s$ and
$\Delta_s^\opp=\tau_s\otimes\id\circ\Delta^\opp$. However, since
$\rho_{V_1,V_2}(s)$ is not scalar--valued in general, it is not clear
whether $\sfR_{V_1,V_2}(s)$ satisfies the QYBE when $V_1=V_2$.

We prove in fact that, even for $\g=\sl_2$, no rational intertwiner
$\sfR_{V_1,V_2}(s)\in\End(V_1\otimes V_2)$ exists which is defined 
for any $V_1,V_2\in\Ryang$, is natural in $V_1$ and $V_2$, and
satisfies the cabling identities \eqref{eq:cable 1}--\eqref{eq:cable 2}.
Equivalently, the tensor category of \fd representations of $\Yhg$
does not admit rational commutativity constraints. In particular,
this raises the question of whether one can consistently define
rational solutions of the QYBE on all \fd representations of $\Yhg
$.\footnote{The Maulik--Okounkov construction mentioned in \ref
{ss:rational soln} provides a partial solution to this question, since
an arbitrary representation of $\Yhg$ may not have a geometric
realisation.}

\subsection{}\label{ss:alt sol}
%--------------

In the present paper, we propose an alternative solution to this
issue, by constructing {\it meromorphic} commutativity constraints
on $\Ryang$, and in particular consistent meromorphic solutions of
the QYBE on all $V\in\Ryang$. Namely, we prove that the universal
$R$--matrix of $\Yhg$, while generally divergent on a pair of \fd
representations $V_1,V_2$, can be canonically {\em resummed},
in two distinct ways. This yields a {\em pair} of meromorphic functions 
\[\RR^{\uparrow}_{V_1,V_2}(s),
\RR^{\downarrow}_{V_1,V_2}(s):
\C\to\End(V_1\otimes V_2)\]
which are natural \wrt $V_1,V_2$, satisfy the intertwining relation
\eqref{eq:delta delta opp}, the cabling identities \eqref{eq:cable 1}--\eqref
{eq:cable 2}, as well as the translation property. The function
$\RR^{\uparrow}_{V_1,V_2}(s)$ (resp. $\RR^ {\downarrow}_{V_1,
V_2}(s)$) is asymptotic to $\RRD_{V_1,V_2}(s)$ as $s\to\infty$
with $\Re(s/\hbar)>0$ (resp. $\Re(s/\hbar)<0$), and is related to
$\RR^ {\downarrow}_{V_1,V_2}(s)$ by
the unitarity relation
\[\RR_{V_1,V_2}^{\uparrow}(s)^{-1}=\RR_{V_2,V_1}^{\downarrow}(-s)^{21}\]

The situation is somewhat analogous to the case of the quantum loop
algebra $\qloop$. In that case, if $\Rq\in\qloop\wh{\otimes}\qloop$ is
the universal $R$--matrix, then 
\[\Rq^\infty(z)=\tau_z\otimes\id (\Rq)\in \hext{\qloop^{\otimes 2}}{z^{-1}}
\quad\text{and}\quad
\Rq^0(z)=\id\otimes\tau_z(\Rq)\in \hext{\qloop^{\otimes 2}}{z}\]
converge, near $z=\infty$ and $z=0$ respectively, to meromorphic
functions of $z\in\C^\times$ on the tensor product $\V_1\otimes\V_2$
of any two \fd representations \cite{etingofmoura-KL,kazhdan-soibelman},
and are related by $\Rq_{\V_1,\V_2}^\infty(z)^{-1}=\Rq_{\V_2,\V_1}^0
(z^{-1})^{21}$. In the case of $\qloop$, however, $\Rq^\infty(z)$ and
$\Rq^0(z)$ are convergent as is, and do not need to be resummed.

\subsection{}
%--------------

Our approach does not rely on Drinfeld's cohomological construction
of $\RRD(s)$ to carry out the resummation. It produces the functions
$\RR^{\uparrow/\downarrow}_{V_1,V_2}(s)$ through a direct, explicit
construction, which shows in particular that they have an asymptotic
expansion as $s\to\infty$. The fact that the latter coincides with $\RR
_{V_1,V_2}(s)$, and therefore {\it a posteriori} that $\RR_{V_1,V_2}(s)$
can be resummed, follows from the fact that the asymptotic expansion
can be lifted to $\hext{\Yhg^{\otimes 2}}{s^{-1}}$, and shown to have
the properties which uniquely determine $\RRD(s)$ by Theorem \ref
{thm:intro1}. In particular, our construction yields an independent, and
constructive proof of the existence of $\RR(s)$.

\subsection{}
%--------------

Our construction can be motivated by the following considerations.
The $R$--matrix of $\Yhg$ is expected to arise as the canonical
element in $\DYg\wh{\otimes}\DYg$, where $\DYg\supset\Yhg$
is the double Yangian of $\g$, which is a quantisation of the graded
Drinfeld double
\[\lp\g[z^{\pm 1}],\g[z],z^{-1}\g[z^{-1}]\rp\]
of $\g[z]$. Although a detailed understanding of $\DYg$ is still
lacking at present (see, however, \cite{khoroshkin-tolstoy} and
\cite{curtis-double}), this
suggests that, given a triangular decomposition $\g=\gn_+\oplus
\h\oplus\gn_-$ of $\g$, $\RRD(s)$ should have a corresponding
Gauss decomposition
\begin{equation}\label{eq:Gauss}
\RRD(s)=\RRD^+(s)\cdot\RRD^0(s)\cdot\RRD^-(s)
\end{equation}
where $\RRD^0(s)$ quantises the canonical
element in $\h[z]\wh{\otimes}z^{-1}\h[z^{-1}]$, and $\RR^\pm(s)$
those in $\gn_\pm[z]\wh{\otimes}z^{-1}\gn_\mp[z^{-1}]$ respectively.
Moreover, the unitarity of $\RR(s)$ suggests that
\[\RRD^0(s)^{-1} = \RRD^{0}(-s)^{21}
\aand
\RRD^+(s)^{-1} = \RRD^{-}(-s)^{21}\]
Accordingly, we construct each factor $\RRD^0(s),\RRD^-(s),
\RRD^+(s)= (\RRD^{-}(-s)^{21})^{-1}$, and their resummation on
\fd representations separately.

\subsection{} % resummation of R0 
%--------------

\KT gave a heuristic formula for $\RRD^0$ \cite{khoroshkin-tolstoy},
as the exponential of an infinite sum in the abelian subalgebra of
$\DYg$ which quantises $\h[z,z^{-1}]$. In \cite{sachin-valerio-III},
the first two named authors gave a precise version of this formula,
where the exponent takes values in the abelian subalgebra $\Ynaughthg$ of $\Yhg$ which
quantises $\h[z]$. We showed moreover that this expression can
be resummed on a tensor product $V_1\otimes V_2$ of \fd representations in
two different ways. This yields two meromorphic functions
$\RRou{V_1,V_2}(s),\RRod{V_1,V_2}(s)$, which have the same
asymptotic expansion on $\pm\Re(s/\hbar)>0$, and are related
by $\RRou{V_1,V_2}(s)^{-1}=\RRod{V_2,V_1}(-s)^{21}$.

\subsection{} % Drinfeld tensor product 
%--------------

An important discovery of \cite{sachin-valerio-III} is that these abelian
$R$--matrices play a similar role to that of the full $R$--matrix of $\Yhg$,
but \wrt the {\it deformed Drinfeld tensor product}. The latter was introduced
in \cite{sachin-valerio-III} by degenerating the Drinfeld tensor product of
the quantum loop algebra introduced by Hernandez
\cite{hernandez-drinfeld}. It gives rise to a family of actions of $\Yhg$ on
the vector space $V_1\otimes V_2$, which is denoted by $V_1\dtensor{s}V_2$ and is a rational
function of a parameter $s\in\C$.  The tensor product $\dtensor{s}$
is associative, in that the identification of vector spaces
\[(V_1\dtensor{s_1}V_2)\dtensor{s_2}V_3 = 
V_1\dtensor{s_1+s_2}(V_2\dtensor{s_2}V_3)\]
intertwines the action of $\Yhg$ for any $s_1,s_2\in\C$, and
endows $\Ryang$ with the structure of a meromorphic tensor
category in the sense of \cite{soibelman-mero,soibelman-meromorphic}. 

The endomorphisms $\RR^{0,\ud}_{V_1,V_2}$ are meromorphic
commutativity constraints \wrt $\dtensor{s}$, that is they satisfy the
representation theoretic version of the identities (1) of Theorem
\ref{thm:intro1}.
In the present paper, %Section \ref{sec: two-tensors} and \ref{sec:fullRd},
we complement the results of \cite{sachin-valerio-III} by lifting $\dtensor{s}$
to a {\it deformed Drinfeld coproduct}
\[\ddelta{s}:\Yhg\to\Linfinity{\lp\Yhg\otimes\Yhg\rp}{s}\]
and the common asymptotic expansion of $\RR^{0,\ud}_{V_1,V_2}(s)$
to an element 
\[\RR^0(s)\in\hext{\lp\Yhg\otimes\Yhg\rp}{s^{-1}}\]
which satisfy the identities (1) of Theorem \ref{thm:intro1}, with $\tau
_s\otimes\id\circ\Delta$ replaced by $\ddelta{s}$, and $\RR(s)$ by $
\RR^0(s)$.

\subsection{} % expected: R- and Drinfeld to KM
%--------------

The central ingredient of the present paper is the construction of $
\RR^\pm(s)$, which is based on the following. The fact that $\RR(s)$
(resp. $\RR^0(s)$) conjugates the standard coproduct $\Delta_{s}=
\tau_s\otimes\id\circ\Delta$ (resp. the deformed Drinfeld coproduct
$\ddelta{s}$) to its opposite, together with the Gauss decomposition
\eqref{eq:Gauss}, suggest that $\RR^-(s)$ should conjugate the standard
coproduct $\Delta_{s}$ to the deformed Drinfeld coproduct $\ddelta
{s}$. This is consistent with the fact that an analogous statement
holds for the quantum loop algebra, and the related construction
of twists conjugating quantum coproducts corresponding to different
polarisations of a Manin triple given in \cite{ekp}. In this case, the
standard (resp. Drinfeld) coproducts on $\qloop$ correspond, respectively,
to the polarisations
\[\g[z]\oplus z^{-1}\g[z^{-1}]
=
\g[z^{\pm 1}]=
\lp\gn_-[z^{\pm 1}]\oplus\h[z]\rp
\oplus
\lp z^{-1}\h[z^{-1}]\oplus\gn_+[z^{\pm 1}]\rp\]

\subsection{} % construct R_- as a KM to D intertwiner
%--------------

We prove that this intertwining property uniquely determines
an element $\RR^-(s)$, provided it is required to lie in $\hext
{\lp\Yminushg\otimes\Yplushg\rp}{s^{-1}}$ and have constant
term 1, where $\Ypmhg\subset\Yhg$ is the subalgebra deforming
$U(\gn_\pm[z])$. We show in fact that, under this triangularity
assumption, $\RR^-(s)$ is uniquely determined by the requirement
that it intertwines the standard and Drinfeld coproducts of the
loop generators $t_{i,0},t_{i,1}$ of $\Yhg$ which deform $\h
\oplus\h\otimes z\subset \h[z]$.

We then show that, for any $V_1,V_2\in\Ryang$, $\RR^-_{V_1,
V_2}(s)$ is a rational function of $s$.
%, and in particular requires no resummation.
Moreover, the following cocycle identity holds
for any $V_1,V_2,V_3\in\Ryang$
\begin{equation}\label{eq:cocycle intro}
\RR^-_{V_1\dtensor{s_1}V_2,V_3}(s_2)
\cdot 
\RR^-_{V_1,V_2}(s_1)
=
\RR^-_{V_1,V_2\dtensor{s_2}V_3}(s_1+s_2)
\cdot \RR^-_{V_2,V_3}(s_2)
\end{equation}
%for any triple $V_1,V_2,V_3\in\Ryang$.
%

Together with the identities satisfied by $\RR^0(s)$, this guarantees
that the product $\RR(s)=\RR^+(s)\cdot\RR^0(s)\cdot\RR^-(s)$,
where $\RR^+(s)=(\RR^-(-s)^{21})^{-1}$, satisfies the identities
\eqref{eq:delta delta opp}--\eqref{eq:cable 2} on any pair of \fd
representations. A separation of points argument %due to Drinfeld
then implies that $\RR(s)$ satisfies Drinfeld's uniqueness criterion of the
universal $R$--matrix of $\Yhg$, and in particular coincides with
it.\footnote{The passage to \fd representations is dictated by the
fact that the cocycle identity \eqref{eq:cocycle intro} does not 
appear to have a natural lift to $\Yhg$. Indeed, when lifted to
$\Yhg$, the \lhs lies in $\hext{\Linfinity
{\Yhg^{\otimes 3}}{s_1}}{s_2^{-1}}$, while the \rhs lies
in $\hext{\Linfinity
{\Yhg^{\otimes 3}}{s_2}}{s_1^{-1}}$.}

Finally, since $\RR^-_{V_1,V_2}(s)$ is a rational function of $s$,
the product
\begin{equation}\label{eq:gauss intro}
\RR^\ud_{V_1,V_2}(s)=
\RR^+_{V_1,V_2}(s)\cdot\RR^{0,\ud}_{V_1,V_2}(s)\cdot\RR^-_{V_1,V_2}(s)
\end{equation}
is a resummation of $\RR_{V_1,V_2}(s)$, as well a meromorphic
commutativity constraint on $\Ryang$ \wrt the standard tensor
product.

\subsection{} % meromorphic tensor categories.
%--------------

Our results may be rephrased as follows. As proved in \cite{sachin-valerio-III},
and mentioned above, \fd representations of $\Yhg$, together with the
deformed Drinfeld tensor product $\dtensor{s}$ and one of the resummed
abelian $R$--matrices $\RR^{0,\ud}(s)$ is a meromorphic braided tensor
category.

Similarly, $\Ryang$ endowed with the deformed standard tensor product
$\otimes_s=\otimes\circ(\tau_s^*\otimes\id)$ is a meromorphic (in fact,
polynomial) tensor category. Our construction of the resummed $R$--matrices 
$\RR^{\uparrow/\downarrow}(s)$
endows this category with a meromorphic braiding.\footnote{
An analogous statement was proved by Kazhdan--Soibelman for the
quantum loop algebra in \cite{kazhdan-soibelman}. As pointed out in
\ref{ss:alt sol}, however, in the case of $\qloop$ no resummation of 
the universal $R$--matrix of $\qloop$ is needed.} Moreover, the
element $\RR^-(s)$ is a rational braided tensor structure on the
identity functor
\[(\Ryang,\dtensor{s},\RR^{0,\ud})\to(\Ryang,\otimes_{s},\RR^{\ud})\]

That is, $\RR^-(s)$ gives rise to a system of natural isomorphisms of
$\Yhg$--modules $\RR^-_{V_1,V_2}(s):V_1\otimes_{s}V_2\to V_1
\dtensor{s} V_2$, which is compatible with the (trivial) associativity
constraints and the meromorphic braidings, \ie such that the following
diagrams commute for any $V_1,V_2,V_3\in\Ryang$
% associativity constraints
\[\xymatrix@R=1.4cm@C=.6cm{
(V_1\otimes_{s_1} V_2) \otimes_{s_2} V_3
\ar[d]_{\RR^-_{V_1,V_2}(s_1)\otimes\id_{V_3}}
\ar@{=}[rr]&&
V_1\otimes_{s_1+s_2} ( V_2 \otimes_{s_2} V_3)
\ar[d]^{\id_{V_1}\otimes \RR^-_{V_2,V_3}(s_2)}\\
( V_1\dtensor{s_1} V_2) \otimes_{s_2} V_3
\ar[d]_{\RR^-_{V_1\dtensor{s_1}V_2,V_3}(s_2)}&&
V_1\otimes_{s_1+s_2} ( V_2 \dtensor{s_2} V_3)
\ar[d]^{\RR^-_{V_1,V_2\dtensor{s_2}V_3}(s_1+s_2)}\\
( V_1\dtensor{s_1} V_2) \dtensor{s_2} V_3
\ar@{=}[rr]&&
V_1\dtensor{s_1+s_2} ( V_2 \dtensor{s_2} V_3)
}\]
as dictated by the cocycle equation \eqref{eq:cocycle intro}, and
% commutativity constraints
\[\xymatrix@R=1.4cm{
V_1(s)\otimes V_2 
\ar[rr]^{\flip\circ\RR^{\ud}_{V_1,V_2}(s)}
\ar[d]_{\RR^-_{V_1,V_2}(s)}
&&
V_2\otimes V_1(s)
\ar[d]^{\RR^-_{V_2,V_1}(-s)}\\
V_1(s)\dtensor{0} V_2
\ar[rr]_{\flip\circ\RR^{0,\ud}_{V_1,V_2}(s)}
&&
V_2\dtensor{0} V_1(s)
}\]
which follows from the Gauss decomposition \eqref{eq:gauss intro},
together with the fact that $(\RR^-(-s)^{21})^{-1}=\RR^+(s)$.

\subsection{Outline of the paper}
%------------------

We review the definition of $\Yhg$ in Section \ref{sec: yangians},
and that of the standard and Drinfeld coproducts in
Section \ref{sec: two-tensors}.
In Section \ref{sec:fusion}, we prove the existence and uniqueness
of $\UF(s)$, and establish its various properties.
In Section \ref{sec:fusion-sl2}, we give
an explicit expressions for $\UF(s)$ when $\g=\sl_2$.
In Section \ref{sec:fullRd}, we review the construction of
$\RR^{0,\ud}(s)$ given in \cite{sachin-valerio-III}. We then explicitly
lift its asymptotic expansion to $\hext{\Ynaughthg^{\otimes 2}}
{s^{-1}}$, and prove that it satisfies properties analogous to
those of Drinfeld's $R$--matrix, but \wrt the deformed Drinfeld
coproduct.
We also prove that there is
no rational commutativity constraint on $\Ryang$.
Combined with the results of Section \ref{sec:fusion}, 
we obtain the same assertions
for the standard tensor product in
Section \ref{sec:fullR}.
We give a proof of 
the uniqueness of the universal $R$--matrix of the Yangian
in Appendix
\ref{asec:UR-unique}, thus completing the proof
that our construction gives rise to Drinfeld's $R$--matrix.
In Section \ref{sec:merocats}, we restate our results
in the language of meromorphic tensor categories.
In the final Section \ref{se:qloop},
we discuss the analogous case of the quantum loop algebra,
and relate the two by means of the meromorphic tensor functor
constructed in \cite{sachin-valerio-III}. Appendix \ref{asec:sep}
contains a proof due to Drinfeld that \fd representations
separate points of $\Yhg$.

\subsection{Acknowledgments}
%---------------------------------------

We would like to thank Pavel Etingof for several helpful discussions
about $\qKZ$ equations and $R$--matrices. We are also grateful to
Maria Angelica Cueto for helping us with the combinatorial aspects
of the paper.
%of Appendix \ref{asec:higher}.
The first author was supported through the Simons foundation
collaboration grant 526947.
The second author was supported through the NSF grant DMS--1802412.
The third author was supported by an NSERC CGS D graduate award and an 
NSERC PDF postdoctoral fellowship.

\section{The Yangian $\Yhg$}\label{sec: yangians}
%=====================

\subsection{}
%------------

Let $\g$ be a complex, semisimple Lie algebra and $(\cdot,\cdot)$ an
invariant, symmetric, non--degenerate bilinear form on $\g$. Let $\h
\subset\g$ be a Cartan subalgebra of $\g$, $\{\alpha_i\}_{i\in\bfI}\subset
\h^*$ a basis of simple roots of $\g$ relative to $\h$ and $a_{ij}=2(\alpha
_i,\alpha_j)/(\alpha_i,\alpha_i)$ the entries of the corresponding Cartan
matrix $\bfA$. Let $\Phi_+\subset \h^*$ be the corresponding set of
positive roots, and $\mathsf{Q}=\Z\Phi_+=\bigoplus_{i\in \bfI}\Z\alpha_i\subset
\h^*$ the root lattice. We assume that $(\cdot,\cdot)$ is normalised so
that the square length of short roots is $2$. Set $d_i=(\alpha_i,\alpha_i)
/2\in\{1,2,3\}$, so that $d_ia_{ij}=d_j a_{ji}$ for any $i,j\in\bfI$. 
In addition, we set $h_i=\nu^{-1}(\alpha_i)/d_i$ and choose root vectors 
$x_i^\pm \in \g_{\pm \alpha_i}$ such that $[x_i^+,x_i^-]=d_ih_i$, where
$\nu:\h\to\h^*$ is the isomorphism determined by $(\cdot,\cdot)$.

\subsection{The Yangian $\Yhg$ \cite{drinfeld-yangian-qaffine}}\label{ssec: yangian}
%----------------------------------------

Let $\hbar\in\C$. The Yangian $\Yhg$ is the $\C$--algebra
generated by elements $\{x^{\pm}_{i,r},\xi_{i,r}\}_{i\in\bfI,
r\in\N}$, subject to the following relations.

\begin{enumerate}[font=\upshape]
\item[(Y1)] For any $i,j\in\bfI$, $r,s\in\N$: $[\xi_{i,r}, \xi_{j,s}] = 0$\\

\item[(Y2)] For $i,j\in\bfI$ and $s\in \N$: $
[\xi_{i,0}, x_{j,s}^{\pm}] = \pm d_ia_{ij} x_{j,s}^{\pm}$\\

\item[(Y3)] For $i,j\in\bfI$ and $r,s\in\N$:
\[[\xi_{i,r+1}, x^{\pm}_{j,s}] - [\xi_{i,r},x^{\pm}_{j,s+1}] =
\pm\hbar\frac{d_ia_{ij}}{2}(\xi_{i,r}x^{\pm}_{j,s} + x^{\pm}_{j,s}\xi_{i,r})\]

\item[(Y4)] For $i,j\in\bfI$ and $r,s\in \N$:
\[
[x^{\pm}_{i,r+1}, x^{\pm}_{j,s}] - [x^{\pm}_{i,r},x^{\pm}_{j,s+1}]=
\pm\hbar\frac{d_ia_{ij}}{2}(x^{\pm}_{i,r}x^{\pm}_{j,s} 
+ x^{\pm}_{j,s}x^{\pm}_{i,r})\]

\item[(Y5)] For $i,j\in\bfI$ and $r,s\in \N$:
$[x^+_{i,r}, x^-_{j,s}] = \delta_{ij} \xi_{i,r+s}$\\

\item[(Y6)] Let $i\not= j\in\bfI$ and set $m = 1-a_{ij}$. For any
$r_1,\cdots, r_m\in \N$ and $s\in \N$
\[\sum_{\pi\in\Sym_m}
\left[x^{\pm}_{i,r_{\pi(1)}},\left[x^{\pm}_{i,r_{\pi(2)}},\left[\cdots,
\left[x^{\pm}_{i,r_{\pi(m)}},x^{\pm}_{j,s}\right]\cdots\right]\right]\right.=0\]
\end{enumerate}

We denote by $\Ynaughthg$ and $\Ypmhg$ the unital subalgebras of 
$\Yhg$ generated by $\{\xi_{i,r}\}_{i\in\bfI, r\in\Z_{\geq 0}}$
and $\{x_{i,r}^{\pm}\}_{i\in\bfI, r\in\Z_{\geq 0}}$, respectively.

\subsection{}\label{ssec: formal-series-y}
%--------------

Assume henceforth that $\hbar\neq 0$, and define $\xi_i(u),x^\pm_i(u)
\in \hext{\Yhg}{u^{-1}}$ by
\[\xi_i(u)=1 + \hbar\sum_{r\geq 0} \xi_{i,r}u^{-r-1}\aand
x^{\pm}_i(u)=\hbar\sum_{r\geq 0} x_{i,r}^{\pm} u^{-r-1}\]

\begin{prop}\cite[Prop. 2.3]{sachin-valerio-2}\label{pr:Y fields} 
The relations (Y1),(Y2)--(Y3),(Y4),(Y5) and (Y6) are respectively
equivalent to the following identities in $\Yhg[u,v;u^{-1},v^{-1}]\negthinspace]$.
\begin{enumerate}[font=\upshape]
\item[($\Y$1)] For any $i,j\in\bfI$,
$[\xi_i(u), \xi_j(v)]=0$\\

\item[($\Y$2)] For any $i,j\in\bfI$, 
$[\xi_{i,0},x^\pm_j(u)]=\pm d_ia_{ij}x^\pm_j(u)$\\

\item[($\Y$3)] For any $i,j\in \bfI$, and $a = \hbar d_ia_{ij}/2$:
\[
(u-v\mp a)\xi_i(u)x_j^{\pm}(v)=
(u-v\pm a)x_j^{\pm}(v)\xi_i(u)\mp 2a x_j^{\pm}(u\mp a)\xi_i(u)\]

\item[($\Y$4)] For any $i,j\in \bfI$, and $a = \hbar d_ia_{ij}/2$:
\begin{multline*}
(u-v\mp a) x_i^{\pm}(u)x_j^{\pm}(v)\\
= (u-v\pm a)x_j^{\pm}(v)x_i^{\pm}(u)
+\hbar\lp [x_{i,0}^{\pm},x_j^{\pm}(v)] - [x_i^{\pm}(u),x_{j,0}^{\pm}]\rp
\end{multline*}

\item[($\Y$5)] For any $i,j\in \bfI$:
\[(u-v)[x_i^+(u),x_j^-(v)]=-\delta_{ij}\hbar\left(\xi_i(u)-\xi_i(v)\right)\]

\item[($\Y$6)] For any $i\neq j\in\bfI$, $m=1-a_{ij}$, $r_1,\cdots, r_m\in
\N$, and $s\in \N$:
\[\sum_{\pi\in\Sym_m}
\left[x^{\pm}_i(u_{\pi_1}),\left[x^{\pm}_i(u_{\pi(2)}),\left[\cdots,
\left[x^{\pm}_i(u_{\pi(m)}),x^{\pm}_j(v)\right]\cdots\right]\right]\right.=0\]
\end{enumerate}
\end{prop}

\remark
When $\g=\sl_2$, we will write $\xi_r$, $x_{r}^\pm$, $\xi(u)$ and $x^\pm(u)$ in place of $\xi_{i,r}$, $x_{i,r}^\pm$, $\xi_i(u)$ and $x_i^\pm(u)$, respectively. 

\subsection{Alternative generators of $\Ynaughthg$}\label{ss:ti}
%------------------------------------------------------------------

Let $\{t_{i,r}\}_{i\in\bfI,r\in\N}\subset\Ynaughthg$ be the
generators defined by
\[t_i(u)=\hbar\sum_{r\geq 0}t_{i,r}u^{-r-1} :=\log(\xi_i(u))\]
In particular, $t_{i,0}=\xi_{i,0}$ and
\begin{equation}\label{eq:tione}
%t_{i,0}=\xi_{i,0}
%\aand
t_{i,1}=\xi_{i,1} -\frac{\hbar}{2} \xi_{i,0}^2
\end{equation}
The relations (Y2)--(Y3) of $\Yhg$ imply that
for any $i,j\in\bfI$ and $r\in\N$,
\begin{equation}\label{eq:h1-zero}
[t_{i,1},x_{j,r}^{\pm}] = \pm d_ia_{ij}x_{j,r+1}^{\pm}
\end{equation}
so that $t_{i,1}$ act as shift operators on the generators
$x_{j,r}^\pm$.

The relation \eqref{eq:h1-zero} also implies that
$\{\xi_{i,0},x_{i,0}^\pm, t_{i,1}\}_{i\in \bfI}$ generate $\Yhg$ as
an algebra. We refer the reader to \cite[Thm.~1.2]{levendorskii}
for a presentation of $\Yhg$ given in terms of these generators,
and to \cite[Thm.~2.13]{guay-nakajima-wendlandt} for a refinement
of this result.

\subsection{Shift automorphism}\label{ssec: shift-yangian}
%--------------------------------------

The group of translations of the complex plane acts on
$\Yhg$ by
\[\tau_a(y_r) = \sum_{s=0}^r
\left(\begin{array}{c}r\\s\end{array}\right)
a^{r-s}y_s\]
where $a\in\C$ and $y$ is one of $\xi_i,t_i,x_i^\pm$. In terms of
the generating series introduced in \ref{ssec: formal-series-y}
and \ref{ss:ti}, we have
\[\tau_a(y(u)) = y(u-a)\]
%where $y\in\{\xi_i,x_i^{\pm},t_i\}$.
Given a representation $V$ of $\Yhg$ and $a\in \C$, set
$V(a)=\tau_a^*(V)$.
 
\subsection{PBW theorem}\label{ssec:filt-tau}
%---------------------------------

Consider the loop filtration $\filt{\bullet}{\Yhg}$ on $\Yhg$ defined
by $\deg(y_r) = r$ for each of the generators $y=\xi_i,x^{\pm}_i$.
Note that $\deg(t_{i,r})=r$. The Hopf algebra structure on $\Yhg$
preserves this filtration, and endows $\gr(\Yhg)$ with the structure
of a graded Hopf algebra. The PBW Theorem for $\Yhg$ \cite
{levendorskii-PBW} (see also \cite[Thm. B.6]{FiTsWe19} and
\cite[Prop. 2.2]{GRWEquiv}) is equivalent to the assertion that
the assignments
\[x_i^\pm.z^r\mapsto \bar x_{i,r}^\pm
\aand
d_i h_i.z^r\mapsto \bar \xi_{i,r} 
\]
%$i\in \bfI,\; r\in\Z_{\geq 0}
uniquely extend to an isomorphism of graded Hopf algebras 
\begin{equation}\label{eq:PBW}
U(\g[z])\isom \gr(\Yhg),
\end{equation}
where, for any fixed $k\in \Z_{\geq 0}$ and element $y_k\in \filt{k}{\Yhg}$,
\[
 \bar{y}_k\in \gr_k(\Yhg):=\filt{k}{\Yhg}/\filt{k-1}{\Yhg}
\]
is defined to be the image of $y_k$ in the $k$-th graded component $\gr_k(\Yhg)$ 
of the associated graded algebra $\gr(\Yhg)$.

Henceforth, we shall freely make use of the above identification without further 
comment. Similarly, we will exploit the fact that it allows 
us to identify  $U(\g[z])\otimes U(\g[w])$ with $\gr(\Yhg\otimes\Yhg)\cong \gr(\Yhg)\otimes \gr(\Yhg)$,
the associated graded algebra of $\Yhg\otimes\Yhg$ with respect to the 
tensor product filtration $\filt{\bullet}{\Yhg\otimes \Yhg}$ induced
by $\filt{\bullet}{\Yhg}$.

\subsection{The embedding $U(\g)\subset\Yhg$}\label{ssec:2-embeddings}
%------------------------------------------------------------

Since $\gr_0(\Yhg)=\filt{0}{\Yhg}\subset \Yhg$, the isomorphism \eqref{eq:PBW}
restricts to an embedding of $U(\g)$ into $\Yhg$, given by
\[
x_i^\pm \mapsto x_{i,0}^\pm
\aand
d_ih_i\mapsto \xi_{i,0}
\]
We shall henceforth identify $U(\g)\subset \Yhg$, with the above embedding implicitly understood.
Viewed as a module over $\h\subset \Yhg$, we then have 
$\Yhg = \bigoplus_{\beta\in \mathsf{Q}}\Yhg_{\beta}$

A second embedding $\oT:\h\to\Yhg$ is given by setting $\oT(d_ih_i)
=t_{i,1}$ for all $i\in \bfI$, where $t_{i,1}$ is defined by \eqref{eq:tione}.
The relation \eqref{eq:h1-zero} then reads
\begin{equation}\label{eq:h1-zeroo}
[\Top{h},x_{i,r}^{\pm}] = \pm \alpha_i(h)x_{i,r+1}^{\pm}% \quad \forall \quad h\in \h,  r\in \Z_{\geq 0}. 
\end{equation}
and implies in particular that, for any $h\in\alpha_i^{\perp}$ and $r\geq 0$,
\begin{equation}\label{eq:h1-higher}
[\Top{h},x_{i,r}^{\pm}] = 0
\end{equation}

\subsection{Formal series filtration}
%--------------------------------------------

For any $k\in \Z$, set
\begin{equation}\label{eq:filt-s}
\begin{split}
\filt{k}{\Linfinity{\Yhg^{\otimes 2}}{s})}
&=
s^k\prod_{n\geq 0}\mathcal{F}_n(\Yhg^{\otimes 2})s^{-n}\\
&=
\left\{\left.\sum_{m\leq M} y_ms^m\in
\Linfinity{\Yhg^{\otimes 2}}{s}\right |\deg(y_m)\leq k-m\right\}
\end{split}
\end{equation}
For $k\leq 0$, we shall write this as $\filt{k}{\hext{\Yhg^{\otimes 2}}{s^{-1}})}$,
for obvious reasons.

The above spaces generate a $\Z$-filtered algebra
\[
 \bigcup_{k\in \Z}
 \filt{k}{\Linfinity{\Yhg^{\otimes 2}}{s}}
 \subset \Linfinity{\Yhg^{\otimes 2}}{s}
\]
with associated graded algebra that can (and will) be identified with the $\C[s
^{\pm 1}]$--submodule of  $\Linfinity{(U(\g[z])\otimes U(\g[w])}{s}$ generated
by
\[\prod_{n\geq 0}(U(\g[z])\otimes U(\g[w]))_n s^{-n}\]
where $(U(\g[z])\otimes U(\g[w]))_n$ is the $n$-th graded component of
$U(\g[z])\otimes U(\g[w])$.

If $\lX(s)\in\filt{k}{\Linfinity{\Yhg^{\otimes 2}}{s}}$, we denote by
\[\begin{split}
\ol{\lX(s)}
&=
\lX(s)\mod\filt{k-1}{\Linfinity{\Yhg^{\otimes 2}}{s}}\\
&\in
s^k\prod_{n\geq 0}(U(\g[z])\otimes U(\g[w]))_n s^{-n}\\
&\subset 
\Linfinity{(U(\g[z])\otimes U(\g[w])}{s}
\end{split}\]
the image of $\lX(s)$ in the $k$--th graded component of the associated
graded algebra. Note in passing that if $\lX(s)\in\mathcal{F}_k(\Linfinity
{\Yhg^{\otimes 2}}{s})$, with $k<0$, then $\exp(\lX(s))-1\in\mathcal{F}_
k(\Linfinity{\Yhg^{\otimes 2}}{s})$, and 
\begin{equation}\label{log-deform}
\ol{\exp(\lX(s))-1}=\ol{\lX(s)}.
\end{equation}

\subsection{Rationality}\label{ss:rationality}
%---------------------------

The following rationality property is due to Beck--Kac \cite{beck-kac}
and Hernandez \cite{hernandez-drinfeld} for the analogous case
of the quantum loop algebra, and to the first two authors for $\Yhg$.
In the form below, the result appears in \cite[Prop. 3.6]{sachin-valerio-2}.

\begin{prop}\label{prop: rationality}
Let $V$ be a $\Yhg$--module on which $\{\xi_{i,0}\}_{i\in\bfI}$ act semisimply
with \fd weight spaces. Then, for every weight $\mu$ of $V$, the
generating series
\[\xi_i(u)\in\hext{\End(V_{\mu})}{u^{-1}}
\aand
x_i^{\pm}(u)\in\hext{\Hom(V_{\mu},V_{\mu\pm\alpha_i})}{u^{-1}}\]
are the Laurent expansions at $\infty$ of rational functions of $u$. Specifically,
\[x_i^{\pm}(u)= \hbar  u^{-1}\lp 1 \mp \frac{\ad(t_{i,1})}{2d_iu}\rp^{-1}x_{i,0}^{\pm}\]
and
\[\xi_i(u) = 1 + [x_i^+(u),x_{i,0}^-]\]

\end{prop}

If $V$ is a \fd $\Yhg$--module, we define $\spec(V)\subset \C$ to be the (finite)
set of poles of the rational $\End(V)$--valued functions 
$\{\xi_i(u),x_i^{\pm}(u)\}_{i\in\bfI}$.\\

\section{The standard and Drinfeld coproducts}\label{sec: two-tensors}
%==================================

We review the definition of the standard coproduct on $\Yhg$ following
\cite{guay-nakajima-wendlandt}, and the deformed Drinfeld tensor product
on its \fd representations introduced in \cite{sachin-valerio-III}. We then
lift the latter to a deformed Drinfeld coproduct $\ddelta{s}:\Yhg\to\Linfinity
{\Yhg^{\otimes 2}}{s}$.

\subsection{Standard coproduct}\label{ssec: tensor-ord}
%----------------------------------------

Set
\begin{equation}\label{eq:sfr}
\sfr=\sum_{\beta\in\Phi_+} x^-_{\beta,0}\otimes x^+_{\beta,0}
\end{equation}
where $x^{\pm}_{\beta,0}\in \g_{\pm\beta}\subset\Yhg$ are root vectors
such that $(x^-_{\beta,0},x^+_{\beta,0})=1$. The coproduct $\kmDelta:
\Yhg\to\Yhg\otimes\Yhg$ is defined by the following formulae
\begin{alignat*}{2}
\kmDelta(\xi_{i,0}) 
&= \xi_{i,0}\otimes 1   	&&+ 1\otimes \xi_{i,0}\\[1.1ex]
\kmDelta(x_{i,0}^{\pm}) 
&= x_{i,0}^{\pm}\otimes 1	&&+ 1\otimes x_{i,0}^{\pm}\\[1.1ex]
\kmDelta(t_{i,1})
&= 
t_{i,1}\otimes 1			&&+ 1\otimes t_{i,1} +\hbar\ad(\xi_{i,0}\otimes 1)\sfr \\
&= t_{i,1}\otimes 1		&&+ 1\otimes t_{i,1} - \hbar
\sum_{\beta\in\Phi_+} (\beta,\alpha_i) x^-_{\beta,0}\otimes x^+_{\beta,0}
\end{alignat*}
We refer the reader to \cite[\S 4.2]{guay-nakajima-wendlandt} for a
proof that $\kmDelta$ is an algebra homomorphism. It is immediate
that $\kmDelta$ is coassociative (see \cite[\S 4.5]{guay-nakajima-wendlandt}).

\subsection{Deformed Drinfeld tensor product}\label{ssec: tensor-dr}
%-----------------------------------------------------------

We review below the definition of the deformed Drinfeld tensor product
introduced in \cite[Section 4.4]{sachin-valerio-III}. Let $V,W\in\Ryang$,
and $\spec(V),\spec(W)\subset\C$ their sets of poles. Let $s\in\C$ be
such that $\spec(V)+s$ and $\spec(W)$ are disjoint, and define an
action of the generators of $\Yhg$ on $V\otimes W$ by
\begin{align*}
\ddelta{s}(\xi_i(u)) &= \xi_i(u-s)\otimes\xi_i(u)\\
\ddelta{s}(x_i^+(u)) &= x_i^+(u-s)\otimes 1 + \oint_{C_2} \frac{1}{u-v}
\xi_i(v-s)\otimes x_i^+(v) dv\\
\ddelta{s}(x_i^-(u)) &= \oint_{C_1} \frac{1}{u-v} x_i^-(v-s)\otimes \xi_i(v) dv
+ 1\otimes x_i^-(u)
\end{align*}
where
\begin{itemize}
\item $C_1$ encloses $\spec(V)+s$ and none of the points in $\spec(W)$.
\item $C_2$ encloses $\spec(W)$ and none of the points in $\spec(V)+s$.
\item The integral $\oint_{C_1}$ (resp. $\oint_{C_2}$) %defining $\Delta_{s}(x_i^+(u))$ (resp. $\Delta_{s}(x_i^-(u))$)
is understood to mean the holomorphic function of $u$ it defines for $u$
outside of $C_1$ (resp. $C_2$).
\end{itemize}
Note that in terms of the generators $t_{i,r}$ of $\Ynaughthg$,
\[\ddelta{s}(t_i(u)) = t_{i}(u-s)\otimes 1+1\otimes t_i(u)\]

\begin{thm}\cite[Thm. 4.6]{sachin-valerio-III}\label{th:D Y}\hfill
\begin{enumerate}[font=\upshape]
\item The formulae above define an action of $\Yhg$ on $V\otimes W$.
The corresponding representation is denoted by $V\dtensor{s}W$.
\item The action of $\Yhg$ on $V\dtensor{s} W$ is a rational
function of $s$, with poles contained in $\spec(W)-\spec(V)$.
\item The identification of vector spaces
\[(V_1\dtensor{s_1}V_2)\dtensor{s_2}V_3 = 
V_1\dtensor{s_1+s_2}(V_2\dtensor{s_2}V_3)\]
intertwines the action of $\Yhg$.
\item If $V\cong\C$ is the trivial representation of $\Yhg$, then
\[V\dtensor{s}W=W\aand W\dtensor{s}V=W(s)\]
\item The following holds for any $s,t\in\C$,
\[V(t)\dtensor{s}W(t)=(V\dtensor{s}W)(t)
\aand
V(t)\dtensor{s}W=V\dtensor{s+t}W
\]
In particular, $V\dtensor{s}W(t)=(V\dtensor{s-t}W)(t)$.
%\[V\dtensor{s+s'}W=V(s)\dtensor{s'}W\]
%and $V(s')\dtensor{s}W(s')=(V\dtensor{s}W)(s')$.
%
\item The following holds for any $s\in\C$,
\[\spec(V\dtensor{s}W)\subset (s+\spec(V))\cup\spec(W)\]
\end{enumerate}
\end{thm}

\subsection{Laurent expansion of the deformed Drinfeld tensor product}\label{ss:Laurent}
%--------------------------------------------------------------------------------------------

\begin{prop}\label{pr:laurent drinfeld}
The Laurent expansion at $s=\infty$ of the formulae of \ref{ssec: tensor-dr}
is given by
\[\ddelta{s}(t_{i,r})=
\tau_s(t_{i,r})\otimes 1+
1\otimes t_{i,r}\]
and
\begin{align*}
\ddelta{s}(x_{i,r}^+) &= 
\tau_s(x_{i,r}^+)\otimes 1+1\otimes x_{i,r}^+ \\
&\qquad +
\hbar
\sum_{N\geq 0}\lp \sum_{n=0}^N (-1)^{n+1} \lp \begin{array}{c} N\\ n\end{array}\rp
\xi_{i,n}\otimes x^+_{i,r+N-n}\rp s^{-N-1}\\[1.1ex]
\ddelta{s}(x_{i,r}^-) &= 
\tau_s(x_{i,r}^-)\otimes 1+1\otimes x_{i,r}^-\\
&\qquad +
\hbar
\sum_{N\geq 0}\lp \sum_{n=0}^N (-1)^{n+1} \lp \begin{array}{c} N\\ n\end{array}\rp
x^-_{i,r+n}\otimes\xi_{i,N-n}  \rp s^{-N-1}
\end{align*}
\end{prop}
\begin{pf}
The expansion of $\ddelta{s}(t_{i,r})$ follows from $\ddelta{s}(t_i(u))=\tau_s(t_i(u))
\otimes 1+1\otimes t_i(u)$.
Expanding in $u^{-1}$ yields
\[\ddelta{s}(x_{i,r}^+)
=
\tau_s(x_{i,r}^+)\otimes 1+
\frac{1}{\hbar}
\oint_{C_2}v^r\xi_i(v-s)\otimes x_i^+(v)\, dv\]
Expanding now $\xi_i(v-s)$ \wrt $s^{-1}$ by using $(1-x)^{-p-1}=\sum_{m\geq 0}
\bin{p+m}{p}x^m$ yields:
\[\begin{split}
\xi_i(v-s)
&=
1+\hbar\sum_{p\geq 0}\xi_{i,p}(v-s)^{-p-1}\\
&=
1+\hbar\sum_{p\geq 0}\xi_{i,p}(-s)^{-p-1}\sum_{q\geq 0}\bin{p+q}{p}v^qs^{-q}\\
&=
1+\hbar\sum_{m\geq 0}s^{-m-1}\sum_{\substack{p,q\geq 0\\p+q=m}}(-1)^{p+1}\bin{m}{p}v^q\xi_{i,p}
\end{split}\]
Substituting gives
\[\begin{split}
\ddelta{s}(x_{i,r}^+)
&=
\tau_s(x_{i,r}^+)\otimes 1
+
\frac{1}{\hbar}\oint_{C_2}
v^r 1\otimes x_i^+(v)\, dv\\
&+
\sum_{m\geq 0}
s^{-m-1}\sum_{\substack{p,q\geq 0\\p+q=m}}(-1)^{p+1}\bin{m}{p}
\oint_{C_2}v^{r+q}\xi_{i,p}\otimes x_i^+(v)\, dv
\end{split}\]
which is the claimed result since, for any $a\geq 0$, $\oint_{C_2}v^ax_i^+(v)\,dv=
\hbar x_{i,a}^+$. The expansion of $\ddelta{s}(x_{i,r}^-)$ is obtained in the same way.
\end{pf}

\subsection{Deformed Drinfeld coproduct}\label{ss:def drinf cop}
%-----------------------------------------------------

We now lift the deformed Drinfeld tensor product to an algebra
homomorphism
\[\ddelta{s} : \Yhg\to \Linfinity{\lp\Yhg\otimes\Yhg\rp}{s}.
\]

\begin{thm}\label{th:drdelta}
The Laurent expansions of Section \ref{ss:Laurent} give rise to 
an algebra homomorphism
\[\ddelta{s}: \Yhg\to\Linfinity{\lp\Yhg\otimes\Yhg\rp}{s}\]
The deformed Drinfeld coproduct $\ddelta{s}$ has the following
properties.
\begin{enumerate}[font=\upshape]
% counit
\item It is compatible with the counit $\epsilon$, that is
\[\epsilon\otimes\id\circ\ddelta{s}=\id
\aand
\id\otimes\epsilon\circ\ddelta{s}=\tau_s\]
\item For every $x\in\Yhg$, the following holds in $\Linfinity
{(\Yhg\otimes\Yhg [a])}{s}$
\[ \tau_a\otimes \tau_a\circ\ddelta{s}(x)=\ddelta{s}\circ \tau_a(x)
\aand
\tau_a\otimes\id\circ\ddelta{s}(x) = \ddelta{s+a}(x)
\]
\item $\ddelta{s}$ is a filtered homomorphism, that is
\[
 \ddelta{s}(\mathcal{F}_k(\Yhg))\subset \mathcal{F}_k(\Linfinity{\Yhg^{\otimes 2}}{s}) 
\]
for each $k\geq 0$, where $\mathcal{F}_\bullet(\Linfinity{\Yhg^{\otimes 2}}{s})$
is the filtration defined in \eqref{eq:filt-s}.
\end{enumerate}
\end{thm}
\begin{pf}
Using Theorem \ref{thm:awesome1}, and the fact
that $\kmdelta{s}$ is an algebra homomorphism, we conclude the
same for
$\ddelta{s}$. Properties (1) and (2) above are easy to verify directly
from the definition. Property (3) follows immediately from the explicit
formulas  given in Proposition \ref{pr:laurent drinfeld}. 
\end{pf}

\begin{rem}\label{rm:coassoc-catch}
The coassociativity property of the Drinfeld tensor product $\dtensor{s}$
does not appear to have a natural lift to $\ddelta{s}$. The candidate identity
\begin{equation}\label{eq:candidate}
\ddelta{s_1}\otimes\id
\circ \ddelta{s_2}(x)
=
\id\otimes \ddelta{s_2} 
\circ \ddelta{s_1+s_2}(x)
\end{equation}
holds if $x$ is one of the commuting generators $t_{i,r}$. However, if $x$
is an arbitrary element of $\Yhg$, the \lhs and \rhs lie, respectively, in
\[\Ltwo{\Yhg^{\otimes 3}}{s_1}{s_2}
\aand
\Ltwo{\Yhg^{\otimes 3}}{s_2}{s_1}\]
and cannot therefore be directly compared. The coassociativity of $\dtensor
{s}$ implies, however, that the evaluation of the left-- and right--hand sides
of \eqref{eq:candidate} on a tensor product $V_1\otimes V_2\otimes V_3$
of \fd representations are, respectively, the expansions at $|s_2|\gg |s_1|$
and $|s_1|\gg |s_2|$ of the same rational function.
\end{rem}

%%%%%%%%%%%%%%%%%%%%%%%%%%
%%%%%%%%%%%%%%%%%%%%%%%%%%
\section{The element $\UF(s)$}\label{sec:fusion}%%
%%%%%%%%%%%%%%%%%%%%%%%%%%
%%%%%%%%%%%%%%%%%%%%%%%%%%

\subsection{}\label{ssec:fusion-mainthm}
%--------------

Set $\sfQ_+=\bigoplus_{i\in\bfI}\Z_{\geq 0}\,\alpha_i\subset\h^*$.
The following is one of the main result of this paper.
\begin{thm}\label{thm:awesome1}
\hfill
\begin{enumerate}[font=\upshape]
\item There is a unique element 
\[\UF(s)\in\hext{(\Yhg\otimes\Yhg)}{s^{-1}}\]
which is 
\begin{itemize}
\item strictly lower triangular, that is
$\UF(s)=\sum_{\beta,\gamma\in\sfQ_+}\UF(s)_{\beta,\gamma}$,
with 
\[\UF(s)_{\beta,\gamma}\in
\hext{(\Yhg_{-\beta}\otimes \Yhg_{\gamma})}{s^{-1}}
\quad\text{and}\quad
\UF(s)_{0} = 1\otimes 1\]
\item of weight zero
\item such that, for any $i\in\bfI$,
\begin{equation}\label{eq:A-yangian}
\UF(s)\cdot \kmdelta{s}\lp t_{i,1} \rp = \ddelta{s}\lp 
t_{i,1}\rp \cdot\UF(s)
\end{equation}
\end{itemize}
\noindent
\item The element $\UF(s)$ lies in $\hext{(\Yminushg\otimes\Yplushg)}
{s^{-1}}$, and has the following additional properties.
\begin{itemize}%[resume, font=\upshape]
% x_{i,0}
\item For any $i\in\bfI$,
\begin{equation}\label{eq:X-yangian}
\UF(s)\cdot \kmdelta{s}\lp x_{i,0}^\pm \rp = \ddelta{s}\lp x_{i,0}^\pm \rp \cdot\UF(s)
\end{equation}
% translation
\item For any $a,b\in \C$, 
\[\tau_{a}\otimes \tau_b \left( \UF(s) \right)=\UF(s+a-b)\]
%

% sc limit
\item 
$\UF(s)-1\in \mathcal{F}_{-1}(\hext{\Yhg^{\otimes 2}}{s^{-1}})$, with semiclassical limit given by
\[\overline{\UF(s)-1} = 
\frac{\hbar \fop}{z+s-w}
\in \hext{(U(\g[z])\otimes U(\g[w]))}{s^{-1}}\]
In particular, $\ds \UF(s) = 1+\hbar\fop/s+ O(s^{-2})$.
\end{itemize}

% Rationality and cocycle
\item Let $V_1,V_2\in\Ryang$, and $\UF_{V_1,V_2}(s)\in\hext{\End
(V_1\otimes V_2)}{s^{-1}}$ the corresponding evaluation of $\UF(s)$.
Then,
\begin{itemize}
\item $\UF_{V_1,V_2}(s)$ is the Taylor expansion at $s=\infty$ of a
rational function.
\item The following cocycle equation holds for any $V_1,V_2,V_3
\in\Ryang$:
\begin{equation}\label{eq:V cocycle}
\Fop{V_1\dtensor{s_1} V_2, V_3} (s_2)
\cdot
\Fop{V_1,V_2}(s_1)%\otimes \id_{V_3}
=
\Fop{V_1, V_2\dtensor{s_2} V_3} (s_1+s_2)
\cdot
%\id_{V_1}\otimes
\Fop{V_2,V_3}(s_2)
\end{equation}
\end{itemize}
\end{enumerate}
\end{thm}

\begin{rem}\label{rm:cocycle-catch}
Analogously to Remark \ref{rm:coassoc-catch}, the cocycle equation
\eqref{eq:V cocycle} only holds as an identity of rational functions with values
in $\End(V_1\otimes V_2\otimes V_3)$, and does not seem to possess
a natural lift to $\Yhg$. Indeed, the candidate identity
\[
\ddelta{s_1}\otimes\id(\UF(s_2))
\cdot \UF_{12}(s_1)
=
\id\otimes \ddelta{s_2}\lp \UF(s_1+s_2)\rp
\cdot \UF_{23}(s_2)
\]
does not make sense since the \lhs lies in $\hext{\Linfinity
{\Yhg^{\otimes 3}}{s_1}}{s_2^{-1}}$, while the \rhs lies
in $\hext{\Linfinity
{\Yhg^{\otimes 3}}{s_2}}{s_1^{-1}}$.
\end{rem}

The rest of the section is devoted to the proof of Theorem \ref{thm:awesome1}.

% uniqueness and existence

\subsection{Existence and uniqueness of $\UF(s)$}
\label{ssec:pf-awesome-main}
%-----------------------------------------------------------------

% strategy
The triangularity and zero--weight assumptions are equivalent to
the requirement that $\UF(s)$ have the form
\begin{equation}\label{eq:triangular}
\UF(s)=\sum_{\beta\in \sfQ_+} \UF(s)_{\beta}
\quad\text{with}\quad
\UF(s)_{\beta}\in \hext{(\Yhg_{-\beta}\otimes\Yhg_{\beta})}{s^{-1}}
\end{equation}
and $\UF(s)_0=1\otimes 1$. We shall construct $\UF(s)_\beta$ recursively
in $\beta$, and prove that the sum over $\beta$ converges in the
$s$--adic topology. In fact, define $\nu:\sfQ_+\to \Z_{\geq 0}$ by 
\[\nu(\beta)=\min\left\{k\in\Z_{\geq 0}|
\beta = \alpha^{(1)}+\cdots+\alpha^{(k)},
\alpha^{(1)},\ldots,\alpha^{(k)}\in\Phi_+ \right\}\]
with $\nu(0)=0$. Then, we shall prove that
$\UF(s)_{\beta}\in s^{-\nu(\beta)}\Yhg^{\otimes 2}[\![s^{-1}]\!]$.
%This also implies that $\sum_\beta \UF(s)_\beta$ lies in $\Yhg^{\otimes 2}[\![s^
%{-1}]\!]$ since, for any $m\geq 0$, the set $\{\beta\in\sfQ_+|\,
%\nu(\beta)\leq m\}$ is finite.

% equation spelled out
Let $\sfr$ be given by \eqref{eq:sfr} and, for any $h\in\h$, set
\[\sfr(h)=
\ad(h\otimes 1)(\sfr)=
-\sum_{\gamma\in \Phi_+}
\gamma(h)\, x^-_{\gamma,0}\otimes x^+_{\gamma,0}.
\]
By \ref{ssec: tensor-ord} and \ref{ss:Laurent},
\begin{align*}
\ddelta{s}(\Top{h})
&=
\sum_{a=1}^2 \Top{h}^{(a)}+s h^{(1)}\\
\Delta_{s}(\Top{h})
&=
\sum_{a=1}^2 \Top{h}^{(a)}+s h^{(1)}+\hbar\sfr(h),
\end{align*}
where we use the standard notation: $X^{(1)} = X\otimes 1$ and
$X^{(2)} = 1\otimes X$.
%By \ref{ssec: tensor-ord} and \ref{ss:Laurent}, t
The intertwining equation \eqref{eq:A-yangian} therefore reads as $A(h)
\UF(s)=0$ where, for any $h\in\h$,
\[A(h)=
\lambda(\ddelta{s}(\Top{h}))-\rho\left(\Delta_{s}(\Top{h})\right)=
\ad\left(\Top{h}^{(1)}+\Top{h}^{(2)}+s h^{(1)}\right)
-
\hbar\rho(\sfr(h))\]
and $\lambda,\rho$ denote left and right multiplication, respectively.
In components, this reads:
\[
\lp \adT(h) -s\beta(h)\rp \UF(s)_{\beta}= 
-\hbar \sum_{\alpha\in \Phi_+} \alpha(h) \UF(s)_{\beta-\alpha}\,
x^-_{\alpha,0}\otimes x^+_{\alpha,0}
\]
where $\adT(h)=\ad(\Top{h}^{(1)}+\Top{h}^{(2)})$ and $\UF(s)_{\gamma}=0$
if $\gamma\notin \sfQ_+$.

If $h\in\h$ is such that $\beta(h)\neq 0$ for any nonzero $\beta\in \sfQ_+$, this
yields 
\begin{align}
\UF(s)_{\beta} 
&= \frac{\hbar}{s\beta(h)} \lp 1- \frac{\adT(h)}{s\beta(h)}\rp^{-1}
\sum_{\alpha\in \Phi_+} \alpha(h) \UF(s)_{\beta-\alpha}
x^-_{\alpha,0}\otimes x^+_{\alpha,0}
\label{eq:pf-awesome}\\
&=
\hbar\sum_{k\geq 0} \frac{\mathcal{T}(h)^k}{(s\beta(h))^{k+1}}
\sum_{\alpha\in \Phi_+}\alpha(h)\UF(s)_{\beta-\alpha} x_{\alpha,0}^-
\otimes x_{\alpha,0}^+
\label{eq:pf-awesome2}
\end{align}
This shows that $\UF(s)$ is uniquely determined by $\UF(s)_0$, and
that it lies in $\hext{(\Yminushg\otimes\Yplushg)}{s^{-1}}$ if $\UF(s)_
0$ does, since $\Yminushg,\Yplushg$ are invariant under $\ad \Top{h}$.

Fix now $h\in\h\setminus\bigcup_{\beta\in (\sfQ_+\setminus\{0\})}\Ker\beta$. The above
equations can be used to define elements $\UF(s)_\beta$ recursively
on %$\nu(\beta)$, 
the height of $\beta$,
starting from $\UF(s)_0=1$, which lie in $s^{-\nu(\beta)}\Yhg
^{\otimes 2}[\![s^{-1}]\!]$. The corresponding sum $\UF(s)=\sum_{\beta
\in \sfQ_+}\UF(s)_\beta$ is therefore well--defined, and satisfies $A(h)
\UF(s)=0$. We claim that it also satisfies $A(h')\UF(s)=0$ for any $h'\in\h$.
Note that
\[[A(h),A(h')]
=
\lambda[\ddelta{s}(\Top{h}),\ddelta{s}(\Top{h'})]
-
\rho[\Delta_{s}(\Top{h}),\Delta_{s}(\Top{h'})]
\]
which vanishes because $\Delta_s$ is an algebra homomorphism, and
$\ddelta{s}(\Top{h}),\ddelta{s}(\Top{h'})\in\Ynaughthg^{\otimes 2}$. Thus, $A(h')
\UF(s)$ satisfies
\[A(h)A(h')\UF(s)=A(h')A(h)\UF(s)=0\]
Since $A(h')\UF(s)$ is also triangular, with
\[\left(A(h')\UF(s)\right)_0=
\ad\left(\Top{h'}^{(1)}+\Top{h'}^{(2)}+s{h'}^{(1)}\right)\UF(s)_0=0\]
it follows by uniqueness that $A(h')\UF(s)=0$ as claimed.

\subsection{Translation invariance}
%-------------------------------------------

The identity $\tau_{a}\otimes \tau_b \left( \UF(s) \right)=\UF(s+a-b)$
follows by uniqueness, since both sides are strictly lower triangular,
intertwine
\[\tau_{a}\otimes \tau_b\circ\kmdelta{s}(t_{i,m})=
\tau_{a-b}\otimes\id\circ\kmdelta{s}(\tau_b(t_{i,m}))=
\kmdelta{s+a-b}(\tau_b(t_{i,m}))\]
and $\tau_{a}\otimes \tau_b\circ\ddelta{s}(t_{i,m})=\ddelta{s+a-b}(\tau
_b(t_{i,m}))$ for $i\in\bfI$ and $m=0,1$, and the span of $\{t_{i,m}\}_
{i\in\bfI,m=0,1}$ is invariant under $\tau_b$.

\subsection{Semiclassical limit}\label{ssec:pf-awesome456}
%---------------------------------------

Since $\adT(h)$ is a filtered operator of degree $1$, the recursion
\eqref{eq:pf-awesome2} shows that $\UF(s)_\beta\in\mathcal{F}_
{-\nu(\beta)}(\Yhg^{\otimes 2}[\![s^{-1}]\!])$ for any $\beta\in\sfQ_+$.
In particular, $\UF(s)-1\in\mathcal{F}_{-1}(\Yhg^{\otimes 2}[\![s^{-1}]\!])$
and, mod $\mathcal{F}_{-2}(\Yhg^{\otimes 2}[\![s^{-1}]\!])$,
\[\UF(s)-1
=
\sum_{\beta:\nu(\beta)=1}\UF(s)_\beta\\
=
\frac{\hbar}{s}
\sum_{\alpha\in\Phi_+}
\sum_{k\geq 0} \frac{\mathcal{T}(h)^k}{(s\alpha(h))^k}
x_{\alpha,0}^-\otimes x_{\alpha,0}^+
\]
whose image in $\mathcal{F}_{-1}(\Yhg^{\otimes 2}[\![s^{-1}]\!])/
\mathcal{F}_{-2}(\Yhg^{\otimes 2}[\![s^{-1}]\!])$ is $\hbar\sfr/(z+s-w)$.

\subsection{Rationality}\label{ssec:pf-awesome-rationality}
%----------------------------

The argument given in \ref{ssec:pf-awesome-main} can be carried
out in $\End(V_1\otimes V_2)$ rather than $\Yhg^{\otimes 2}$, and
shows the existence and uniqueness of an element
\[\Fop{V_1,V_2}(s)
=
\sum_{\beta\in \sfQ_+}\Fop{V_1,V_2}(s)_\beta
\in
\hext
{\End(V_1\otimes V_2)}{s^{-1}}\]
with
\[[h\otimes 1,\Fop{V_1,V_2}(s)_\beta]=-\beta(h)\,\Fop{V_1,V_2}(s)_\beta
\qquad
[1\otimes h,\Fop{V_1,V_2}(s)_\beta]=\beta(h)\,\Fop{V_1,V_2}(s)_\beta\]
for any $h\in\h$, $\Fop{V_1,V_2}(s)_0=1$, and 
\[\Fop{V_1,V_2}(s)\cdot \kmdelta{s}\lp t_{i,1}\rp = 
\ddelta{s}\lp t_{i,1}\rp \cdot\Fop{V_1,V_2}(s)\]

The recursive construction of $\Fop{V_1,V_2}(s)$ given by \eqref{eq:pf-awesome}
shows that each $\Fop{V_1,V_2}(s)_\beta$ is a rational function of $s$,
regular at $s=\infty$. Since only finitely many $\Fop{V_1,V_2}(s)_\beta$ 
are non--zero by finite--dimensionality, $\Fop{V_1,V_2}(s)$ is therefore
rational. It follows by uniqueness that the Taylor expansion of $\Fop{V_1,
V_2}(s)$ is the evaluation of $\UF(s)$ on $V_1\otimes V_2$.

\subsection{Cocycle equation}
%=====================

% preamble: a different cocycle identity
Let $V_1,V_2,V_3\in\Ryang$. We obtain below an alternative version
of the cocycle equation, namely 
\begin{equation}\label{eq:K cocycle}
\Fop{V_1,V_2}(s_1)
\cdot
\Fop{V_1\kmtensor{s_1} V_2, V_3} (s_2)
=
\Fop{V_2,V_3}(s_2)
\cdot
\Fop{V_1, V_2\kmtensor{s_2} V_3} (s_1+s_2)
\end{equation}
Note that \eqref{eq:K cocycle} and \eqref{eq:V cocycle} are equivalent,
provided the intertwining equation \eqref{eq:X-yangian} is established.
The latter will be proved in \ref{ss:rank 1 inter}, by relying in part on \eqref
{eq:K cocycle}.\footnote{The cocycle equation \eqref{eq:V cocycle}
arises from the definition of a tensor structure on a functor $F:\cC\to\D$
as a system of natural isomorphisms $J_{U,V}:F(U)\otimes_\D F(V)\to
F(U\otimes_\cC V)$, by interpreting $\UF(s)$ as a tensor structure on
the identity functor $(\Ryang,\dtensor{s})\to(\Ryang,\kmtensor{s})$.
Similarly, \eqref{eq:K cocycle} arises by adopting the opposite
convention of a tensor structure as a system of isomorphisms
$K_{U,V}:F(U\otimes_\cC V)\to F(U)\otimes_\D F(V)$, and interpreting
$\UF(s)$ as a tensor structure on $(\Ryang,\kmtensor{s})\to(\Ryang,
\dtensor{s})$.}

%proof of it now
The intertwining property of $\UF(s)$ implies that
\begin{multline*}
\Fop{V_1,V_2}(s_1)
\cdot
\Fop{V_1\kmtensor{s_1} V_2, V_3}(s_2)
\cdot
\pi_{(V_1\kmtensor{s_1}V_2)\kmtensor{s_2}V_3}(t_{i,1})\\
=
\Fop{V_1,V_2}(s_1)
\cdot
\pi_{(V_1\kmtensor{s_1}V_2)\dtensor{s_2}V_3}(t_{i,1})
\cdot
\Fop{V_1\dtensor{s_1} V_2, V_3} (s_2)\\
=
\pi_{(V_1\dtensor{s_1}V_2)\dtensor{s_2}V_3}(t_{i,1})
\cdot
\Fop{V_1,V_2}(s_1)
\cdot
\Fop{V_1\kmtensor{s_1} V_2, V_3} (s_2),
\end{multline*}
where the second equality stems from the fact that 
\[\pi_{(V_1\kmtensor{s_1}V_2)\dtensor{s_2}V_3}(t_{i,1})=
\pi_{(V_1\kmtensor{s_1}V_2)}(\tau_{s_2}(t_{i,1}))+\pi_{V_3}(t_{i,1})\]
Similarly,
\begin{multline*}
\Fop{V_2,V_3}(s_2)
\cdot
\Fop{V_1, V_2\kmtensor{s_2} V_3} (s_1+s_2)
\cdot
\pi_{V_1\kmtensor{s_1+s_2}(V\kmtensor{s_2}V_3)}(t_{i,1})\\
=
\Fop{V_2,V_3}(s_2)
\cdot
\pi_{V_1\dtensor{s_1+s_2}(V\kmtensor{s_2}V_3)}(t_{i,1})
\cdot
\Fop{V_1, V_2\kmtensor{s_2} V_3} (s_1+s_2)\\
=
\pi_{V_1\dtensor{s_1+s_2}(V\dtensor{s_2}V_3)}(t_{i,1})
\cdot
\Fop{V_2,V_3}(s_2)
\cdot
\Fop{V_1, V_2\kmtensor{s_2} V_3} (s_1+s_2)
\end{multline*}
Since $\kmtensor{s},\dtensor{s}$ are coassociative, both sides of
the cocycle equation \eqref{eq:K cocycle} are therefore solutions
of
\begin{equation}\label{eq:cocycle X}
X \cdot
\Delta_{s_1}\otimes\id\circ\Delta_{s_2}(t_{i,1})\\
=
\ddelta{s_1}\otimes\id\circ\ddelta{s_2}(t_{i,1})
\cdot
X
\end{equation}

By \ref{ssec: tensor-ord} and \ref{ss:Laurent},
\begin{align*}
\ddelta{s_1}\otimes\id\circ\ddelta{s_2}(\Top{h})
&=
\sum_{a=1}^3 \Top{h}^{(a)}+(s_1+s_2)\, h^{(1)}+s_2\,h^{(2)}\\
\Delta_{s_1}\otimes\id\circ\Delta_{s_2}(\Top{h})
&=
\sum_{a=1}^3 \Top{h}^{(a)}+(s_1+s_2)\, h^{(1)}+s_2\,h^{(2)}+
\hbar\left(\sfr(h)_{13}+\sfr(h)_{23}\right)
\end{align*}
The conclusion now follows by noting that, analogously to \ref{ssec:pf-awesome-main},
\eqref{eq:cocycle X} admits at most one rational solution $X(s_1,s_2)$ with values in
$\End(V_1\otimes V_2\otimes V_3)$, provided it is strictly lower triangular, that is of
the form
$X=\sum_{\beta,\gamma\in \sfQ_+}X_{\beta,\gamma}$, where
\[X_{\beta,\gamma}\in
\End(V_1)[-\beta]\otimes \End(V_2)[\beta-\gamma]\otimes
\End(V_3)[\gamma]\]
and $X_{0,0}=\id$.

\subsection{Rank 1 reduction}
%-------------------------------------

Consider the intertwining identity
\begin{equation}\label{eq:x proof}
\UF(s)\cdot \kmdelta{s}(x_{i,0}^\pm)
=
\ddelta{s}(x_{i,0}^\pm)\cdot\UF(s)
\end{equation}
We claim that it holds for any $\g$ and $i\in\bfI$ if, and only if it holds for
$\g=\sl_2$. Set $\sfQ^{(i)}=\sfQ/\Z\alpha_i\subset\h^*/\C\alpha_i=(\alpha
_i^\perp)^*$, and let $\sfQ^{(i)}_+$ be the image of $\sfQ_+$ in $\sfQ^{(i)}$.
Both sides of \eqref{eq:x proof} are lower triangular and of weight zero \wrt
the adjoint action of $\alpha_i^\perp\subset\h$ \ie lie in $\bigoplus_{\beta\in
\sfQ^{(i)}_+}\hext{(\Yhg^{(i)}_{-\beta}\otimes\Yhg^{(i)}_{\beta})}{s^{-1}}$ where,
for $\gamma\in\sfQ^{(i)}$,
\[\Yhg^{(i)}_\gamma=\{x\in\Yhg|\,[h,x]=\gamma(h)x,\,h\in\alpha_i^\perp\}\]
Since both intertwine $\kmdelta{s}(\Top{h})$ and $\ddelta{s}(\Top{h})$ for $h\in
\alpha_i^\perp$, it follows by uniqueness that they are equal if, and only if,
their projections on $\hext{(\Yhg^{(i)}_0\otimes\Yhg^{(i)}_0)}{s^{-1}}$
coincide.

Let $\pi_i\in\End(\Yhg)$ be the projection onto $\Yhg^{(i)}_0=\bigoplus_
{m\in\Z}\Yhg_{m\alpha_i}$. Then,
\begin{align*}
\pi_i\otimes\pi_i(\UF(s)\cdot \kmdelta{s}(x_{i,0}^\pm))
&=
\pi_i\otimes\pi_i(\UF(s))\cdot \kmdelta{s}(x_{i,0}^\pm)\\[1.2ex]
\pi_i\otimes\pi_i(\ddelta{s}(x_{i,0}^\pm)\cdot\UF(s))
&=
\ddelta{s}(x_{i,0}^\pm)\cdot\pi_i\otimes\pi_i(\UF(s))
\end{align*}
Let $\varphi_i:Y_{d_i\hbar}(\mathfrak{sl}_2)\to\Yhg$ be the
algebra homomorphism given by
$x_r^\pm\mapsto d_i^{-1/2}x_{i,r}^\pm$ and
$\xi_r\mapsto d_i^{-1}\xi_{i,r}$. Then,
\[\kmdelta{s}(x_{i,0}^\pm)
=
\sqrt{d_i}\cdot\varphi_i^{\otimes 2}\circ\kmdelta{s}(x_0^\pm)
\aand
\ddelta{s}(x_{i,0}^\pm)
=
\sqrt{d_i}\cdot\varphi_i^{\otimes 2}\circ\ddelta{s}(x_0^\pm).
\]

We claim that $\pi_i\otimes\pi_i(\UF(s))=\varphi_i\otimes\varphi_i
(\UF_{(\sl_2)}(s))$ so that, by the foregoing, \eqref{eq:x proof} holds
for any $\g$ and $i\in\bfI$ if, and only if it holds for $\g=\sl_2$. The
claim follows by uniqueness, since both $\pi_i\otimes\pi_i(\UF(s))$
and $\varphi_i\otimes\varphi_i(\UF_{(\sl_2)}(s))$ are strictly lower
triangular and of weight zero \wrt the action of $\xi_{i,0}$, and both
intertwine
\[\varphi_i^{\otimes 2}\circ\kmdelta{s}(t_i)
=
d_i^{-1}\pi_i^{\otimes 2}\circ\kmdelta{s}(t_{i,1})
\quad\text{and}\quad
\varphi_i^{\otimes 2}\circ\ddelta{s}(t_i)
=d_i^{-1}\,\ddelta{s}(t_{i,1})
=d_i^{-1}\,\pi_i^{\otimes 2}\circ\ddelta{s}(t_{i,1}).
\]

\subsection{Rank 1 intertwining relations}\label{ss:rank 1 inter}
%----------------------------------------------------

Assume $\g=\sl_2$, and consider the identity
\begin{equation}\label{eq:toprove-rank1}
\UF(s)\cdot\kmdelta{s}(x^-_0) = \ddelta{s}(x_0^-)\cdot\UF(s)
\end{equation}
The latter can be proved by a lengthy, direct calculation\footnote
{The calculation is carried out in \S 5 of the earlier version of this
paper, {\sf arXiv:1907.03525 v1}.}. We give below an alternative proof, which relies on the
cocycle identity \eqref{eq:V cocycle} satisfied by $\UF(s)$ to reduce
it to the case when $\UF(s)$ is acting on $\C^2\otimes V$, where
$V$ is an arbitrary \fd representation.

% Separation of points reduction
By Appendix \ref{asec:sep}, it is sufficient to prove that \eqref{eq:toprove-rank1}
holds on the tensor product $V_1\otimes V_2$, where $V_1$ (resp. $V_2$)
is chosen from a collection $\V_1$ (resp. $\V_2$) of \fd representations of
$\Yhg$ which is stable under tensor product, contains the trivial
representation, and a representation whose restriction to $\g$ is faithful.
We choose $\V_2$ to consist of all \fd representations of $\Yhg$, while
$\V_1$ consists of arbitrary tensor products $\C^2(a_1)\otimes\cdots\otimes
\C^2(a_m)$ of evaluation representations.

% Further reduction to C^2(x)V via cocycle identity
For any $V_1,V_2,V_3\in\Ryang$, the cocycle identity \eqref{eq:K cocycle}
implies that \eqref{eq:toprove-rank1} holds on $(V_1(s_1)\otimes V_2)\otimes V_3$
if it holds on $V_1\otimes V_2$, $V_2\otimes V_3$, and $V_1\otimes(V_2(s_2)\otimes
V_3)$. Indeed, 
\[\begin{split}
&\UF_{V_1\kmtensor{s_1}V_2,V_3}(s_2)
\cdot
\pi_{(V_1\kmtensor{s_1}V_2)\kmtensor{s_2}V_3}(x_0^-)\\
%1
&=
\Fop{V_1,V_2}(s_1)^{-1}
\cdot
\Fop{V_2,V_3}(s_2)
\cdot
\Fop{V_1, V_2\kmtensor{s_2} V_3} (s_1+s_2)
\cdot
\pi_{V_1\kmtensor{s_1+s_2}(V_2\kmtensor{s_2}V_3)}(x_0^-)\\
%2
&=
\Fop{V_1,V_2}(s_1)^{-1}
\cdot
\Fop{V_2,V_3}(s_2)
\cdot
\pi_{V_1\dtensor{s_1+s_2}(V_2\kmtensor{s_2}V_3)}(x_0^-)
\cdot
\Fop{V_1, V_2\kmtensor{s_2} V_3} (s_1+s_2)\\
%3
&=
\Fop{V_1,V_2}(s_1)^{-1}
\cdot
\pi_{V_1\dtensor{s_1+s_2}(V_2\dtensor{s_2}V_3)}(x_0^-)
\cdot
\Fop{V_2,V_3}(s_2)
\cdot
\Fop{V_1, V_2\kmtensor{s_2} V_3} (s_1+s_2)
\\
%4
&=
\Fop{V_1,V_2}(s_1)^{-1}
\cdot
\pi_{(V_1\dtensor{s_1}V_2)\dtensor{s_2}V_3}(x_0^-)
\cdot
\Fop{V_2,V_3}(s_2)
\cdot
\Fop{V_1, V_2\kmtensor{s_2} V_3} (s_1+s_2)\\
%5
&=
\pi_{(V_1\kmtensor{s_1}V_2)\dtensor{s_2}V_3}(x_0^-)
\cdot
\Fop{V_1,V_2}(s_1)^{-1}
\cdot
\Fop{V_2,V_3}(s_2)
\cdot
\Fop{V_1, V_2\kmtensor{s_2} V_3} (s_1+s_2)\\
%6
&=
\pi_{(V_1\kmtensor{s_1}V_2)\dtensor{s_2}V_3}(x_0^-)
\cdot
\UF_{V_1\kmtensor{s_1}V_2,V_3}(s_2)
\end{split}\]
It is therefore sufficient to check \eqref{eq:toprove-rank1} on $\C^2(a)
\otimes V$, where $V$ is an arbitrary \fd representation of $\Yhg$, and
$a\in\C$. Moreover, by using the translation invariance of $\UF(s)$,
$\kmtensor{s}$, and $\dtensor{s}$, it is sufficient to consider the case $a=0$.

% Case of C^2(x)V
Let now $x^\pm,\xi$ be the standard nilpotent and semisimple generators of
$\sl_2$ acting on $\C^2$. Then, the following define an action of $\Yhsl
{2}$ on $\C^2$:
\[
x^\pm(u)=\hbar\frac{x^\pm}{u}
\aand
\xi(u)=1+\hbar\frac{\xi}{u}
\]
By \ref{ssec: tensor-ord} and \ref{ssec: tensor-dr}, $\pi_{\C^2\kmtensor{s} V}(x^-_0)
=x^-\otimes 1+1\otimes x^-_0$, and
\[\pi_{\C^2\dtensor{s} V}(x_0^-) = 
\oint_{C_1} \frac{x^-}{v-s}\otimes \xi_i(v)\, dv
+1\otimes x_0^-
=
x^-\otimes \xi(s)+1\otimes x_0^- 
\]
%\oint_{C_2}v^r\xi_i(v-s)\otimes x_i^+(v)\, dv\]
where $C_1$ encloses $\spec(\C^2)+s=\{s\}$ and none of the points in $\spec(V)$.

On the other hand, since $x^-(u)x^-(v)=0$ on $\C^2$, \ref{ss:rat log} and Theorem
\ref{th:rank 1 R} yield
\[\UF_{\C^2,V}(s)
=
1
-
\oint_{C_2}\frac{x^-}{u-s}\otimes x^+(u)\,du
=
1
+
\oint_{C_1}\frac{x^-}{u-s}\otimes x^+(u)\,du=
1
+
x^-\otimes x^+(s)\]
where $C_2$ encloses $\sigma(V)$ and none of the points in $\sigma(\C^2)+s=s$,
and the second equality follows because the residue of the integrand at infinity is $0$.
Therefore,
\begin{align*}
\UF_{\C^2,V}(s)\cdot\pi_{\C^2\kmtensor{s} V}(x^-_0)
&=
x^-\otimes 1+1\otimes x_0^-
+x^-\otimes x^+(s)x_0^- \\
\pi_{\C^2\dtensor{s} V}(x_0^-)\cdot\UF_{\C^2,V}(s)
&=
x^-\otimes \xi(s)+1\otimes x_0^-
+
x^-\otimes x_0^-x^+(s)
\end{align*}
so that
\[\UF_{\C^2,V}(s)\cdot\pi_{\C^2\kmtensor{s} V}(x^-_0)
-\pi_{\C^2\dtensor{s} V}(x_0^-)\cdot\UF_{\C^2,V}(s)
=
x^-\otimes
\left(1
+[x^+(s),x_0^-]
-\xi(s)\right)
\]
which is equal to zero since $[x^+(s),x_0^-]=\xi(s)-1$.

The identity $\UF(s)\cdot\kmdelta{s}(x^+_0) = \ddelta{s}(x_0^+)\cdot
\UF(s)$ is proved in a similar way, by taking $\V_1$ to consist of all
\fd representations, $\V_2$ of tensor products $\C^2(a_1)\otimes\cdots
\otimes\C^2(a_m)$, and using the cocycle identity to reduce this to a
check on $V\otimes\C^2$, for an arbitrary $V\in\Ryang$.

\section{The element $\UF(s)$ for $\sl_2$}\label{sec:fusion-sl2}
%===============================

In this section, we give an explicit formula for the element $\UF(s)$ when
$\g=\sl_2$.

\subsection{}\label{ssec:notation-sl} % setup
%--------------

It will be convenient to consider the generating series
\[\UF(s,z)=
\sum_{n\geq 0}
\UF_{n\alpha}(s)z^n\in\hext{(\Yhsl{2}\otimes\Yhsl{2})[z]}{s^{-1}}
\]
where $\alpha$ is the positive root of $\g$. Let $G=PSL_2(\C)$ be the
complex Lie group of adjoint type corresponding to $\g$, and $H\subset
G$ its maximal torus with Lie algebra $\h$. We identify $H$ with $\C^
\times$ via the character corresponding to $-\alpha$.

In particular, $\UF(s,z)=\Ad(z^{(1)})\UF(s)$. Moreover, $\UF(s,z)$ satisfies
\[\UF(s,z)\cdot \kmdelta{s}^z\lp t_1\rp = \ddelta{s}^z(t_1)\cdot\UF(s,z)\]
where $\kmdelta{s}^z=\Ad(z^{(1)})\circ\kmdelta{s}^z$ and $\ddelta{s}=\Ad
(z^{(1)})\circ\ddelta{s}$, so that
\[\ddelta{s}^z(t_1)
=
\sum_{a=1}^2 t_1^{(a)}+s h^{(1)}
\aand
\kmdelta{s}^z(t_1)
=
\sum_{a=1}^2 t_1^{(a)}+s h^{(1)}+\hbar z\sfr
\]

\subsection{}\label{ssec:rank1-proof}
%---------------

The intertwining equation for $\UF(s,z)$ may be written as the following
ODE, together with the initial condition $\UF(s,0)=1\otimes 1$
\begin{equation}\label{eq:diffeq-def-R}
\lp sz\partial_z - \frac{1}{2} \ad(t_1\otimes 1 + 1\otimes t_1)\rp
\UF(s,z) = \UF(s,z) \cdot \hbar z\fop 
\end{equation}
where $\ds \fop = x_0^-\otimes x_0^+$.

\begin{lem}\label{le:omega}
Set $\logFop(s,z) = \UF(s,z)^{-1}\negthickspace\cdot z\partial_z\UF(s,z)$.
Then, $\logFop(s,z)$ satisfies 
\begin{equation}\label{eq:ODE w}
\lp sz\partial_z - \adT + \hbar z \ad(\fop)\rp\logFop(s,z) 
= \hbar\fop z
\end{equation}
where $\displaystyle \adT = \frac{1}{2} \ad(t_1\otimes 1 + 1\otimes t_1)$.
\end{lem}
\begin{pf}
Since $\lp sz\partial_z - \adT\rp$ is a derivation, we have
\[\begin{split}
\lp sz\partial_z - \adT\rp
\logFop(s,z) &= -\hbar z\fop \logFop(s,z)
+ \UF(s,z)^{-1}z\partial_z(\UF(s,z)\hbar z\fop) \\
&= -\hbar z [\fop,\logFop(s,z)] + \hbar z \fop
\end{split}\]
\end{pf}

\subsection{Formula for $\logFop(s,z)$}\label{ssec:explicitomega}
%--------------------------------------------------

Given a vector space $\HH$, let
\[\ForRes{-}{u}:\Linfinity{\HH}{u}\to\HH\]
be the formal residue, given by taking the coefficient of $u^{-1}$.

\begin{prop}\label{pr:fullomega}
The series $\logFop(s,z)$ is given by
\begin{equation}
\logFop(s,z) = \sum_{k\geq 1} z^k \frac{(-1)^k}{k\hbar}
\ForRes{x^-(u-s)^k\otimes x^+(u)^k}{u}
\label{eq:fullomega}
\end{equation}
\end{prop}
\begin{pf}
Write $\logFop(s,z)= \sum_{k\geq 1}\logFop_k(s)z^k$. For $k=1$,
\eqref{eq:ODE w} yields
\[\begin{split}
\logFop_1(s)
&=
(s-\adT)^{-1}\hbar\sfr
=
\hbar\sum_{m\geq 0}s^{-m-1}\adT^m\sfr\\
&=
\hbar\sum_{m\geq 0}s^{-m-1}\sum_{n=0}^m (-1)^n 
\begin{pmatrix}m\\n\end{pmatrix}x^-_n\otimes x^+_{m-n}
\end{split}\]
where the last equality follows from $[t_1/2,x^\pm_k]=\pm x^\pm_{k+1}$.
On the other hand, 
\[\begin{split}
x^-(u-s)
&=
\hbar\sum_{n\geq 0}x^-_n(u-s)^{-n-1}\\
&=
-\hbar\sum_{n\geq 0}(-1)^n x^-_n s^{-n-1}
\sum_{a\geq 0}\begin{pmatrix}n+a\\n\end{pmatrix}u^as^{-a}\\
&=
-\hbar\sum_{m\geq 0}s^{-m-1}\sum_{n=0}^m(-1)^n
\begin{pmatrix}m\\n\end{pmatrix} x^-_n u^{m-n}
\end{split}\]
which shows that $\logFop_1$ is given by \eqref{eq:fullomega}.

Let now $k\geq 2$, and set $\ds I_k(v,s)=\frac{(-1)^k}{k\hbar} x^-(v-s)^k
\otimes x^+(v)^k$. We claim that 
\[(sk - \adT)I_k(v,s)=-[\hbar\fop, I_{k-1}(v,s)]\]
which proves in particular that \eqref{eq:fullomega} satisfies \eqref{eq:ODE w}.
We shall need the commutation relation $[x^{\pm}_0,x^{\pm}(u)] = \mp x^{\pm}
(u)^2$. The latter follows from relation ($\Y 4$) of Proposition \ref{pr:Y fields},
namely
\[
(u-v\mp \hbar)x^{\pm}(u)x^{\pm}(v) - (u-v\pm \hbar) x^{\pm}(u)x^{\pm}(v)
= \hbar ([x^\pm_0,x^\pm(v)] - [x^{\pm}(u),x^{\pm}_0]),
\]
by taking $u=v$. Thus, for every $k\geq 1$, $\ad(x^{\pm}_0)\cdot x^{\pm}(u)^k
= \mp k x^{\pm}(u)^{k+1}$.

Using this, and $\ds{[\frac{t_1}{2},x^{\pm}(u)]= \pm(ux^{\pm}(u)-\hbar x^{\pm}_0)}$, yields
\[\begin{split}
\left[\frac{t_1}{2},x^{\pm}(u)^k\right]
&=
\pm \sum_{j=1}^k x^{\pm}(u)^{j-1}(ux^{\pm}(u)-\hbar x^{\pm}_0)x^{\pm}(u)^{k-j}\\
&=
\pm k u x^{\pm}(u)^k \mp \hbar\sum_{j=1}^k x^{\pm}(u)^{j-1}x^{\pm}_0x^{\pm}(u)^{k-j}\\
&=
\pm k u x^{\pm}(u)^k \mp \hbar r x^{\pm}(u)^{k-1}x^{\pm}_0+\hbar 
\frac{k(k-1)}{2} x^{\pm}(u)^k
\end{split}\]
Thus,
\begin{multline*}
(sr-\adT)\,x^-(v-s)^k\otimes x^+(v)^k =
-\hbar k(k-1)x^-(v-s)^k\otimes x^+(v)^k\\
-\hbar k x^-(v-s)^{k-1}x^-_0\otimes x^+(v)^k
+\hbar k x^-(v-s)^k \otimes x^+(v)^{k-1}x^+_0.
\end{multline*}

Together with $[A\otimes B, C\otimes D] = [A,C]\otimes [B,D] + [A,C]\otimes DB + CA\otimes [B,D]$,
this yields
\begin{align*}
[x^-_0\otimes x^+_0, x^-(u)^k\otimes x^+(u)^k] &= [x^-_0,x^-(u)^k]\otimes [x^+_0,x^+(u)^k]
+[x^-_0,x^-(u)^k]\otimes x^+(u)^kx^+_0\\
& \hspace*{1cm}+ x^-(u)^kx^-_0\otimes [x^+_0,x^+(u)^k] \\
&= -k^2 x^-(u)^{k+1}\otimes x^+(u)^{k+1} + kx^-(u)^{k+1}\otimes x^+(u)^kx^+_0\\
&\hspace*{1cm}-k x^-(u)^kx^-_0\otimes x^+(u)^{k+1}
\end{align*}
as claimed.
\end{pf}

\subsection{Formula for $\logFop_{V_1,V_2}(s,z)$}\label{ss:rat log}
%-----------------------------------------------------------------

If $V_1,V_2$ are \fd representations of $\Yhg$, Theorem \ref{thm:awesome1}
(3) implies that $\logFop_{V_1,V_2}(s,z)=\pi_{V_1}\otimes\pi_{V_2}(\logFop(s,
z))$ is a rational function of $s$. It follows from the lemma below that
\[\logFop_{V_1,V_2}(s,z) = \sum_{k\geq 1} \frac{(-1)^kz^k}{k\hbar}
\oint_{C_2}x^-(u-s)^k\otimes x^+(u)^k\,du\]
where $C_2$ encloses $\spec(V_2)$ and none of the points in $\spec(V_1)+s$.
Note that in this case the sum over $k$ is finite since $x^\pm(u)$ are nilpotent
on $V_1,V_2$.

\begin{lem}
Let $A$ be a \fd dimensional algebra over $\C$, and $f,g:\C\to A$ rational
functions which are regular at $\infty$. Consider the integral
\[I(s)=\oint_{C_2}f(u-s)g(u)\,du\]
where $s\in\C$, and $C_2$ is a contour enclosing all poles of $g(u)$ and
none of those of $f(u-s)$. Then, the Taylor expansion $\wh{I}(s)$ of $I(s)$
at $s=\infty$ is equal to
\[\wh{I}(s)=\oint\wh{f}(u-s)\wh{g}(u)\,du\]
where $\ds{\oint}-du$ is the formal residue defined in \ref{ssec:explicitomega},
$\wh{f},\wh{g}$ are the Taylor series of $f,g$ at $\infty$, and $\wh{f}(u-s)$
is expanded in $\hext{A[u]}{s^{-1}}$.
\end{lem}
\begin{pf}
$\wh{I}(s)$ is equal to $\ds{\oint_{C_2}\wh{f}(u-s)g(u)du}$. Since $\wh{f}(u-s)
\in\hext{A[u]}{s^{-1}}$, it suffices to prove that $\ds{\oint_{C_2}p(u)g(u)du=
\oint p(u)\wh{g}(u)}du$ for any $p\in A[u]$, which follows by deformation of contour.
\end{pf}

\subsection{Formula for $\UF(s)$}\label{ssec:rank1R-}
%--------------------------------------------

Integrating $z\partial_z\UF(s,z)=\UF(s,z)\cdot\logFop(s,z)$ and setting $z=1$
yields the following corollary of Proposition \ref{pr:fullomega}.

\begin{thm}\label{th:rank 1 R}
The element $\UF(s)\in\hext{\Yhsl{2}^{\otimes 2}}{s^{-1}}$ is given by
\[\UF(s)=
1^{\otimes 2}+\sum_{n\geq 1}
\sum_{\begin{subarray}{c} 
k_1,\ldots,k_r \in\Z_{\geq 1}\\
k_1+\cdots+k_r = n\end{subarray}}
\frac{1}{k_1(k_1+k_2)\cdots(k_1+\cdots+k_r)}
%p(k_1,\ldots,k_r)
\logFop(s)_{k_1}\cdots \logFop(s)_{k_r}
\]
where
\[\logFop(s)_k = \frac{(-1)^k}{k\hbar}
\ForRes{x^-(u-s)^k\otimes x^+(u)^k}{u}\]
\end{thm}

\begin{rem}
It is an interesting problem to give an explicit formula for $\UF(s)$ for $\g
\ncong\sl_2$.
\end{rem}

%%%%%%%%%%%%%%%%%%%%%%%%%%%%%%%%%%%%%%
%%%%%%%%%%%%%%%%%%%%%%%%%%%%%%%%%%%%%%
\section{The universal and the meromorphic abelian $R$--matrices of $\Yhg$}\label{sec:fullRd}
%%%%%%%%%%%%%%%%%%%%%%%%%%%%%%%%%%%%%%
%%%%%%%%%%%%%%%%%%%%%%%%%%%%%%%%%%%%%%

In this section, we review the construction of the meromorphic
abelian $R$--matrix of $\Yhg$ given in \cite{sachin-valerio-III}.
We then show that it gives rise to rational intertwiners $V_1(s)
\dtensor{0}V_2\to V_2\dtensor{0}V_1(s)$ for any \fd representations
$V_1,V_2$. We also prove that these cannot be chosen to be
both natural and compatible with the Drinfeld tensor product.
Finally, we lift the meromorphic abelian $R$--matrix to obtain
a universal abelian $R$--matrix for the deformed Drinfeld coproduct.

\subsection{The endomorphism  $\A_{V_1,V_2}(s)$ \cite{sachin-valerio-III}}
\label{ssec:A Y}
%-------------------------------------------------------------------------------------------------

%\subsubsection{}
%-------------------

Let $V\in\Ryang$. For any $i\in\bfI$, let $\sigma_V(\xi_i)$ be the set
of poles of $\xi_i(u)^{\pm 1}$ acting on $V$, and
\[\mathsf{X}_V(\xi_i) = \bigcup_{a\in\sigma_V(\xi_i)} [0,a]\]
the union of straight line segments joining $0$ to points
in $\sigma_V(\xi_i)$. %, where $[0,a]=\{ta: 0\leq t\leq 1\}$.
Then, the
generating series $t_i(v)=\hbar\sum_{r\geq 0}t_{i,r}v^{-r-1}$ introduced in \ref{ss:ti} converges to a
holomorphic function $t_i(v):\C\setminus\mathsf{X}_V(\xi_i)\to\End
_\C(V)$, which is uniquely determined by $\exp(t_i(v)) = \xi_i(v)$
and $t_i(\infty)=0$ \cite[\S 5.4]{sachin-valerio-III}.

%\subsubsection{}
%-------------------

Let $V_1,V_2\in\Ryang$, $s\in\C$, and define $\A_{V_1,V_2}(s)\in\End_\C
(V_1\otimes V_2)$ by
\[\A_{V_1,V_2}(s)
=
\exp\lp
-\sum_{\begin{subarray}{c} i,j\in\bfI \\ r\in\Z \end{subarray}}
c_{ij}^{(r)}
\oint_{C_1}
%\xi_i(v)^{-1}\cdot\xi_i'(v)
\frac{dt_i(v)}{dv}\otimes t_j\lp
v+s+\frac{(\ell+r)\hbar}{2}
\rp \ dv \rp
\]
where
\begin{itemize}
\item $C_1$ is a contour enclosing $\sigma_{V_1}(\xi_i)$.
\item $\ell=m h^\vee$, with $h^\vee$ the dual Coxeter number of $\g$ and $m=\frac{(\theta,\theta)}{2}$ for $\theta\in \Phi_+$ the highest root.
\item The non--negative integers $c_{ij}^{(r)}$ are the entries
of the following matrix \cite[Appendix A]{sachin-valerio-III}\footnote
{It was proved in \cite[Appendix A]{sachin-valerio-III} that
$c_{ij}(q)\in\N[q,q^{-1}]$. It is clear from the definition that
$c_{ij}(q) = c_{ji}(q) = c_{ij}(q^{-1})$, and the matrix
identity can be expanded as
\begin{equation}\label{eq:cij-s}
\sum_{k\in\bfI} c_{ik}(q)[d_ka_{kj}]_q = \delta_{ij}[\ellg]_q
\ \ \forall i,j\in\bfI
\end{equation}}
\[
\lp c_{ij}(q) = \sum_{r\in\Z} c_{ij}^{(r)} q^r \rp = 
[\ellg]_q \cdot \lp [d_ia_{ij}]_q \rp^{-1}
\]
\item $s$ is large enough so that 
$t_j(v+s+\hbar(\ellg+r)/2)$
is an analytic function of $v$ within $C_1$, for every
$j\in\bfI$ and $r\in \Z$ such that $c_{ij}^{(r)}\neq 0$
for some $i\in\bfI$.
\end{itemize}
Then, $\A_{V_1,V_2}(s)$ is a rational function of $s$, regular
at $\infty$, with an expansion of the form $\ds{\id - \frac{\ellg\hbar}{s^2}\Omega_\h + O(s^{-3})}$,
where $\Omega_\h\in\h\otimes\h$ is the Cartan part of the Casimir
tensor \cite[Thm. 5.5]{sachin-valerio-III}. Moreover, $[\A_
{V_1,V_2}(s),\A_{V_1,V_2}(s')]=0$ for any $s,s'\in\C$.

\subsection{The meromorphic abelian $R$--matrix of $\Yhg$ \cite{sachin-valerio-III}}
\label{ssec:diagonalR}

Consider the additive difference equation determined by $\A_{V_1,V_2}(s)$
\begin{equation}\label{eq:diff A} 
\RR^0_{V_1,V_2}(s+\ellg\hbar) = \A_{V_1,V_2}(s)\cdot \RR^0_{V_1,V_2}(s)
\end{equation}
It admits two meromorphic solutions $\RR^{0,\uparrow/\downarrow}
_{V_1,V_2}(s)$, which are uniquely determined by the requirement
that $\RR^{0,\uparrow}_{V_1,V_2}(s)$ (resp. $\RR^{0,\downarrow}
_{V_1,V_2}(s)$) is holomorphic and invertible for $\Re(s/\ell\hbar)
\gg 0$ (resp. $\Re(s/\ell\hbar)\ll 0$), and possesses an asymptotic
expansion of the form $1+O(s^{-1})$ as $s\to\infty$ in any halfplane
of the form $\Re(s/\ell\hbar)>m$ (resp. $\Re(s/\ell\hbar)<m$). These
solutions are explicitly given by the infinite products
\[
\Rupd_{V_1,V_2}(s)=
\stackrel{\longrightarrow}{\prod_{n\geq 0}}\A_{V_1,V_2}(s+n\ellg\hbar)^{-1}
\aand
\Rdownd_{V_1,V_2}(s)=
\stackrel{\longrightarrow}{\prod_{n\geq 1}}\A_{V_1,V_2}(s-n\ellg\hbar)
\]
The product defining $\Rupd_{V_1,V_2}(s)$ (resp. $\Rdownd_
{V_1,V_2}(s)$) converges uniformly on compact subsets of the complement
of $\calZ-\ell\hbar\Z_{\geq 0}$ (resp. $\calP+\ell\hbar\Z_{>0}$), where
$\calZ$ (resp. $\calP$) is the set of poles of $\A(s)^{-1}$ (resp. $\A(s)$). 

Then, the following holds \cite[Thm 5.9]{sachin-valerio-III}.
\begin{thm}\label{thm:mero-braid-d}
Fix $\varepsilon \in \{\uparrow,\downarrow\}$. Then, the following
holds for any $V_1,V_2,V_3\in \Ryang$
\begin{enumerate}[font=\upshape]
\item %For any $V_1,V_2\in \Ryang$, t
The map 
\[
\flip\circ \URd_{V_1,V_2}(s) : V_1(s)\dtensor{0} V_2
\to V_2\dtensor{0} V_1(s)
\]
%where $\flip$ is the flip of tensor factors,
is a morphism of $\Yhg$--modules, which is natural in $V_1$ and $V_2$.
\item The following cabling identities hold:
\begin{align*}
\URd_{V_1\dtensor{s_1} V_2, V_3} (s_2) &=
\URd_{V_1,V_3}(s_1+s_2)\cdot
\URd_{V_2,V_3}(s_2)\\
\URd_{V_1, V_2\dtensor{s_2} V_3} (s_1+s_2) &=
\URd_{V_1,V_3}(s_1+s_2)\cdot
\URd_{V_1,V_2}(s_1)
\end{align*}

\item For any $a,b\in \C$, 
\[
\URd_{V_1(a),V_2(b)}(s) = \URd_{V_1,V_2}(s+a-b)
\]

\item The following unitary condition holds:
\[
\flip \circ \Rupd_{V_1,V_2}(-s)\circ \flip^{-1} = 
\Rdownd_{V_2,V_1}(s)^{-1}
\]

\item $\RR^{0,\uparrow/\downarrow}_{V_1,V_2}(s)$ have the same asymptotic
expansion, as $s\to\infty$ in any halfplane of the form $\pm\Re(s/\ell\hbar)>m$,
which is of the form
\[\RR^{0,\uparrow/\downarrow}_{V_1,V_2}(s)
\sim
1 + \hbar \Omega_{\h} s^{-1} + O(s^{-2})\]
\end{enumerate}
\end{thm}

\subsection{Existence of a rational intertwiner} 
%-----------------------------------------------------------

Let $V_1,V_2\in\Ryang$.

\begin{thm}\label{th:ratl fact R0}
%For any $V_1,V_2\in\Ryang$, t
There is a rational map $\Rratd_{V_1,V_2}:\C\to\Aut_\C(V_1\otimes
V_2)$, which is normalised by $\Rratd_{V_1,V_2}(\infty)=\id$ and such
that
\[\flip\circ\Rratd_{V_1,V_2}(s):
V_1(s)\dtensor{0}V_2\to V_2\dtensor{0}V_1(s)\]
intertwines the action of $\Yhg$. In particular, $V_1(s)\dtensor{0}V_2$
and $V_2\dtensor{0}V_1(s)$ are isomorphic as $\Yhg$--modules for 
all but finitely many values of $s$.
\end{thm}
\begin{pf}
The map $\Rratd_{V_1,V_2}(s)$ will be obtained from the meromorphic
intertwiners $\RR_{V_1,V_2}^{0,\updown}(s)$ as follows. Consider the
difference equation \eqref{eq:diff A} satisfied by $\RR_{V_1,V_2}^{0,
\updown}(s)$. Its monodromy is given by
\[\eta^0_{V_1,V_2}(s) = \Rupd_{V_1,V_2}(s)^{-1}\cdot\Rdownd_{V_1,V_2}(s)\]

By construction, $\eta^0_{V_1,V_2}$ is an $\ell\hbar$--periodic function
of $s$, and in fact a rational function of $z=\exp\left(\frac{2\pi\iota s}{\ell
\hbar}\right)$  which takes the value $\id$ at $z=0,\infty$ \cite[\S 4.8]
{sachin-valerio-2}. Moreover, by Theorem \ref{thm:mero-braid-d} (1),
$\eta^0_{V_1,V_2}(s)$ commutes with the action of $\Yhg$ on $V_1
\dtensor{s}V_2$. In fact, %for any $X\in\Yhg$,
the periodicity of $\eta^0_{V_1,V_2}$ implies that $\eta^0_{V_1,V_2}
(s)$ commutes with the action of $\Yhg$ on $V_1\dtensor{s'}V_2$ for
any $s'$ of the form $s+\ell\hbar m$, $m\in\Z$. Since that action is
rational in $s'$, it follows that $\eta^0_{V_1,V_2}$ takes values in 
the subalgebra $\calZ\subset\End_\C(V_1\otimes V_2)$ which consists
of elements commuting with the action of $\Yhg$ on $V_1\dtensor
{s'}V_2$ for any $s'\in\C$.

Consider now the following factorisation problem. Find two meromorphic
functions $\lX^{0,\updown}_{V_1,V_2}(s) : \C \to \calZ $
such that
\begin{enumerate}[font=\upshape]
\item $\lX^{0,\updown}_{V_1,V_2}(s)$ is holomorphic and invertible
for $\Re(\pm s/\hbar)\gg 0$.
\item $\lX^{0,\updown}_{V_1,V_2}(s)$ possesses an asymptotic expansion
of the form $1+O(s^{-1})$, valid in any half--plane of the form $\Re
(s/\ell\hbar)\gtrless m$.
\item $\eta^0_{V_1,V_2}(s) = \lX_{V_1,V_2}^{0,\uparrow}(s)^{-1}\cdot 
\lX_{V_1,V_2}^{0,\downarrow}(s)$.
\end{enumerate}
Since $[\eta^0_{V_1,V_2}(s),\eta^0_{V_1,V_2}(s')]=0$ for any $s,s'$, 
such a factorisation exists, and can be obtained explicitly, once a choice
of representatives of poles of $\eta^0(s)$ modulo translations by $\ellg
\hbar\Z$ is made \cite[\S 4.14]{sachin-valerio-2}. We remark that the
consistency equation, required upon such a choice in \cite[\S 4.14]
{sachin-valerio-2}, is vacuous in our case, since $\A_{V_1,
V_2}(s) = 1 + O(s^{-2})$.

Summarising, $\eta^0_{V_1,V_2}(s)$ admits two factorisations in
$\End_\C(V_1\otimes V_2)$, namely
\[
\lX^{0,\uparrow}_{V_1,V_2}(s)^{-1}\cdot 
\lX^{0,\downarrow}_{V_1,V_2}(s)
=
\eta^0_{V_1,V_2}(s)
=
\Rupd_{V_1,V_2}(s)^{-1}\cdot\Rdownd_{V_1,V_2}(s)
\]
Set now
\[\Rratd_{V_1,V_2}(s)=
\Rupd_{V_1,V_2}(s)\cdot\lX^{0,\uparrow}_{V_1,V_2}(s)^{-1}=
\Rdownd_{V_1,V_2}(s)\cdot\lX^{0,\downarrow}_{V_1,V_2}(s)^{-1}
\]
Then, $(1\,2)\circ\Rratd_{V_1,V_2}(s):V_1(s)\dtensor{0}V_2\rightarrow
V_2\dtensor{0}V_1(s)$ intertwines the action of $\Yhg$, and $\Rratd
_{V_1,V_2}$ is a rational function of $s$ equal to $\id$ at $s=\infty$
\cite[\S 4.11]{sachin-valerio-2}. 
\end{pf}

\subsection{Non--existence of rational commutativity constraints}\label{ssec:norational}
%-----------------------------------------------------------------------------------

Theorem \ref{th:ratl fact R0} raises the question of whether the
rational factor $\Rratd_{V_1,V_2}(s)$ may be chosen consistently
for any pair of representations $V_1,V_2$ so as to satisfy the
cabling identities of Theorem \ref{thm:mero-braid-d}. The following
shows that not to be the case.

\begin{thm}\label{th:norational-d}
There is no function $\Rratd_{V_1,V_2}:\C\to\Aut_{\C}(V_1\otimes V_2)$
which is rational, defined for any $V_1,V_2\in\Ryang$, and such that the
following conditions hold.
\begin{enumerate}[font=\upshape]
\item $\flip\circ \Rratd_{V_1,V_2}(s) : V_1(s)\dtensor{0} V_2
\to V_2\dtensor{0} V_1(s)$ is $\Yhg$--linear, and natural
in $V_1$ and $V_2$.
\item For any $V_1,V_2,V_3\in\Ryang$
\begin{align}
\Rratd_{V_1\dtensor{s_1} V_2, V_3} (s_2) &= 
\Rratd_{V_1,V_3}(s_1+s_2)\cdot\Rratd_{V_2,V_3}(s_2)
\label{eq:C01}\\
\Rratd_{V_1, V_2\dtensor{s_2} V_3} (s_1+s_2) &= 
\Rratd_{V_1,V_3}(s_1+s_2)\cdot\Rratd_{V_1,V_2}(s_2)
\label{eq:C02}
\end{align}
\end{enumerate}
\end{thm}
\begin{pf}
% normalisation and translation
Note first if $V_1$ is the trivial one--dimensional representation
$\Triv$ of $\Yhg$, the first cabling identity \eqref{eq:C01} and
part (4) of Theorem \ref{th:D Y} imply that $\Rratd_{\Triv,V_3}(s)
= \Id_{V_3}$. Setting
$V_2=\Triv$ in \eqref{eq:C01} then yields $\Rratd_{V_1(a),V_3}
(s)=\Rratd_{V_1,V_3}(s+a)$.
Similarly, upon setting $V_2=\Triv$ first, and then $V_3=\Triv$,
the second cabling identity \eqref{eq:C02} %and Theorem \ref{th:D Y}
implies that $\Rratd_{V_1,\Triv}(s) = \Id_V$, and that
$\Rratd_{V_1,V_2(b)}(s)=\Rratd_{V_1,V_2}(s-b)$.

% g=sl_2
We now take $\g=\sl_2$ and proceed by contradiction, assuming the
existence of a rational $\Rratd_{V_1,V_2}(s)$ with the stated properties.
We will use the following facts about $\Yhsl{2}$. There is a two--dimensional
representation $\C^2$ of $\Yhsl{2}$ (see Section
\ref{ss:rank 1 inter}. Explicitly, in a fixed basis $\vup{},\vdown{}$,
the action of $\Yhsl{2}$ is given
by the following $2\times 2$ matrices.
\[
\xi(u) = \id+\frac{\hbar}{u}
\left(\begin{array}{rr}
1&0\\0&-1
\end{array}\right)
\aand
x^+(u) = \frac{\hbar}{u}
\begin{pmatrix}0&1\\0&0\end{pmatrix}
= x^-(u)^T
\]
Further, there is a $\Yhsl{2}$--linear map $\Triv \to \C^2\dtensor
{\hbar} \C^2$ given by $1 \mapsto \vup{}\otimes\vdown{}$.

Given $V\in\Ryang$, we consider $\Rratd_{\C^2,V}(s)\in\End_\C
(\C^2\otimes V)$ as a $2\times 2$--matrix with entries in $\End(V)$
given by
\[\Rratd_{\C^2,V}(s)= 
\begin{pmatrix}
\alpha(s) & \beta(s) \\ \gamma(s) & \delta(s)
\end{pmatrix}
\]
The intertwining property of $\Rratd_{\C^2,V}(s)$ reads
\begin{equation}\label{eq:R rat inter}
\Rratd_{\C^2,V}(s)\cdot
\pi_{\C^2(s)\dtensor{0}V}(x)=
(\pi_{V\dtensor{0}\C^2(s)}(x))^{21}
\cdot\Rratd_{\C^2,V}(s)
\end{equation}
for any $x\in\Yhg$.
For $x=\xi_0$,
\[\pi_{\C^2(s)\dtensor{0}V}(\xi_0)=
(\pi_{V\dtensor{0}\C^2(s)}(\xi_0))^{21}=
\xi_0\otimes 1+1\otimes\xi_0=
\begin{pmatrix}
1+\xi_0&0\\0&-1+\xi_0
\end{pmatrix}\]
the intertwining relation \eqref{eq:R rat inter} yields
\[[\xi_0,\alpha]=0\qquad
[\xi_0,\delta]=0\qquad
[\xi_0,\beta]=-2\beta\qquad
[\xi_0,\gamma]=2\gamma\]

% x=t_1
For $x=t_1$, $\pi_{V_1\dtensor{0}V_2}(t_1)=\pi_{V_1}(t_1)+\pi_{V_2}(t_1)$,
so that
\[\pi_{\C^2(s)\dtensor{0}V}(t_1)=
(\pi_{V\dtensor{0}\C^2(s)}(t_1))^{21}=
(t_1+s\xi_0)\otimes 1+1\otimes t_1=
(\xi_1-\frac{\hbar}{2}\xi_0^2+s\xi_0)\otimes 1+1\otimes t_1
\]
Since $\xi_1$ acts by $0$ on $\C^2$, $\pi_{\C^2(s)\dtensor{0}V}(t_1)$ acts
by the matrix
\[\begin{pmatrix}
-\hbar/2+s+t_1&0\\
0&\hbar/2-s+t_1
\end{pmatrix}\]
The relation \eqref{eq:R rat inter} then implies that $[t_1,\alpha]=0$,
$[t_1,\delta]=0$, $[t_1,\beta]=-2s\beta$ and $[t_1,\gamma]=2s\gamma$.
The last two relations imply in turn that $\beta=\gamma=0$ since
$\ad(t_1)\pm s$ is invertible on $\End(V)$ for all but finitely many
values of $s$.

% x=x_0^-
For $x=x_0^-$, \ref{ssec: tensor-dr} implies that the $-+$ coefficient
of $\pi_{V\dtensor{0}\C^2(s)}(x_0^-)$ is $\id_V$, while that of
$\pi_{\C^2(s)\dtensor{0}V}(x_0^-)$ is equal to $\oint_{C_1}\frac{\xi(v)}
{v-s}=\xi(s)$, where $C_1$ encloses $\sigma(\C^2(s))=\{s\}$ and none
of the points in $\sigma(V)$. Taking $-+$ coefficients in \eqref{eq:R rat inter} 
then yields the relation $\alpha(s)=\delta(s)\cdot\xi(s)$. Similarly,
taking the $+-$ coefficients in the intertwining relation \eqref{eq:R rat inter} 
with $x=x^+_0$ yields $\alpha=\xi\cdot\delta$. Summarising, 
\eqref{eq:R rat inter} implies that
\[\Rratd_{\C^2,V}(s)= 
\begin{pmatrix}
\alpha(s) & 0 \\ 0 & \delta(s)
\end{pmatrix}\]
where $\alpha,\delta\in\End(V)(s)$ commute with $\xi_0,t_1$, and
satisfy $\alpha=\delta\cdot\xi=\xi\cdot\delta$.

We now use the cabling relation \eqref{eq:C01}
\[\Rratd_{\C^2_{(1)}\dtensor{\hbar}\C^2_{(2)}, V}(s) = 
\Rratd_{\C^2_{(1)},V}(s+\hbar)\cdot\Rratd_{\C^2_{(2)},V}(s)\]
in $\End_\C(\C^2\otimes\C^2\otimes V)$, where the subscripts are
added to emphasize the order of the tensors. The \rhs applied to
$v_+\otimes v_-\otimes v$ yields $\alpha(s+\hbar)\cdot\delta(s)\,v$.
On the other hand, by naturality, 
\[\Rratd_{\C^2_{(1)}\dtensor{\hbar}\C^2_{(2)}, V}(s)\,v_+\otimes v_-\otimes v
=
\Rratd_{\C, V}(s) 1\otimes v
=
v\]

Combining these relations yields $\alpha(s+\hbar)\cdot\alpha(s)=\xi(s)$.
Taking now $V=\C^2$, this equation implies that the coefficient $a(s)$
of $\vup{}$ in $\alpha(s)\vup{}$ satisfies the additive difference equation
\[\frac{a(s+2\hbar)}{a(s)}=
\frac{c(s+\hbar)}{c(s)}=
\frac{s(s+2\hbar)}{(s+\hbar)^2}\]
where $c(s)=(s+\hbar)/s$ is the matrix coefficient of $\xi(s)$ corresponding
to $v_+$. This equation has a unique solution $\varphi(s)$ which is holomorphic
and non--zero for $\Re(s/2\hbar)\gg 0$ and is asymptotic to $1+O(s^{-1})$ in
that domain (see, \eg \cite[\S 4]{sachin-valerio-2}). Clearly,
\[\varphi(s)=
\frac{\Gamma\left(\frac{s}{2\hbar}\right)\cdot\Gamma\left(\frac{s+2\hbar}{2\hbar}\right)}
{\Gamma\left(\frac{s+\hbar}{2\hbar}\right)^2}
\]
which is not a rational function.
\end{pf}

\subsection{Taylor expansion of $\A_{V_1,V_2}(s)$} % at $\infty$}
\label{ssec:A-formal}
%------------------------------------------------------------------

Let $V_1,V_2\in\Ryang$. We determine in \ref{ss:R0 asy} the asymptotic
expansion of $\URdupdown_{V_1,V_2}(s)$ as $s\to\infty$. As a preliminary
step, we compute the Taylor series of the endomorphism $\A_{V_1,V_2}
(s)$ introduced in \ref{ssec:diagonalR}. Let
\[\LL_{V_1,V_2}(s)
=
-\sum_{\begin{subarray}{c} i,j\in\bfI \\ r\in\Z \end{subarray}}
c_{ij}^{(r)}
\oint_{C_1}
%\xi_i(v)^{-1}\cdot\xi_i'(v)
\frac{dt_i(v)}{dv}\otimes t_j\lp
v+s+\frac{(\ell+r)\hbar}{2}
\rp \ dv\]
be the logarithm of $\A_{V_1,V_2}$.

\begin{lem}
The Taylor series of $\LL_{V_1,V_2}$ at $s=\infty$ is given by the action
on $V_1\otimes V_2$ of the element $\LL(s)\in s^{-2}\hext{(\Ynaughthg
\otimes\Ynaughthg)}{s^{-1}}$ defined by
\[\LL(s)=-\hbar^2\sum_{\begin{subarray}{c} i,j\in\bfI \\ r\in\Z \end{subarray}}
c_{ij}^{(r)}\,
T_{\frac{\ell+r}{2}}\,
\sum_{m\geq n\geq 0}
(-1)^n(m+1)! s^{-m-2}\,
\frac{t_{i,n}}{n!}\otimes\frac{t_{j,m-n}}{(m-n)!}\]
where $T_m f(s)=f(s+m\hbar)$.
\end{lem}
\begin{pf}
Given that
\[\begin{split}
t_j(v+s)
&=
\hbar\sum_{m\geq 0}t_{j,m}(v+s)^{-m-1}\\
&=
\hbar\sum_{m\geq 0}
s^{-m-1}t_{j,m}\sum_{n\geq 0}(-1)^n\begin{pmatrix}m+n\\n\end{pmatrix}
v^ns^{-n}\\
&=
\hbar\sum_{m\geq 0}
s^{-m-1}\sum_{n=0}^m(-1)^n\begin{pmatrix}m\\n\end{pmatrix}
t_{j,n-m}v^n
\end{split}\]
and that $t_i(v)'=-\hbar\sum_{a\geq 0}(a+1)t_{i,a}v^{-a-2}$, the expansion
of the summand of $\LL_{V_1,V_2}$ corresponding to a triple $i,j,r$ is equal
to
\[\hbar^2 c_{ij}^{(r)}T_{\frac{\ell+r}{2}}
\sum_{m\geq 0}
s^{-m-1}\sum_{n=1}^m(-1)^n\begin{pmatrix}m\\n\end{pmatrix}
n t_{i,n-1}\otimes t_{j,n-m}\]
\end{pf}

\begin{rem}\label{rk:formal L}
Since $(m+1)!s^{-m-2}=(-1)^m\partial_s^m s^{-2}$, the above reads
\begin{equation}\label{eq:L Borel}
\LL(s)
=
\sum_{\begin{subarray}{c} i,j\in\bfI \\ r\in\Z \end{subarray}}
c_{ij}^{(r)}\,
T_{\frac{\ell+r}{2}}\,
B_i(\partial_s)\otimes B_j(-\partial_s)\, (-s^{-2})
\end{equation}
where $B_k(z)=\hbar\sum_{n\geq 0}\frac{t_{k,n}}{n!}z^n$ is the inverse Borel
transform of $t_k(v)$.
\end{rem}

\subsection{Asymptotic expansion of $\RR^{0,\uparrow/\downarrow}_{V_1,V_2}(s)$}
%-------------------------------------------------------------------------------------------------------------
\label{ss:R0 asy}

Let\footnote{$g(x)$ is equal to $\sum_{k\geq 0}B_k x^{-k-1}$, and is the asymptotic
expansion at $x=\infty$ of the trigamma function $\Psi_1(x)=\frac{d^2}{dx^2}\ln\Gamma
(x)=\sum_{n\geq 0}(x+n)^{-2}$.} $g(x)\in x^{-1}\hext{\C}{x^{-1}}$ be the unique solution
of the difference equation $g(x+1)=g(x)-x^{-2}$.

Define $\RR^0(s)\in 1+s^{-1}\hext{(\Ynaughthg\otimes\Ynaughthg)}{s^{-1}}$ by\\
\begin{equation}\label{eq:log R0}
\log(\RR^0(s))
= 
\frac{1}{\ellg^2\hbar^2}
\sum_{\begin{subarray}{c} i,j\in\bfI \\ r\in\Z \end{subarray}}
c_{ij}^{(r)} T_{\frac{\ellg+r}{2}}\,
B_i(\partial_s)\otimes B_j(-\partial_s)\,
 g\lp\frac{s}{\ellg\hbar}\rp
 \end{equation}
Let $V_1,V_2\in\Ryang$. By Theorem \ref{thm:mero-braid-d}, $\URdupdown
_{V_1,V_2}(s)$ has an asymptotic expansion as $s\to\infty$ in any halfplane
$\Re(s/\ell\hbar)\gtrless m$, which is of the form $1+O(s^{-1})$.

\begin{prop}\label{pr:asym-Rd}
The asymptotic expansion of $\URdupdown_{V_1,V_2}(s)$ as $s\to\infty$
is given by
\[
\RR^{0}_{V_1,V_2}(s)=\pi_{V_1}\!\otimes\pi_{V_2}(\RR^0(s))
\]
\end{prop}
\begin{pf}
By definition of $\log(\RR^0(s))$ and \eqref{eq:L Borel}, we have
\[\begin{split}
(T_{\ellg}-1)\log(\RR^0(s)) 
&= 
\frac{1}{\ellg^2\hbar^2}
\sum_{\begin{subarray}{c} i,j\in\bfI \\ r\in\Z \end{subarray}}
c_{ij}^{(r)} T_{\frac{\ellg+r}{2}}\,
B_i(\partial_s)\otimes B_j(-\partial_s)\,
\lp g\lp\frac{s}{\ellg\hbar}+1\rp-g\lp\frac{s}{\ellg\hbar}\rp\rp\\
&=
\sum_{\begin{subarray}{c} i,j\in\bfI \\ r\in\Z \end{subarray}}
c_{ij}^{(r)} T_{\frac{\ellg+r}{2}}\,
B_i(\partial_s)\otimes B_j(-\partial_s)\,
\lp -\frac{1}{s^2} \rp\\
&=
\LL(s)
\end{split}\]
Thus, $\RR^0_{V_1,V_2}(s)$ is the unique formal solution of
\[
 \RR^0_{V_1,V_2}(s+\ellg\hbar)
=\A_{V_1,V_2}(s)\cdot \RR^0_{V_1,V_2}(s)
\]
and therefore equals the asymptotic expansion of $\URdupdown
_{V_1,V_2}(s)$.
\end{pf}

\subsection{Properties of $\RR^0(s)$}
\label{ssec:R0-formal-thm}
%------------------------------------------------

The following is the universal analogue of Theorem \ref{thm:mero-braid-d}.

\begin{thm}\label{thm:R0-formal}
% $\RR^0(s)\in 1+s^{-1}\hext{(\Ynaughthg\otimes\Ynaughthg)}{s^{-1}}$
%The element
$\RR^0(s)\in\hext{(\Ynaughthg\otimes\Ynaughthg)}{s^{-1}}$ has the following properties.
\begin{enumerate}[font=\upshape]
\item For every $x\in \Yhg$, the following holds in $\Linfinity
{\Yhg^{\otimes 2}}{s}$
\[\RR^0(s)\cdot\ddelta{s}(x) =
\ddelta{-s}^{\scriptstyle{(21)}}
(\tau_s(x))\cdot \RR^0(s)\]
\item The cabling identities
\begin{align*}
\ddelta{s_1}\otimes\id\lp\RR^0(s_2)\rp&=
\RR^0_{13}(s_1+s_2)\cdot\RR^0_{23}(s_2)\\
\id\otimes \ddelta{s_2}\lp\RR^0(s_1+s_2)\rp&=
\RR^0_{13}(s_1+s_2)\cdot\RR^0_{12}(s_1)
\end{align*}
hold in $\hext{\Ynaughthg^{\otimes 3}[s_1]}{s_2^{-1}}$
and $\hext{\Ynaughthg^{\otimes 3}[s_2]}{s_1^{-1}}$
respectively.

\item $\RR^0(s)$ is unitary
\[\RR^0(s)^{-1}=
\RR^0_{21}(-s)\]

\item For any $a,b\in \C$
\[
(\tau_a\otimes\tau_b)(\RR^0(s)) = \RR^0(s+a-b)
\]
\item $\RR^0(s)-1\in \mathcal{F}_{-1}(\hext{\Yhg^{\otimes 2}}{s^{-1}})$, with semiclassical limit
given by
\[\overline{\RR^0(s)-1} = 
\frac{\hbar \Omega_{\h}}{z+s-w}
\in \hext{(U(\g[z])\otimes U(\g[w]))}{s^{-1}}
\]
\end{enumerate}
\end{thm}

Statements (1)--(4) follow from that fact that $\RR^0_{V_1,V_2}(s)$ is the
asymptotic expansion of $\URdupdown_{V_1,V_2}(s)$ and Theorem \ref
{thm:mero-braid-d}, since \fd representations separate points in $\Yhg$
(Proposition \ref{P:sep}). For completeness, we give a direct proof of
Theorem \ref{thm:R0-formal} below, which does not rely on this fact.

\subsection{Direct proof of Theorem \ref{thm:R0-formal}}
\label{ssec:R0-formal-proof}
%------------------------------------------------------------------------
\hfill\break

%2
\noindent {\bf Proof of (2).} Let $P^0\subset\Ynaughthg$ be
the $\C$--linear span of $\{t_{i,r}\}_{i\in\bfI, r\in\N}$, so that $
\log(\RR^0(s))\in\hext{(P^0\otimes P^0)}{s^{-1}}$. The cabling
identities follow from the fact that each $t_{i,r}$ is primitive \wrt
the Drinfeld coproduct, that is satisfies
\[\ddelta{s}(t_{i,r}) = \tau_s(t_{i,r})\otimes 1 + 1\otimes t_{i,r}\]
and the fact that $\Ynaughthg$ is a commutative subalgebra.\\

%3
\noindent {\bf Proof of (3).}  Write
\[\begin{split}
\log(\RR^0(s))
&= 
\frac{1}{\ellg^2\hbar^2}
\sum_{\begin{subarray}{c} i,j\in\bfI \\ r\in\Z \end{subarray}}
c_{ij}^{(r)} T_{\frac{r}{2}}\,%e^{\frac{r\hbar}{2}\partial_s}\,%
B_i(\partial_s)\otimes B_j(-\partial_s)\,
g\lp\frac{s}{\ellg\hbar}+\frac{1}{2}\rp\\
&=
\frac{1}{\ellg^2\hbar^2}
\sum_{i,j\in\bfI}
c_{ij}(e^{\frac{\hbar}{2}\partial_s})\,%T_{\frac{r}{2}}\,
B_i(\partial_s)\otimes B_j(-\partial_s)\,
g\lp\frac{s}{\ellg\hbar}+\frac{1}{2}\rp\\
&=
\frac{1}{\ellg^2\hbar^2}\,
\Omega(\partial_s)\,
g\lp\frac{s}{\ellg\hbar}+\frac{1}{2}\rp
\end{split}\]
where $\Omega(z)=\sum_{i,j\in\bfI}c_{ij}(z)B_i(z)\otimes B_j(-z)$.
The unitary condition follows from
\[\Omega^{21}(z) = \Omega(-z) \aand 
g\lp \frac{1}{2} + x\rp = - g\lp\frac{1}{2} -x\rp
\]
The first identity holds because $c_{ij}(q) = c_{ij}(q^{-1}) = c_{ji}(q)$,
and the second because $g(x)=-g(1-x)$, since both sides are solutions
of the same difference equation.\\

%4
\noindent {\bf Proof of (4).} Since $\tau_a B_i(z) = e^{az}B_i(z)$,
\[\begin{split}
\tau_a\otimes\tau_b\,\log(\RR^0(s))
&= 
\frac{1}{\ellg^2\hbar^2}
\sum_{\begin{subarray}{c} i,j\in\bfI \\ r\in\Z \end{subarray}}
c_{ij}^{(r)} T_{\frac{\ellg+r}{2}}\,
B_i(\partial_s)\otimes B_j(-\partial_s)e^{(a-b)\partial_s}\,
g\lp\frac{s}{\ellg\hbar}\rp\\
&=
\frac{1}{\ellg^2\hbar^2}
\sum_{\begin{subarray}{c} i,j\in\bfI \\ r\in\Z \end{subarray}}
c_{ij}^{(r)} T_{\frac{\ellg+r}{2}}\,
B_i(\partial_s)\otimes B_j(-\partial_s)\,
g\lp\frac{s+a-b}{\ellg\hbar}\rp\\
&=
\log(\RR^0(s+a-b))
\end{split}\]

%5
\noindent {\bf Proof of (5).} By \eqref{log-deform}, it suffices to prove
that $\log(\RR^0(s))\in\filt{-1}{\hext{\Yhg^{\otimes 2}}{s^{-1}}}$, and
that
\begin{equation*}
\overline{\log(\RR^0(s))}=
\hbar\frac{\Omega_{\h}}{z+s-w}
\end{equation*}
Note first that
\[t_{i,n}\otimes t_{j,m}\,(-1)^m\partial^m_s g\lp\frac{s}{\ell\hbar}\rp=
t_{i,n}\otimes t_{j,m}\,\lp\frac{\ell\hbar\cdot m!}{s^{m+1}}+O(s^{-m-2})\rp\]
lies in $\filt{-1}{\hext{\Yhg^{\otimes 2}}{s^{-1}}}$, and has symbol
$\ell\hbar\cdot m!\,d_ih_iz^n\otimes d_jh_jw^m\,s^{-m-1}$. Since
the shifts $T_x$ preserve the filtration $\filt{\bullet}{\hext{\Yhg
^{\otimes 2}}{s^{-1}}}$ and act as the identity on its associated graded
space, it follows
that $\log(\RR^0(s))\in\filt{-1}{\hext{\Yhg^{\otimes 2}}{s^{-1}}}$, and
that
\[\begin{split}
\overline{\log(\RR^0(s))}
&=
\frac{\hbar}{\ellg}
\sum_{i,j\in\bfI}
c_{ij}(1)\,
\sum_{m\geq n\geq 0}
(-1)^n
\begin{pmatrix}m\\n\end{pmatrix}
d_ih_iz^n\otimes d_jh_jw^{m-n}\,s^{-m-1}
\\
&=
\hbar\frac{\Omega_{\h}}{z+s-w}
\end{split}\]
where the last equality follows from the fact that
\[\sum_{j\in \bfI} c_{ij}(1)d_jh_j=\ellg\sum_{j\in \bfI} 
(\mathbf{B}^{-1})_{ij} d_j h_j=\ellg \varpi_i^{\vee}
\]
with $\mathbf{B}=(d_i a_{ij})$, and $\varpi_i^{\vee}\in \h$ the fundamental
coweights.\\

%1
\noindent {\bf Proof of (1).}
The intertwining relation is obvious for $x\in\Ynaughthg$, since
the latter is the commutative algebra generated by the elements
$t_{i,r}$, which satisfy
\[\ddelta{s}(t_{i,r}) = 
\tau_s(t_{i,r})\otimes 1 + 1\otimes t_{i,r}=
\ddelta{-s}^{(21)}(\tau_s t_{i,r})\]

Thus, it suffices to prove that, for any $k\in\bfI$
\[
\RR^0(s) \ddelta{s}(x^{\pm}_{k,0}) = 
\ddelta{-s}^{(21)}(x^{\pm}_{k,0}) \RR^0(s)
\]
We verify this identity for the $+$ case only. By Proposition
\ref{pr:laurent drinfeld}, $\ddelta{s}(x^+_{k,0})$ is equal to
$\primitive{x^+_{k,0}}+\lX_k(s)$, where
\[\lX_k(s)=
\hbar\sum_{N\geq 0} s^{-N-1} 
\sum_{n=0}^N (-1)^{n+1}\cbin{N}{n}
\xi_{k,n}\otimes x^+_{k,N-n}\]
We therefore have to prove that
\[\Ad(\RR^0(s))\cdot \lp\primitive{x^+_{k,0}}\rp
=
\primitive{x^+_{k,0}} + \lX^{(21)}_k(-s)
-\Ad(\RR^0(s))\cdot \lX_k(s)\]

We claim that $\ds \Ad(\RR^0(s))\cdot (x^+_{k,0}\otimes 1)
= x^+_{k,0}\otimes 1 + \lX^{(21)}_k(-s)$. Given this, the unitary
condition (3) then implies that $\Ad(\RR^0(s))^{-1}\cdot (1
\otimes x^+_{k,0})= 1\otimes x^+_{k,0} + \lX_k(s)$ which,
combined with the claim yields the required intertwining
equation for $x^+_{k,0}$.

To prove the claim, we rely on the following commutation
relation, which was obtained in \cite[\S 2.9]{sachin-valerio-1}
\[
[B_i(z),x^+_{k,n}] = \frac{e^{\frac{d_ia_{ik}\hbar}{2}z}
-e^{-\frac{d_ia_{ik}\hbar}{2}z}}{z} \cdot
\lp\sum_{p\geq 0} x^+_{k,n+p} \frac{z^p}{p!}\rp
\]
Combining with the definition of $\Omega(z)$ given above,
we can carry out the following computation,
for each $k\in\bfI$, $n\in\N$, and $y\in\Ynaughthg$.
\begin{align*}
[\Omega(z),x^+_{k,n}\otimes y] &=
\frac{1}{\ellg^2\hbar^2}\sum_{j\in\bfI}
\lp
\sum_{i\in\bfI} c_{ij}\lp e^{\frac{\hbar}{2}z}\rp
\lp 
\frac{e^{\frac{d_ia_{ik}\hbar}{2}z}
-e^{-\frac{d_ia_{ik}\hbar}{2}z}}{z}
\rp
\rp \\
& \ \ \ \cdot
\lp\sum_{p\geq 0} x^+_{k,n+p} \frac{(-z)^p}{p!}\rp \otimes
B_j(z)y \\
&= \frac{1}{\ellg^2\hbar^2}
\frac{e^{\frac{\ellg\hbar}{2}z}-e^{-\frac{\ellg\hbar}{2}z}}{z}\cdot
\lp
\sum_{p\geq 0} x^+_{k,n+p} \frac{(-z)^p}{p!}
\rp \otimes B_k(z)y
\end{align*}
Note that we used the equation \eqref{eq:cij-s} satisfied
by $(c_{ij}(q))$ above.

This calculation, combined with 
\[
\frac{e^{\frac{\ellg\hbar}{2}\partial_s}
-e^{-\frac{\ellg\hbar}{2}\partial_s}}{\partial_s} \cdot
g\lp \frac{s}{\ellg\hbar} + \frac{1}{2}\rp
 = \frac{\ellg^2\hbar^2}{s}
\]
yields the commutation relation
\begin{align*}
[\log(\RR^0(s)),x^+_{k,n}\otimes y] &=
\sum_{p\geq 0} (-1)^p x^+_{k,n+p}\otimes
\lp
\hbar \sum_{r\geq 0} \cbin{r+p}{r}
t_{k,r} s^{-r-p-1}
\rp y
\end{align*}
The claim now follows from
\[\Ad(\RR^0(s)) = \exp(\ad(\log(\RR^0(s))))\]
where both sides are acting on $V_k^+\otimes \Ynaughthg$, where $V_k^+$
is the $\C$--linear span of $\{x_{k,n}\}_{n\geq 0}$.
\newpage

%%%%%%%%%%%%%%%%%%%%%%%%%%%%%%%%%%%%%%%%%
%%%%%%%%%%%%%%%%%%%%%%%%%%%%%%%%%%%%%%%%%
\section{The universal and the meromorphic $R$--matrices of $\Yhg$}\label{sec:fullR}
%%%%%%%%%%%%%%%%%%%%%%%%%%%%%%%%%%%%%%%%%
%%%%%%%%%%%%%%%%%%%%%%%%%%%%%%%%%%%%%%%%%

In this section, we construct the meromorphic and universal
$R$--matrices of $\Yhg$.

\subsection{The meromorphic $R$--matrix} \label{ss:fullR}
%------------------------------------------------------

Given $V_1,V_2\in\Ryang$ and $\veps\in\{\uparrow,\downarrow\}$,
define $\UR_{V_1,V_2}:\C\to\End(V_1\otimes V_2)$ by
\[
\UR_{V_1,V_2}(s) =
\RR^+_{V_1,V_2}(s)\cdot
\URd_{V_1,V_2}(s)\cdot
\UF_{V_1,V_2}(s),
\]
where $\RR_{V_1,V_2}^+(s)=\flip\circ\RR_{V_2,V_1}^-(-s)^{-1}\circ
\flip$.

\begin{thm}\label{thm:mero-braid}
The meromorphic function $\UR_{V_1,V_2}(s)$ has the following
properties.
\begin{enumerate}[font=\upshape]
% flip
\item The map 
\[
\flip\circ \UR_{V_1,V_2}(s) : V_1(s)\kmtensor{} V_2
\to V_2\kmtensor{} V_1(s)\]
is a morphism of $\Yhg$--modules, which is natural in $V_1,V_2$.
% cabling
\item For any $V_1,V_2,V_3\in\Ryang$,
\begin{align*}
\UR_{V_1\kmtensor{s_1} V_2, V_3} (s_2) &=
\UR_{V_1,V_3}(s_1+s_2)
\cdot
\UR_{V_2,V_3}(s_2)\\
\UR_{V_1, V_2\kmtensor{s_2} V_3} (s_1+s_2) &=
\UR_{V_1,V_3}(s_1+s_2)
\cdot
\UR_{V_1,V_2}(s_1).
\end{align*}
% QYBE
In particular, the QYBE holds on $V_1\otimes V_2\otimes V_3$:
\[
\UR_{V_1,V_2}(s_1)\UR_{V_1,V_3}(s_1+s_2)\UR_{V_2,V_3}(s_2)
= 
\UR_{V_2,V_3}(s_2)\UR_{V_1,V_3}(s_1+s_2)\UR_{V_1,V_2}(s_1).
\]
% translation
\item For any $a,b\in \C$, 
\[
\UR_{V_1(a),V_2(b)}(s) = \UR_{V_1,V_2}(s+a-b).
\]
% unitarity
\item $\Rup_{V_1,V_2}(s)$ and $\Rdown_{V_2,V_1}(s)$ are related by the
unitarity relation:
\[
\flip \circ \Rup_{V_1,V_2}(-s)\circ \flip=
\Rdown_{V_2,V_1}(s)^{-1}.
\]
% asymptotic expansion
\item $\RR^{\updown}_{V_1,V_2}(s)$ have the same asymptotic expansion,
which is of the form
\[
\RR^{\updown}_{V_1,V_2}(s)
\sim 
1 + \hbar\Omega_{\g} s^{-1} + O(s^{-2})
%1 + \hbar (
%\pi_1\otimes\pi_2)(\Omega_{\g}) s^{-1} + O(s^{-2})
%= (\pi_{V_1}\otimes\pi_{V_2})(\RR(s)),
\]
as $s\to\infty$ in any halfplane of the form $\Re(s/\hbar)\gtrless m$.
\end{enumerate}
\end{thm}

\begin{pf}
(1) By definition, 
\[
\flip\circ \UR_{V_1,V_2}(s) = 
\UF_{V_2,V_1}(-s)^{-1}\cdot\lp \flip\circ \URd_{V_1,V_2}(s) \rp\cdot
\UF_{V_1,V_2}(s).
\]
The result therefore follows from the fact that $\UF_{V_1,V_2}(s)$
is a morphism of $\Yhg$--modules $V_1\kmtensor{s} V_2 \to V_1
\dtensor{s} V_2$ (Theorem \ref{thm:awesome1}), and Theorem
\ref{thm:mero-braid-d} (1).\\

(2) We will prove the following equivalent version of the first
cabling identity
\begin{equation}\label{eq:toprove-cable1}
(1\ 2\ 3)\circ \UR_{V_1\kmtensor{s_1}V_2,V_3}(s_2)
=
(1\ 2) \lp \UR_{V_1,V_3}(s_1+s_2)\otimes \id \rp
(2\ 3) \lp \id\otimes \UR_{V_2,V_3}(s_2) \rp
\end{equation}

By definition, the \lhs is equal to
\begin{align*}
%(1\ 2\ 3)\circ \UR_{V_1\kmtensor{s_1}V_2,V_3}(s_2)
%\text{L.H.S.}
& \UF_{V_3,V_1\kmtensor{s_1}V_2}(-s_2)^{-1}\cdot
\lp (1\ 2\ 3)\circ \URd_{V_1\kmtensor{s_1} V_2,V_3}(s_2)\rp\cdot
\UF_{V_1\kmtensor{s_1}V_2,V_3}(s_2)\\
=& \UF_{V_3,V_1\kmtensor{s_1}V_2}(-s_2)^{-1}\cdot
\lp \id\otimes \UF_{V_1,V_2}(s_1)^{-1} \rp \cdot
\lp (1\ 2\ 3) \circ \URd_{V_1\dtensor{s_1} V_2,V_3}(s_2) \rp\cdot \\
& \qquad \lp \UF_{V_1,V_2}(s_1)\otimes \id \rp\cdot
\UF_{V_1\kmtensor{s_1}V_2,V_3}(s_2) \\
=& \UF_{V_3,V_1\kmtensor{s_1}V_2}(-s_2)^{-1}\cdot
\lp \id\otimes \UF_{V_1,V_2}(s_1)^{-1} \rp \cdot
 (1\ 2) \lp \URd_{V_1,V_3}(s_1+s_2)\otimes \id \rp \cdot \\
& \qquad (2\ 3) \lp \id\otimes \URd_{V_2,V_3}(s_2)\rp\cdot 
\lp \UF_{V_1,V_2}(s_1)\otimes \id \rp\cdot
\UF_{V_1\kmtensor{s_1}V_2,V_3}(s_2)
\end{align*}
In the first equality, we used Theorem \ref{thm:awesome1}
in order to change $V_1\kmtensor{s_1}V_2$ to $V_1\dtensor{s_1} V_2$,
while the second equality follows from the cabling identity satisfied by
$\URd(s)$ (Theorem \ref{thm:mero-braid-d} (2)).

Note that we have the following identity, which follows from the
cocycle equation \eqref{eq:K cocycle} after renaming variables
\begin{gather*}
\lp \UF_{V_1,V_3}(s_1+s_2)\otimes \id \rp \cdot
\UF_{V_1\kmtensor{s_1+s_2} V_3,V_2}(-s_2) = 
\lp \id\otimes \UF_{V_3,V_2}(-s_2)\rp \cdot
\UF_{V_1,V_3\kmtensor{-s_2} V_3}(s_1)
\end{gather*}

Inserting this operator and its inverse in the last line
of the computation above allows us to write the \lhs of
\eqref{eq:toprove-cable1} as $A(s_1,s_2)\cdot B(s_1,s_2)$, where
\begin{align*}
A(s_1,s_2) &= \UF_{V_3,V_1\kmtensor{s_1}V_2}(-s_2)^{-1}\cdot
\lp \id\otimes \UF_{V_1,V_2}(s_1)^{-1} \rp \cdot
 (1\ 2) \lp \URd_{V_1,V_3}(s_1+s_2)\otimes \id \rp \cdot \\
& \qquad \lp \UF_{V_1,V_3}(s_1+s_2)\otimes \id \rp \cdot
\UF_{V_1\kmtensor{s_1+s_2} V_3,V_2}(-s_2)\\
B(s_1,s_2) &= \lp\lp \id\otimes \UF_{V_3,V_2}(-s_2)\rp \cdot
\UF_{V_1,V_3\kmtensor{-s_2} V_3}(s_1)\rp^{-1}\cdot \\
& \qquad (2\ 3) \lp \id\otimes \URd_{V_2,V_3}(s_2)\rp\cdot 
\lp \UF_{V_1,V_2}(s_1)\otimes \id \rp\cdot
\UF_{V_1\kmtensor{s_1}V_2,V_3}(s_2)
\end{align*}
Thus, in order to prove \eqref{eq:toprove-cable1}, it is 
enough to show that 
\[
A(s_1,s_2) = (1\ 2) \circ \UR_{V_1,V_3}(s_1+s_2)\otimes \id \aand
B(s_1,s_2) = (2\ 3) \circ \id\otimes \UR_{V_2,V_3}(s_2)
\]
We verify the latter below, the proof of the former, being
entirely analogous, is omitted.\\

Using the cocycle equation \eqref{eq:K cocycle} again, we
have:
\begin{align*}
B(s_1,s_2) &= \UF_{V_1,V_3\kmtensor{-s_2} V_2}(s_1)^{-1}
\cdot \lp \id\otimes \UF_{V_3,V_2}(-s_2)^{-1}\rp \cdot
\lp (2\ 3) \circ \id \otimes \URd_{V_2,V_3}(s_2) \rp \cdot \\
& \qquad \lp \id\otimes \UF_{V_2,V_3}(s_2) \rp \cdot
\UF_{V_1,V_2\kmtensor{s_2}V_3}(s_1+s_2)\\
&= \UF_{V_1,V_3\kmtensor{-s_2} V_2}(s_1)^{-1}\cdot
\lp (2\ 3) \circ \id \otimes \UR_{V_2,V_3}(s_2) \rp \cdot
\UF_{V_1,V_2\kmtensor{s_2}V_3}(s_1+s_2)\\
&= \UF_{V_1,V_3\kmtensor{-s_2} V_2}(s_1)^{-1}\cdot
\UF_{V_1,V_3\kmtensor{-s_2} V_2}(s_1)\cdot
\lp (2\ 3) \circ \id \otimes \UR_{V_2,V_3}(s_2) \rp\\
&= (2\ 3) \circ \id \otimes \UR_{V_2,V_3}(s_2)
\end{align*}
In the second line, we have used the definition of $\UR(s)$,
and in the third, Part (1) of this theorem.\\

(3) follows from Theorem \ref{thm:awesome1} (2), and Theorem
\ref{thm:mero-braid-d} (3).\\

(4) By definition,
\begin{align*}
\flip\circ\Rup_{V_1,V_2}(-s)\circ\flip &=
\UF_{V_2,V_1}(s)^{-1}\cdot \lp
\flip\circ\Rupd_{V_1,V_2}(-s)\circ\flip\rp
\cdot \RR^+_{V_2,V_1}(s)^{-1}\\
&= \UF_{V_2,V_1}(s)^{-1}\cdot \lp
\Rdownd_{V_2,V_1}(s)^{-1}\rp
\cdot \RR^+_{V_2,V_1}(s)^{-1}\\
&= \Rdown_{V_2,V_1}(s)^{-1}
\end{align*}
where the second equality uses Theorem \ref{thm:mero-braid-d} (4).\\

(5) follows from Theorem \ref{thm:mero-braid-d} (5) and
the Taylor series expansion of $\UF_{V_1,V_2}(s) = 1+\hbar \fop/s+
O(s^{-2})$ given in Theorem \ref{thm:awesome1} (3).
\end{pf}

\subsection{Existence of a rational intertwiner} %factorisation}
%----------------------------------------------------------

The following extends to an arbitrary pair of representations $V_1,V_2
\in\Ryang$ a result due to Drinfeld, which is valid when $V_1,V_2$ are
irreducible \cite[Thm. 4]{drinfeld-qybe}, and Maulik--Okounkov, which 
is valid when $\g$ is simply--laced, and $V_1,V_2$ arise from geometry
\cite{maulik-okounkov-qgqc}.

\begin{thm}\label{th:ratl fact R}
There is a rational map $\Rrat_{V_1,V_2}(s):\C\to\End_\C
(V_1\otimes V_2)$, which is normalised by $\Rrat_{V_1,V_2}(\infty)=
\id$ and such that
\[\flip\circ \Rrat_{V_1,V_2}(s):
V_1(s)\kmtensor{}V_2\to V_2\kmtensor{}V_1(s)\]
intertwines the action of $\Yhg$. In particular, $V_1(s)\kmtensor{}V_2$
and $V_2\kmtensor{}V_1(s)$ are isomorphic as $\Yhg$--modules for 
all but finitely many values of $s$.
\end{thm}
\begin{pf}
Let $\Rratd_{V_1,V_2}(s)$ be the rational operator such that $(1\,2)\circ
\Rratd_{V_1,V_2}(s)$ intertwines the action on $V_1(s)\dtensor{0}V_2$
and $V_2\dtensor{0}V_1(s)$, as obtained in Theorem \ref{th:ratl fact R0}.
Then, 
\[\Rrat_{V_1,V_2}(s)
=
\RR^+_{V_1,V_2}(s)\cdot
\Rratd_{V_1,V_2}(s) \cdot \UF_{V_1,V_2}(s)\]
yields the required map.
\end{pf}

\subsection{Non existence of rational commutativity constraints}
%----------------------------------------------------------------------------------

\begin{thm}\label{th:norational}
There is no function $\Rrat_{V_1,V_2}:\C\to\Aut_{\C}(V_1\otimes V_2)$
which is rational, defined for any $V_1,V_2\in\Ryang$, and such that the
following holds
\begin{enumerate}[font=\upshape]
\item $\flip\circ \Rrat_{V_1,V_2}(s) : V_1(s)\kmtensor{} V_2\to V_2
\kmtensor{} V_1(s)$ intertwines the action of $\Yhg$, and is natural
in $V_1$ and $V_2$.
\item For any $V_1,V_2,V_3\in\Ryang$, %
\begin{align*}
\Rrat_{V_1\kmtensor{s_1} V_2, V_3} (s_2) &= 
\Rrat_{V_1,V_3}(s_1+s_2)\cdot \Rrat_{V_2,V_3}(s_2)\\
\Rrat_{V_1, V_2\kmtensor{s_2} V_3} (s_1+s_2) &= 
\Rrat_{V_1,V_3}(s_1+s_2)\cdot \Rrat_{V_1,V_2}(s_2)
\end{align*}
\end{enumerate}
\end{thm}

\begin{pf}
Assume that such a rational $\Rrat_{V_1,V_2}(s)$ exists. 
Set
\[
\Rratd_{V_1,V_2}(s) = \RR^+_{V_1,V_2}(s)^{-1}
\cdot \Rrat_{V_1,V_2}(s) \cdot \UF_{V_1,V_2}(s)^{-1}
\]
Then, the rational operator $\Rratd_{V_1,V_2}(s)$ contradicts
Theorem \ref{th:norational-d}, and the claim follows.
\end{pf}

\subsection{The universal $R$-matrix}\label{ssec:Uni-R}
%------------------------------------------------

We now turn our attention to obtaining a formal, universal analogue
of Theorem \ref{thm:mero-braid}. Consider the formal power series
\begin{equation*}
\RR(s)=\RR^+(s)\cdot\RR^0(s)\cdot\RR^-(s)\in \hext{(\Yhg\otimes\Yhg)}{s^{-1}}
\end{equation*}
where $\RR^+(s)=\RR_{21}^-(-s)^{-1}$. This series admits an expansion
\begin{equation*}
\RR(s)=1+\sum_{k=1}^\infty \RR_k s^{-k}, \qquad \RR_k\in \mathcal{F}_{k-1}(\Yhg\otimes\Yhg).  
\end{equation*}
\begin{thm}\label{thm:R-formal}
The formal power series $\RR(s)$ has the following properties.
 \begin{enumerate}[font=\upshape]
  \item\label{UniR:1} For every $x\in \Yhg$, the following holds in $\Linfinity{\Yhg
^{\otimes 2}}{s}$
\[
\tau_s\otimes\id\circ \Delta\sop(x)
= 
\RRD(s)\cdot
\tau_s\otimes\id\circ \Delta(x)
\cdot
\RRD(s)^{-1}
\]
  \item\label{UniR:2} $\RR(s)$ satisfies the cabling identities
  \begin{align*}
   \Delta\otimes \id(\RR(s))&=\RR_{13}(s)\RR_{23}(s)\\
   \id\otimes \Delta(\RR(s))&=\RR_{13}(s)\RR_{12}(s)
  \end{align*}
  \item\label{UniR:3}  $\RR(s)$ is unitary
  \begin{equation*}
   \RR(s)^{-1}=\RR_{21}(-s). 
  \end{equation*}
  \item\label{UniR:4} For any $a,b\in \C$, we have
  \begin{equation*}
   (\tau_a\otimes \tau_b)\RR(s)=\RR(s+a-b)
  \end{equation*}  
  \item\label{UniR:5} $\RR(s)-1\in \mathcal{F}_{-1}(\hext{\Yhg^{\otimes 2}}{s^{-1}})$, with semiclassical limit
given by
  \begin{equation*}
   \overline{\RR(s)-1}=
\frac{\hbar \Omega_{\g}}{s+z-w} \in \hext{(U(\g[z])\otimes U(\g[w]))}{s^{-1}}
  \end{equation*}
  In particular, $\RR(s)=1+\hbar s^{-1}\Omega_{\g}+O(s^{-2})$.
  
  \item\label{UniR:6} For any  $V_1,V_2\in \Ryang$ and $\veps\in\{\uparrow,\downarrow\}$,
 we have
  \begin{equation*}
   \UR_{V_1,V_2}(s) \sim \RR_{V_1,V_2}(s)   
\end{equation*}
as $s\to\infty$ in any halfplane of the form $\Re(s/\hbar)\gtrless m$.
Here $\RR_{V_1,V_2}(s)=\pi_{V_1}\!\otimes\pi_{V_2}(\RR(s))$. That is, 
$\RR_{V_1,V_2}(s)$ is equal to the asymptotic expansion of 
$\UR_{V_1,V_2}(s)$ from (5) of Theorem \ref{thm:mero-braid}.
 \end{enumerate}
\end{thm}
\begin{pf}
 Parts (1) and (3)--(6) are deduced directly from the definition of 
$\RR(s)$ using the properties of $\RR^-(s)$ and $\RR^0(s)$ established in Theorems 
\ref{thm:awesome1} and \ref{thm:R0-formal}, respectively. 

\medskip 

We note, however, that the cabling identities (2) do not follow 
directly from Theorems \ref{thm:awesome1} and \ref{thm:R0-formal}, 
as we do not have access to a formal, universal version of the 
cocycle equation \eqref{eq:V cocycle} of Theorem \ref{thm:awesome1}. 
To remedy this, we shall instead make use of the fact that 
$\Ryang$ is sufficiently large to distinguish elements of $\Yhg$. More precisely, by Proposition \ref{P:sep} (with $n=3$), it is enough to prove that the following identities hold on $V_1\otimes V_2\otimes V_3$, for every $V_1,V_2,V_3\in \Ryang$: 
  \begin{equation}\label{R-cablingV}
  \begin{aligned}
  \RR_{V_1\otimes V_2, V_3}(s)&=\RR_{V_1,V_3}(s)\RR_{V_2,V_3}(s)\\
   \RR_{V_1,V_2\otimes V_3}(s)&=\RR_{V_1,V_3}(s)\RR_{V_1,V_2}(s)
  \end{aligned}
  \end{equation}
Fix $\veps\in\{\uparrow,\downarrow\}$ and consider the first equality. By (6), the left-hand side (resp. right-hand side) is equal to the uniquely determined asymptotic expansion of $\UR_{V_1\otimes V_2, V_3}(s)$ (resp. $\UR_{V_1,V_3}(s)\UR_{V_2,V_3}(s)$) as $s\to \infty$ in 
any halfplane of the form $\Re(s/\hbar)\gtrless m$. Moreover, by the first equality in (2) of Theorem \ref{thm:mero-braid}, taken with $s_1=0$ and $s_2=s$, we have 
\begin{equation*}
 \UR_{V_1\otimes V_2, V_3}(s)=\UR_{V_1,V_3}(s)\UR_{V_2,V_3}(s).
\end{equation*}
Thus, we can conclude that the first cabling identity of \eqref{R-cablingV} necessarily holds. An identical argument establishes the second identity. 
\end{pf}

As an immediate consequence of the above theorem and the uniqueness
assertion of Theorem \ref{thm:intro1} (see also Appendix \ref{asec:UR-unique}),
we obtain the following corollary.
\begin{cor}
$\RR(s)$ is equal to Drinfeld's universal $R$--matrix. %$\RR^{(D)}(s)$.
\end{cor}
In particular, Theorem \ref{ssec:Uni-R} provides an independent, and
constructive proof of the existence of Drinfeld's universal $R$--matrix.

\section{Meromorphic tensor structures}\label{sec:merocats}
%=============================

In this section, we reinterpret our results in the language
of meromorphic tensor categories. We refer to \cite
{soibelman-mero, soibelman-meromorphic} for a more
abstract and general treatment of meromorphic tensor
categories. We caution the reader, however, that the
framework developed in \cite{soibelman-mero,soibelman-meromorphic}
does not include examples where the tensor product 
depends non--trivially on a parameter, as is the case
for the deformed Drinfeld tensor product. 
The setup of \cite{soibelman-mero,soibelman-meromorphic}
is also more general than needed for our purposes, in
that it is adapted to {\em pseudo--tensor categories},
where the tensor product need not be defined for all
pairs of representations, or be representable.

\subsection{Drinfeld tensor product}\label{ssec:reps-merocat-d}
%--------------

\begin{prop}\label{pr:merocat-d}
\hfill
\begin{enumerate}[font=\upshape]
\item The category $(\Ryang,\dtensor{s})$ is a meromorphic
(in fact, rational) tensor category over $(\C,+)$.
\item Each of the resummed abelian $R$--matrices $\URdupdown(s)$
is a meromorphic braiding on $(\Ryang,\dtensor{s})$.
\item $(\Ryang,\dtensor{s})$ does not admit a rational braiding.
\end{enumerate}
\end{prop}

\begin{pf}
(1) $\Ryang$ admits an action of the additive group $(\C,+)$
given by $V\mapsto V(s)$. As recalled in \ref{ssec: tensor-dr}, for
every $V,W\in \Ryang$, there is a rational action of $\Yhg$ on $V\otimes
W$ given by the deformed Drinfeld tensor product. The properties (1)--(5)
of Theorem \ref{th:D Y} mean exactly that $(\Ryang, \dtensor{s})$
is a rational tensor category over $(\C,+)$.

(2) is the content of  Theorem \ref{thm:mero-braid-d} (1)--(3).

(3) is a rephrasing of Theorem \ref{th:norational-d}.
\end{pf}

\subsection{Standard tensor product}\label{ssec:reps-merocat}
%----------------------------------------------

\begin{prop}\label{pr:merocat}
\hfill
\begin{enumerate}[font=\upshape]
\item The category $(\Ryang,\kmtensor{s})$ is a meromorphic
(in fact, polynomial) tensor category over $(\C,+)$.
\item Each of the resummed $R$--matrices $\URupdown(s)$
is a meromorphic braiding on $(\Ryang, \kmtensor{s})$.
\item $(\Ryang,\kmtensor{s})$ does not admit a rational braiding.
\end{enumerate}
\end{prop}
\begin{pf}
\noindent (1) This is a consequence of the fact that
$V_1\kmtensor{s} V_2$ arises from the algebra homomorphism
$\kmdelta{s} : \Yhg \to (\Yhg\otimes\Yhg)[s]$. 
This tensor product satisfies the properties analogous
to (3)--(5) of Theorem \ref{th:D Y}.

(2) is the content of (1)--(3) of Theorem \ref{thm:mero-braid}.

(3) is a rephrasing of Theorem \ref{th:norational}.
\end{pf}

\subsection{Meromorphic tensor structures}\label{ss:mero tens str}
%-------------------------------------------------------

\begin{prop}\label{pr:mero tens str}
$\UF(s)$ is a rational braided tensor structure on the identity functor
$$\lp \Ryang, \dtensor{s}, \URdupdown(s) \rp
\longrightarrow
\lp \Ryang, \kmtensor{s}, \URupdown(s) \rp$$
\end{prop}
\begin{pf}
By definition of a tensor structure
on a functor, the statement means that, for every
$V_1,V_2\in\Ryang$, there is a rational $\End(V_1\otimes V_2)$--valued
function of $s$, $\UF_{V_1,V_2}(s)$, which satisfies
(1)--(3) of Theorem \ref{thm:awesome1}. Namely,
\[\UF_{V_1,V_2}(s) : V_1\kmtensor{s} V_2
\to V_1\dtensor{s} V_2\]
is a $\Yhg$--intertwiner such that $\UF_{V_1(a),V_2(b)}(s)
=\UF_{V_1,V_2}(s+a-b)$, and the following diagram commutes,
for every $V_1,V_2,V_3\in\Ryang$
% associativity constraints
\[\xymatrix@R=1.4cm@C=.6cm{
(V_1\otimes_{s_1} V_2) \otimes_{s_2} V_3
\ar[d]_{\RR^-_{V_1,V_2}(s_1)\otimes\id_{V_3}}
\ar@{=}[rr]&&
V_1\otimes_{s_1+s_2} ( V_2 \otimes_{s_2} V_3)
\ar[d]^{\id_{V_1}\otimes \RR^-_{V_2,V_3}(s_1)}\\
( V_1\dtensor{s_1} V_2) \otimes_{s_2} V_3
\ar[d]_{\RR^-_{V_1\dtensor{s_1}V_2,V_3}(s_2)}&&
V_1\otimes_{s_1+s_2} ( V_2 \dtensor{s_2} V_3)
\ar[d]^{\RR^-_{V_1,V_2\dtensor{s_2}V_3}(s_1+s_2)}\\
( V_1\dtensor{s_1} V_2) \dtensor{s_2} V_3
\ar@{=}[rr]&&
V_1\dtensor{s_1+s_2} ( V_2 \dtensor{s_2} V_3)
}\]

% commutativity constraints
Lastly, it is claimed in (3) that $\UF(s)$ is compatible with the
braidings on the two categories, that is satisfies
\[\xymatrix@R=1.6cm{
V_1(s)\otimes V_2 
\ar[rr]^{\flip\circ\URupdown_{V_1,V_2}(s)}
\ar[d]_{\RR^-_{V_1,V_2}(s)}
&&
V_2\otimes V_1(s)
\ar[d]^{\RR^-_{V_2,V_1}(-s)}\\
V_1(s)\dtensor{0} V_2
\ar[rr]_{\flip\circ\URdupdown_{V_1,V_2}(s)}
&&
V_2\dtensor{0} V_1(s)
}\]
The commutativity of the diagram follows from the fact that, by
definition
\[\RR^{\uparrow/\downarrow}_{V_1,V_2}(s) =
\RR^+_{V_1,V_2}(s)\cdot
\URdupdown_{V_1,V_2}(s)\cdot
\UF_{V_1,V_2}(s)\]
where $\RR^+_{V_1,V_2}(s) = \flip \circ \UF_{V_2,V_1}(-s)^{-1}
\circ \flip$.
\end{pf}

%%%%%%%%%%%%%%%%%%%%%%%%%%%%%%%%%%%
%%%%%%%%%%%%%%%%%%%%%%%%%%%%%%%%%%%
\section{Relation to the quantum loop algebra $\qloop$}
\label{se:qloop}
%%%%%%%%%%%%%%%%%%%%%%%%%%%%%%%%%%%
%%%%%%%%%%%%%%%%%%%%%%%%%%%%%%%%%%%

In this section, we review the construction of the tensor functor between
\fd representations of $\Yhg$ and the quantum loop algebra $\qloop$
obtained in \cite{sachin-valerio-2,sachin-valerio-III}. We then disprove
a conjecture stated in \cite{sachin-valerio-III}, and relate the meromorphic
$R$--matrices of $\Yhg$ and $\qloop$.

\subsection{The functor $\sfGamma$  \cite{sachin-valerio-2}}
\label{ssec:mero-tensor-D}
%-------------------------------------------------------------------------------------------------------------------------------------------------------

% qloop
Set $q=\exp(\pi\iota\hbar)$ and assume that $\hbar\in\C\setminus\Q$,
so that $q$ is not a root of unity. Let $\qloop$ be the quantum loop algebra
of $\g$. We refer to \cite[Ch. 12]{chari-pressley} and references therein for
the definition and basic properties of $\qloop$.

% abelian functor
A \fd representation $V\in \Ryang$ is said to be  {\em non--congruent} if,
for every $a,b\in \sigma(V)$, $a-b\not\in \Z_{\neq 0}$. The full subcategory
of such representations is denoted by $\Rync\subset\Ryang$.

In \cite[\S 5]{sachin-valerio-2}, an exact, essentially surjective and faithful functor
\[\sfGamma : \Rync\to\Rloopintro\]
is defined. $\sfGamma$ is such that $\sfGamma(V)=V$ as vector spaces
for any $V\in\Rync$, and restricts to an isomorphism of categories on an
explicit subcategory of $\Rync$ determined by a choice of $\log$. 

% rough form of Gamma
Given $V\in\Rync$, and $i\in\bfI$, consider the additive difference
equation
\[\phi_i(u+1) = \xi_i(u)\phi_i(u)\]
determined by the action of the commuting current $\xi_i(u)$ of $\Yhg$ on
$V$. The action of the commuting current $\psi_i(z)$ of $\qloop$ on $V$ is
given by the monodromy of this equation, that is by
\[\psi_i(z)=\left. 
\lim_{n\to\infty}\xi_i(u+n)\cdots\xi_i(u-n)
\right|_{z=e^{2\pi\iota u}}\]
The action of the raising and lowering operators of
$\qloop$ on $V$ requires the non--congruence hypothesis.
It is not relevant for our current discussion, and we refer to
\cite[\S 5]{sachin-valerio-2} for details.

\subsection{Meromorphic tensor structure on $\sfGamma$ \cite{sachin-valerio-III}}
%----------------------------------------------------------------------------------------------------------

% tensor structure
Let $V_1,V_2\in\Ryang$, and consider the abelian qKZ equation
determined by $\calR^{0,\updown}_{V_1,V_2}(s)$, that is the
difference equation
\begin{equation}\label{eq:ab qKZ}
\Phi(s+1) = \RR^{0,\updown}_{V_1,V_2}(s)\cdot\Phi(s)
\end{equation}
Let $\J_{V_1,V_2}^{\updown}(s):\C\to\End_\C(V_1\otimes V_2)$ be
the left canonical solution of \eqref{eq:ab qKZ}, which is uniquely
determined by the requirement that it be holomorphic and invertible
for $\Re(s)\ll 0$, and possess an asymptotic expansion of the form
$(-s)^{\hbar\Omega}\left(1+O(s^{-1})\right)$ as $s\to\infty$ in any
halfplane of the form $\Re(s)<m$.

The deformed Drinfeld tensor product $\dtensor{\zeta}$ on \fd representations
of $\qloop$ was introduced by D. Hernandez in \cite{hernandez-drinfeld},
and further studied in \cite{sachin-valerio-III}.

\begin{thm}\cite[Thm. 7.3]{sachin-valerio-III}\label{th:drin tens str}
$\J^{\updown}(s)$ is a meromorphic tensor structure on the functor
$\sfGamma$, \wrt the deformed Drinfeld tensor products 
\[\left(\sfGamma,\J^{\updown}(s)\right): 
\left(\Rync,\dtensor{s}\right)
\longrightarrow
\left(\Rloopintro,\dtensor{\zeta}\right)\]
where $\zeta=\exp(2\pi\iota s)$.\footnote{In \cite{sachin-valerio-III},
the {\it right} canonical solution of the equations $\phi(s+1)=\calR^
{0,\updown}_{V_1,V_2}(s)\cdot \phi(s)$ is shown to give rise to a
tensor structure on $\sfGamma$. A similar computation shows that
the left solution yields a tensor structure on a variant of $\sfGamma$,
which we denote $\sfGamma$ for simplicity.}
\end{thm}

\subsection{Tensor structure \wrt the standard coproducts}
%--------------------------------------------------------------------------

The Drinfeld coproduct of $\qloop$ is known to be conjugated to the
standard coproduct by the lower triangular part $\scR^-(\zeta)$ of the
universal $R$--matrix of $\qloop$ (see, for example, \cite{ekp}). Thus,
for any $V_1,V_2\in\Rync$, we have the following isomorphisms of
$\qloop$--modules
\[
\xy
(0,0)*{\sfGamma(V_1)\dtensor{\zeta} \sfGamma(V_2)}="a";
(50,0)*{\sfGamma\lp V_1\dtensor{s} V_2\rp}="b";
(0,30)*{\sfGamma(V_1)\kmtensor{\zeta} \sfGamma(V_2)}="c";
(50,30)*{\sfGamma\lp V_1\kmtensor{s} V_2\rp}="d";
{\ar "c"; "a"};
(-10,15)*{\scriptstyle{\scR^-_{\sfGamma(V_1),
\sfGamma(V_2)}(\zeta)}};
{\ar "d"; "b"};
(60,15)*{\scriptstyle{\RR^-_{V_1,V_2}(s)}};
{\ar "a"; "b"};
(25,-5)*{\scriptstyle{\J^{\updown}_{V_1,V_2}(s)}};
{\ar^{J^{\updown}_{V_1,V_2}(s)} "c"; "d"};
\endxy
\]
where $\V_1\kmtensor{\zeta}\V_2=\tau_\zeta^*\V_1\otimes\V_2$
for any $\V_1,\V_2\in\Rloopintro$, and $J^{\updown}_{V_1,V_2}
(s)$ is defined as the composition 
\begin{equation}\label{eq:composition}
J^{\updown}_{V_1,V_2}(s)
=
\UF_{V_1,V_2}(s)^{-1} \cdot \J^{\updown}_{V_1,V_2}(s)
\cdot \scR^-_{\sfGamma(V_1),\sfGamma(V_2)}
(\zeta)
\end{equation}
Theorem \ref{th:drin tens str} and Proposition \ref{pr:mero tens str} 
therefore imply the following
\begin{cor}\label{cor:tensor-str-KM}
$J^{\updown}_{V_1,V_2}(s)$ is a meromorphic tensor structure
on the functor $\sfGamma$, \wrt the standard tensor products
\[\left(\sfGamma,J^{\updown}(s)\right): 
\left(\Rync,\kmtensor{s}\right)
\longrightarrow
\left(\Rloopintro,\kmtensor{\zeta}\right)\]
\end{cor}

\subsection{Non regularity of $J^{\updown}_{V_1,V_2}(s)$}
%---------------------------------------------------------------------------

Since the tensor products $\kmtensor{s}$ and $\kmtensor
{\zeta}$ are polynomial, the first two authors conjectured in \cite[\S
2.13]{sachin-valerio-III} that $J^{\updown}_{V_1,V_2}(s)$ is regular
and invertible at $s=0$. If so, $J^{\updown}_{V_1,V_2}(0)$ would
give rise to a (non--meromorphic) tensor structure on the functor
$\sfGamma$ \wrt the standard (unshifted) tensor products on $\Rync$ and
$\Rloopintro$. The following shows that this is not the case.

\begin{prop}\label{pr:not-regular}
The meromorphic tensor structure $J^{\updown}(s)$ is either
singular or not invertible at $s=0$.
\end{prop}

\begin{pf}
Since $\sfGamma(V(a)) = 
\sfGamma(V)(e^{2\pi\iota a})$ and
each of the factors in the definition \eqref{eq:composition} of $J^
{\updown}$ is compatible with shifts, we have
\[
J^{\updown}_{V_1(a),V_2(b)}(s)
=J^{\updown}_{V_1,V_2}(s+a-b)
\]
Thus, if $J^{\updown}_{V_1,V_2}(s)$ were regular and invertible at
$s=0$ for every $V_1,V_2$,  then for a fixed $V_1,V_2$ it would be
holomorphic and invertible for all $s$. This cannot be true, as the
following argument shows.

Assume that $V_1, V_2$ are irreducible, with highest weight vectors $v_1,v_2$ of
weights $\lambda_1,\lambda_2\in\h^*$, and that $(\lambda_1,\lambda_2)\neq 0$.
Let $W\subset V_1\otimes V_2$ be the subspace spanned by $v_1\otimes v_2$.
The triangularity of $\UF(s)$ and $\scR^-(\zeta)$ implies that 
\[
\left.J_{V_1,V_2}^{\uparrow}(s)\right|_W=
\left.\J_{V_1,V_2}^{\uparrow}(s)\right|_W
\]
If the restriction of $\J_{V_1,V_2}^{\uparrow}(s)$ to $W$ were holomorphic
and invertible for every $s\in \C$, the same would be true for $\RR_{V_1,V_2}
^{0,\uparrow}(s)$, since by \ref{ssec:mero-tensor-D}
\[
\J_{V_1,V_2}^{\uparrow}(s) = 
\RR_{V_1,V_2}^{0,\uparrow}(s+1)
\cdot
\J_{V_1,V_2}^{\uparrow}(s+1).
\]
In turn, $\left.\A_{V_1,V_2}(s)\right|_W$ would also be holomorphic and
invertible, since (see equation \eqref{eq:diff A} above)
\[
\RR_{V_1,V_2}^{0,\uparrow}(s+\ell\hbar) = 
\A_{V_1,V_2}(s)
\RR_{V_1,V_2}^{0,\uparrow}(s).
\]
Since $\A_{V_1,V_2}(s)$ is a rational function of $s$ such that
$\A_{V_1,V_2}(\infty)=\id$, this implies that $\left.\A_
{V_1,V_2}(s)\right|_W\equiv\id$. The expansion $\ds\A
(s) = 1-\frac{\ell\hbar}{s^2}\Omega_{\h}+O(s^{-3})$, then yields
$\left.\Omega_{\h}\right|_W=0$, which contradicts the fact that 
$\Omega_{\h}v_1\otimes v_2=(\lambda_1,\lambda_2)v_1\otimes v_2$.
\end{pf}

\begin{rem}
The above result leaves open the question of whether there
exists a tensor functor between the (non meromorphic)
tensor categories
\[\left(\Ryang,\otimes
\right)
\to 
\left(\Rloopintro,\otimes
\right)\]
\end{rem}

\subsection{The meromorphic abelian $R$--matrix of $\qloop$} % \cite[\S 8]{sachin-valerio-III}}
%--------------------------------------------------------------------------------

Assume now that $|q|\neq 1$. We review below the analogue of Theorem
\ref{thm:mero-braid-d} for $\qloop$ obtained in \cite{sachin-valerio-III},
based on the results of \cite{damiani}.

Let $\V_1,\V_2$ be two \fd representations of $\qloop$. In \cite[\S 8]
{sachin-valerio-III}, a rational function $\mathscr{A}_{\V_1,\V_2}(\zeta):
\IP^1\to\End_\C(\V_1\otimes\V_2)$ is defined via a contour integral formula involving the
commuting generators of $\qloop$, which is analogous to the one given
in \ref{ssec:A Y}. $\mathscr{A}_{\V_1,\V_2}(\zeta)$ is %a rational function of $\zeta$, 
regular at $\zeta=0,\infty$, and such that
\[\mathscr{A}_{\V_1,V_2}(0)=\id=\mathscr{A}_{\V_1,\V_2}(\infty)\]
Moreover, $[\mathscr{A}_{\V_1,
\V_2}(\zeta),\mathscr{A}_{\V_1,\V_2}(\zeta')]=0$ for any $\zeta,\zeta'$.

Consider the (regular) difference equation with step $q^{2\ellg}$ determined by
$\mathscr{A}_{\V_1,\V_2}(\zeta)$
\begin{equation}\label{eq:Aqdiff}
\ol{\scR}(q^{2\ellg}\zeta) =
\mathscr{A}_{\V_1,\V_2}(\zeta)
\cdot
\ol{\scR}(\zeta)
\end{equation}
This equation admits two meromorphic solutions $\ol{\scR}
^{\uparrow}(\zeta),\ol{\scR}^{\downarrow}(\zeta)$, which are
uniquely determined by the requirement that $\ol{\scR}^{\updown}
(\zeta)$ be holomorphic near $z=q^{\pm\infty}$ and such that $\ol
{\scR}^{\updown}(q^{\pm\infty})=1$. These are explicitly
given by
\[
\ol{\scR}^{\uparrow}(\zeta) = 
\stackrel{\longrightarrow}{\prod_{n\geq 0}}\mathscr{A}_{\V_1,\V_2}(q^{2\ell n}\zeta)^{-1}
\aand
\ol{\scR}^{\downarrow}(\zeta) = 
\stackrel{\longrightarrow}{\prod_{n\geq 1}}\mathscr{A}_{\V_1,\V_2}(q^{-2\ell n}\zeta)
\]

Now set
\[
\scR_{\V_1,\V_2}^{0,\updown}(\zeta)
= \left\{ \begin{array}{ll}
q^{\mp\Omega_\h}\cdot \ol{\scR}^{\updown}(\zeta) & \text{if } |q|<1\\[1.1ex]
q^{\pm\Omega_\h}\cdot \ol{\scR}^{\updown}(\zeta) & \text{if } |q|>1
\end{array}
\right.
\]

\begin{thm}
\cite[\S 8]{sachin-valerio-III}
The category $\lp \Rloopintro, \dtensor{\zeta}, \scR^{0,\updown}(\zeta)\rp$
is a meromorphic braided tensor category.
\end{thm}

% asymptotic statement
\begin{rem}\label{rk:ays R0}
Let $\scR^0$ be the abelian part of the universal $R$--matrix of $\qloop$, and
set
\[\scR^0(\zeta)=\tau_\zeta\otimes\id(\scR^0)\in\hext{\qloop^{\otimes 2}}{\zeta}\]
It is easy to see that $\scR^0_{\V_1,\V_2}(\zeta)$ satisfies the difference
equation \eqref{eq:Aqdiff} \cite[\S 8]{sachin-valerio-III}. It follows by uniqueness
that $\scR^0_{\V_1,\V_2}(\zeta)$ is the Taylor expansion at $\zeta=0$ of
$\scR^{0,\uparrow}_{\V_1,\V_2}(\zeta)$ if $|q|<1$, and of $\scR^{0,\downarrow}
_{\V_1.\V_2}(\zeta)$ if $|q|>1$. Similarly, if
\[\scR^{0,\vee}(\zeta)=\id\otimes\tau_\zeta(\scR^0)\in\hext{\qloop^{\otimes 2}}
{\zeta^{-1}}\]
then $\scR^{0,\vee}_{\V_1,\V_2}(\zeta)$ is the Taylor expansion at $\zeta=\infty$
of $\scR^{0,\downarrow}_{\V_1,\V_2}(\zeta)$ if $|q|<1$, and of $\scR^{0,\uparrow}
_{\V_1.\V_2}(\zeta)$ if $|q|>1$.
\end{rem}

\subsection{Abelian $q$\KD theorem}\label{ss:ab qKZ}
%-----------------------------------------------

Let now $V_1,V_2\in\Ryang$, and consider the abelian $q$KZ equations
\eqref{eq:ab qKZ} determined by $\calR^{0,\updown}_{V_1,V_2}(s)$.

\begin{thm}\cite[Thm. 9.3]{sachin-valerio-III}
\label{th:ab qKD}
The monodromy of the abelian $q$KZ equations \eqref{eq:ab qKZ}  is
equal to $\scR^{0,\updown}_{\sfGamma(V_1),\sfGamma(V_2)}(\zeta)$.
In other words, the following holds
\begin{equation}\label{eq:abqKD}
\scR^{0,\updown}
_{\sfGamma(V_1),\sfGamma(V_2)}(\zeta)\\
= 
\left.
\lim_{n\to\infty}
\RR^{0,\updown}_{V_1,V_2}(s+n)
\cdots
\RR^{0,\updown}_{V_1,V_2}(s)
\cdots
\RR^{0,\updown}_{V_1,V_2}(s-n)
\right|_{\zeta=\exp(2\pi\iota s)}
\end{equation}
\end{thm}
In \cite[\S 9.6]{sachin-valerio-III}, this assertion is strengthened to
include the abelian $q$KZ equations with values in $V_1\otimes
\cdots\otimes V_n$, where $V_i\in\Rync$ and $n\geq 3$. Thus,
Theorem \ref{th:ab qKD} is an analogue of the \KD theorem for
the abelian $q$KZ equations determined by $\RR^{0,
\updown}(s)$.

As is the case for the \KD theorem, Theorems \ref{th:drin tens str}
and \ref{th:ab qKD} can be understood as defining a meromorphic
{\it braided} tensor functor
\[\left(\Rync,\dtensor{s},\calR^{0,\updown}(s)\right)
\longrightarrow
\left(\Rloopintro,\dtensor{\zeta},\scR^{0,\updown}(\zeta)\right)\]
Unlike its non--meromorphic counterpart, this notion involves an
ordered {\it pair} $(\K,\check{\K})$ of meromorphic tensor structures
on the functor $\sfGamma$, such that the following holds \cite
[Rem. 9.3]{sachin-valerio-III}\footnote{For the meromorphic
braided tensor structures discussed in \ref{ss:mero tens str},
$\check{\K}=\K$.}
\begin{equation}\label{eq:twist eq}
\scR^{0,\updown}_{\sfGamma(V_1),\sfGamma(V_2)}(\zeta)=
{\check{\K}_{V_2,V_1}(-s)}^{-1}_{21}
\cdot
\calR^{0,\updown}_{V_1,V_2}(s)
\cdot
\K_{V_1,V_2}(s)
\end{equation}
Comparing \eqref{eq:twist eq} with \eqref{eq:abqKD}, and using
the definition of $\J_{V_1,V_2}$ given in \ref{ssec:mero-tensor-D},
we see that one can take $\K_{V_1,V_2}=\J^{\updown}_{V_1,V_2}$,
which is a regularisation of the product
\[\calR^{0,\updown}_{V_1,V_2}(s-1)\cdot\calR^{0,\updown}_{V_1,V_2}(s-2)\cdots\]
and $\check{\K}_{V_1,V_2}=\J^{\downup}_{V_1,V_2}$, which is a regularisation of
\[\calR^{0,\downup}_{V_2,V_1}(s-1)\cdot\calR^{0,\downup}_{V_1,V_2}(s-2)\cdots
=
{\calR^{0,\updown}_{V_2,V_1}(-s+1)}^{-1}_{21}\cdot{\calR^{0,\updown}_{V_2,V_1}(-s+2)}^{-1}_{21}\cdots
\]
where we used the unitarity relation $\calR^{0,\downup}_{V_1,V_2}(s)
={\calR^{0,\updown}_{V_2,V_1}(s)}^{-1}_{21}$.

\subsection{Meromorphic braided tensor equivalence for standard coproducts}
%A non abelian $q$\KD theorem}
%-----------------------------------------------------------------------------------------------------

% Recap on universal R-matrix of \qloop
Let $\scR$ be the universal $R$--matrix of $\qloop$, $\scR=\scR
^+\cdot\scR^0\cdot\scR^-$ its Gauss decomposition, set $\scR
(\zeta)=\tau_\zeta\otimes\id(\scR)$, and $\scR^\pm(\zeta)=\tau
_\zeta\otimes\id(\scR^\pm)$. Then, if $\V_1,\V_2\in\Rloopintro$,
$\scR^\pm_{\V_1,\V_2}(\zeta)$ are rational functions of $\zeta$
such that
\[{\scR^+_{\V_1,\V_2}(\zeta)}=
{\scR^-_{\V_2,\V_1}(\zeta^{-1})}^{-1}_{21}\]

% Meromorphic R-matrix of \qloop
Define the meromorphic $R$--matrix of $\qloop$ by
\begin{equation}\label{eq:mero full R}
\scR^{\updown}_{\V_1,\V_2}(\zeta)=
\scR^{+}_{\V_1,\V_2}(\zeta)\cdot
\scR^{0,\updown}_{\V_1,\V_2}(\zeta)\cdot
\scR^{-}_{\V_1,\V_2}(\zeta)
\end{equation}
By Remark \ref{rk:ays R0}, $\scR_{\V_1,\V_2}(\zeta)$ is the
Taylor expansion at $\zeta=0$ of $\scR^{\uparrow}_{\V_1,\V_2}(\zeta)$
if $|q|<1$, and of $\scR^{\downarrow}_{\V_1,\V_2}(\zeta)$ if $|q|>1$.

\begin{prop}%\label{cor:tensor-str-KM}
The pair $\left(J^{\updown}_{V_1,V_2}(s),J^{\downup}_{V_1,V_2}(s)\right)$
is a meromorphic braided tensor structure on the functor $\sfGamma$
\wrt the standard tensor products and meromorphic $R$--matrices
\[\left(\Rync,\kmtensor{s},\RR^{\updown}(s)\right)
\longrightarrow
\left(\Rloopintro,\kmtensor{\zeta},\scR^{\updown}(\zeta)\right)\]
\end{prop}
\begin{pf}
By Corollary \ref{cor:tensor-str-KM}, and \ref{ss:ab qKZ}, we only need
to check the compatibility of $\left(J^{\updown}_{V_1,V_2}(s),J^{\downup}
_{V_1,V_2}(s)\right)$ with the meromorphic braidings, that is the relation
\begin{equation}\label{eq:non ab qKD}
\scR^{\updown}_{\sfGamma(V_1),
\sfGamma(V_2)}(\zeta)
=
{J^{\downup}_{V_2,V_1}(-s)}^{-1}_{21}
\cdot
\RR^{\updown}_{V_1,V_2}(s)
\cdot
J^{\updown}_{V_1,V_2}(s)
\end{equation}

The Gauss decomposition  \eqref{eq:mero full R} yields
\[\begin{split}
\scR^{\updown}_{\sfGamma(V_1),
\sfGamma(V_2)}(\zeta)
=&\,
\scR^{+}_{\sfGamma(V_1),\sfGamma(V_2)}(\zeta)\cdot
\scR^{0,\updown}_{\sfGamma(V_1),\sfGamma(V_2)}(\zeta)\cdot
\scR^{-}_{\sfGamma(V_1),\sfGamma(V_2)}(\zeta)\\
% by Gauss
=&\,
\scR^{+}_{\sfGamma(V_1),\sfGamma(V_2)}(\zeta)
\cdot
{\J^{\downup}_{V_2,V_1}(-s)}^{-1}_{21}
\cdot
\calR^{0,\updown}_{V_1,V_2}(s)
\cdot
\J^{\updown}_{V_1,V_2}(s)
\cdot
\scR^{-}_{\sfGamma(V_1),\sfGamma(V_2)}(\zeta)\\
% by abelian qDK theorem
=&\,
\scR^{+}_{\sfGamma(V_1),\sfGamma(V_2)}(\zeta)
\cdot
{\J^{\downup}_{V_2,V_1}(-s)}^{-1}_{21}
\cdot
\RR^+_{V_1,V_2}(s)^{-1}
\cdot
\RR^{\updown}_{V_1,V_2}(s)
\\
&\cdot
\UF_{V_1,V_2}(s)^{-1}
\cdot
\J^{\updown}_{V_1,V_2}(s)
\cdot
\scR^{-}_{\sfGamma(V_1),\sfGamma(V_2)}(\zeta)\\
% by definition of full R for Y
=&\,
{J^{\downup}_{V_2,V_1}(-s)}^{-1}_{21}
\cdot
\RR^{\updown}_{V_1,V_2}(s)
\cdot
J^{\updown}_{V_1,V_2}(s)
\end{split}
% by definition of J
\]
where the second equality follows from the twist relation \eqref
{eq:twist eq}, the third from the definition of $\RR_{V_1,V_2}^
{\updown}(s)$ given in \ref{ss:fullR}, and the fourth from the
definition \eqref{eq:composition} of $J^{\updown}_{V_1,V_2}(s)$,
together with the unitarity relations $\scR^+_{\V_1,\V_2}(\zeta)
={\scR^-_{\V_2,\V_1}(\zeta^{-1})}^{-1}_{21}$ and $\RR^+_{V_1,
V_2}(s)={\RR^-_{V_2,V_1}(-s)}^{-1}_{21}$.
\end{pf}
%\newpage 

\begin{rem}\label{rm:depressing}
Much like \eqref{eq:abqKD}, the twist relation \eqref{eq:non ab qKD}
can be regarded as a monodromy relation. Indeed, both
\begin{align*}
J^{\updown}_{V_1,V_2}(s)
&=
\UF_{V_1,V_2}(s)^{-1}
\cdot
\J^{\updown}_{V_1,V_2}(s)
\cdot
\scR^-_{\sfGamma(V_1),\sfGamma(V_2)}
(\zeta)\\
\intertext{and}
\left({J^{\downup}_{V_2,V_1}(-s)}^{-1}_{21}
\cdot
\RR^{\updown}_{V_1,V_2}(s)\right)^{-1}
&=
\UF_{V_1,V_2}(s)^{-1}
\cdot
\RR^{0,\updown}_{V_1,V_2}(s)^{-1}
\cdot
\J^{\downup}_{V_2,V_1}(-s)_{21}\\
&\qquad \cdot 
\scR^-_{\sfGamma(V_2),\sfGamma(V_1)}(\zeta^{-1})_{21}\\
&=
\UF_{V_1,V_2}(s)^{-1}
\cdot
\RR^{0,\downup}_{V_2,V_1}(-s)_{21}
\cdot
\J^{\downup}_{V_2,V_1}(-s)_{21}\\
&\qquad \cdot 
\scR^-_{\sfGamma(V_2),\sfGamma(V_1)}(\zeta^{-1})_{21}\\
\end{align*}
are solutions of the difference equation
\[\Phi(s+1)=
\left(\UF(s+1)_{V_1,V_2}^{-1}\cdot \RR^{0,\updown}_{V_1,V_2}(s)\cdot\UF(s)_{V_1,V_2}\right)
\cdot\Phi(s)\]
though, due to the presence of the factors $\scR^-_{\sfGamma(V_1),\sfGamma(V_2)}
(\zeta)$ and $\scR^-_{\sfGamma(V_2),\sfGamma(V_1)}(\zeta^{-1})_{21}$, neither is
a canonical left or right solution. Unlike \eqref{eq:abqKD}, however, the twist relation
\eqref{eq:non ab qKD} should not be considered as a difference version of the (non--abelian)
\KD Theorem on $n=2$ points for several reasons.
\begin{enumerate}
% not the right equations
\item As just pointed out, the difference equations underlying \eqref{eq:non ab qKD}
are not the $q$KZ equations determined by $\RR^{\updown}_{V_1,V_2}(s)$, but (a
gauge transformation of) the abelian $q$KZ equations determined by $\RR^{0,\updown}
_{V_1,V_2}(s)$.
% no dynamical parameter
\item The $q$KZ equations of Frenkel--Reshetikhin \cite{frenkel-reshetikhin}
include a {\em dynamical parameter} $\lambda\in\h$, and
are given by
\[
\Phi(s+1) = e^{\lambda}\otimes \id\cdot \RR^{\updown}_{V_1,V_2}(s) \cdot \Phi(s)
\]
Since the form of the asymptotics of solutions of this equation depends
on the regularity of $e^\lambda$, its monodromy is a meromorphic function
of $\lambda$, which is conjectured to be equivalent to $\scR^{\varepsilon}
_{\sfGamma(V_1)\sfGamma(V_2)}(\zeta)$, via a $\lambda$--dependent
gauge transformation.
\end{enumerate}
\end{rem}

\appendix

\section{Separation of points}\label{asec:sep}
%=====================

\subsection{}\label{assec:sep-state}
%--------------

Let $\V$ be a collection of \fd representations of $\Yhg$ such that
\begin{enumerate}
\item\label{V:C1} $\C\in \mathcal{V}$, and 
     $V_1 \otimes V_2\in \mathcal{V}$ for all $V_1,V_2\in \mathcal{V}$,
\item\label{V:C2} $\exists$ $V\in \mathcal{V}$ 
     such that $\Ker({\pi_V}|_\g)=\{0\}$. 
\end{enumerate}
The goal of this section is to prove that the elements of $\V$ separate
points in $\Yhg$. More generally, the following holds.
\begin{prop}\label{P:sep} %Fix $n\in \Z_{>0}$, l
Let $\V_1,\ldots,\V_n$ be collections of \fd representations of $\Yhg$
satisfying the conditions (1) and (2) above. Then, the ideal $\J_n\subset \Yhg^{\otimes n}$
defined by  
\begin{equation*}
\J_n= 
\bigcap_{V_i\in\V_i} \Ker(\pi_{V_1}\otimes \cdots \otimes \pi_{V_n})
\end{equation*}
is the zero ideal.
\end{prop}
\begin{rem}
The analogous statement for $U(\g[z])$ fails. Indeed, take $n=1$ and
let  $V_\g$ be any faithful, finite-dimensional $\g$-module. Set $V=
\mathrm{ev}^*(V_\g)$, where $\mathrm{ev}$ is the evaluation morphism 
\begin{equation*}
\mathrm{ev}: \g[z]\to \g, \quad f(z)\to f(0). %\quad \forall \; f(z)\in \g[z]. 
\end{equation*}
Then, $V$ is a $U(\g[z])$-module satisfying \eqref{V:C2}, and $\V=\{V
^{\otimes n}\}_{n\in \N}$ satisfies \eqref{V:C1}--\eqref{V:C2}. However, 
\begin{equation*}
z\g[z]\subset \bigcap_{k\geq 0} \Ker(\pi_{V^{\otimes k}})
\end{equation*}
\end{rem}

The proof of the proposition is given in Sections \ref{assec:sep-red}--\ref
{assec:sep-J}. In Section  \ref{assec:sep-red}, we reduce the task to proving
that $\J_1=\{0\}$. The reduction step is elementary, but is included
for the sake of completeness. That the ideal $\J_1\subset \Yhg$
vanishes is
%follows
%from 
an unpublished result of V. Drinfeld's, whose proof in the $\hbar$--formal setting has
been reproduced in \cite[Prop.~8.8]{sachin-valerio-1}. After recalling relevant background
material on co-Poisson Hopf algebras and Lie bialgebras in Sections \ref{assec:sep-coP}
and \ref{assec:sep-quant}, we explain in Section \ref{assec:sep-J} how to modify the
argument given in \cite{sachin-valerio-1} to deduce that $\J_1=\{0\}$.

\subsection{Reduction step}\label{assec:sep-red}
%----------------------------------

Let $\mathcal{H}_1$ and $\mathcal{H}_2$ be associative
algebras over $\C$. Assume in addition that, 
for each $i\in\{1,2\}$, $\mathcal{H}_i$ is equipped 
with a family of representations $\mathcal{V}_i$ such that 
\begin{equation*}
\J_{\mathcal{V}_i}:=
\bigcap_{V\in \mathcal{V}_i}\Ker(\pi_V)=\{0\}
\end{equation*}

The following general result, coupled with a simple induction on $n$, 
implies that Proposition \ref{P:sep} holds,
provided  $\J_1=\{0\}$. 
\begin{lem} We have 
\begin{equation*}
\bigcap_{V_i\in \mathcal{V}_i}\Ker(\pi_{V_1}\otimes \pi_{V_2})=\{0\}
\end{equation*}
\end{lem}
\begin{pf} The lemma follows easily from the the identity
\[
\Ker(\pi_{V_1}\otimes \pi_{V_2})
=(\id\otimes \pi_{V_2})^{-1}\big(\Ker(\pi_{V_1})\otimes \End(V_2)\big),
\]
the assumption $\J_{V_i}=\{0\}$, and the fact that, 
for any vector spaces $M$, $N$ and collection of subspaces 
$\{M_\lambda\}_{\lambda\in \Lambda}\subset M$, 
we have the following equality in $M\otimes N$:
\begin{equation*}
\bigcap_{\lambda\in \Lambda} (M_\lambda\otimes N) 
=\left(\bigcap_{\lambda\in \Lambda}M_\lambda \right)\otimes N
\qedhere
\end{equation*}
\end{pf}
\subsection{}\label{assec:sep-coP}
%-----------
We now pause to collect pertinent facts from the 
theories of co-Poisson Hopf algebras and Lie bialgebras. 
Fix a Lie algebra $\mathfrak{a}$ over $\C$, and recall 
that a co-Poisson algebra structure on $U(\mathfrak{a})$ 
is given by the additional data of a 
Poisson cobracket $\delta$, i.e. a linear map 
\[
\delta:U(\mathfrak{a})\to U(\mathfrak{a})\wedge U(\mathfrak{a}) 
\subset U(\mathfrak{a})^{{\otimes 2}}
\]
satisfying the co-Jacobi and co-Leibniz identities: 
\begin{gather*}
(\id+(1\,2\,3)+(1\,3\,2))\circ\delta\otimes\id\circ\delta=0\\
\id\otimes\Delta\circ\delta
=\delta\otimes\id\circ\Delta+(1\,2)\circ\id\otimes\delta\circ\Delta
\end{gather*}
If in addition $\delta$ satisfies the compatibility condition 
\begin{equation*}
\delta(xy)=\delta(x)\Delta(y)+\Delta(x)\delta(y) \quad 
%\forall \; x,y\in U(\mathfrak{a}), 
\end{equation*} 
then $(U(\mathfrak{a}),\delta)$ is said to be a co-Poisson 
Hopf algebra. In this case, the cobracket $\delta$ is uniquely
determined by its restriction $\delta|_\mathfrak{a}$ to 
$\mathfrak{a}$, which can be shown to define a Lie bialgebra 
structure on $\mathfrak{a}$. Conversely, every Lie bialgebra 
cocommutator  on $\mathfrak{a}$ uniquely extends to a 
co-Poisson Hopf cobracket on $U(\mathfrak{a})$: 
see \cite[Prop.~6.2.3]{chari-pressley}.

\medskip
 
 A perhaps less well-known correspondence, which we shall exploit, 
 takes place at the level of ideals. Recall that a co-Poisson 
 bialgebra ideal of $(U(\mathfrak{a}),\delta)$ is a two-sided 
 ideal $\mathrm{J}$ of $U(\mathfrak{a})$ satisfying the coideal 
 and co-Poisson conditions
\begin{gather*}
\Delta(\mathrm{J})\subset \mathrm{J}\otimes U(\mathfrak{a}) 
+ U(\mathfrak{a})\otimes\mathrm{J},\quad \epsilon(\mathrm{J})=0 \\
\delta(\mathrm{J})
\subset \mathrm{J}\otimes U(\mathfrak{a}) + U(\mathfrak{a})\otimes\mathrm{J}
\end{gather*}
 where $\epsilon$ is the counit for $U(\mathfrak{a})$. 
 Similarly, a Lie bialgebra ideal of $(\mathfrak{a},\delta|_\mathfrak{a})$ 
 is a Lie algebra ideal $\mathfrak{l}$ of $\mathfrak{a}$ satisfying 
\begin{equation*}
\delta(\mathfrak{l})\subset \mathfrak{l}\otimes\mathfrak{a}
+\mathfrak{a}\otimes \mathfrak{l}. 
\end{equation*} 
\begin{lem}\label{L:Hopfi} Fix a co-Poisson Hopf algebra 
structure on $U(\mathfrak{a})$. Then the assignment
\begin{equation}\label{J-bij:1}
\mathrm{J}\subset U(\mathfrak{a}) \mapsto 
\mathrm{J}\cap \mathfrak{a} \subset \mathfrak{a}
\end{equation}
determines a bijective correspondence between co-Poisson
 bialgebra ideals $\mathrm{J}$ of $U(\mathfrak{a})$ and 
 Lie bialgebra ideals $\mathfrak{l}\subset \mathfrak{a}$, 
 with inverse 
\begin{equation}\label{J-bij:2}
\mathfrak{l}\subset \mathfrak{a}\mapsto 
\langle \mathfrak{l}\rangle \subset U(\mathfrak{a})
\end{equation}
where $\langle \mathfrak{l}\rangle$ is the two-sided 
ideal generated by $\mathfrak{l}$. 
\end{lem}
When the underlying cobracket is taken to be trivial 
(that is, $\delta\equiv 0$), this reduces to the more 
familiar assertion (see \cite[Prop.~4.8]{Pass14}, for example) 
that \eqref{J-bij:1} determines a bijective correspondence
 between bialgebra ideals of $U(\mathfrak{a})$ and Lie 
 algebra ideals in $\mathfrak{a}$, with inverse \eqref{J-bij:2}. 
The lemma follows readily from this special case by a 
straightforward verification that \eqref{J-bij:1} and 
\eqref{J-bij:2} preserve any additional co-Poisson structure.
\subsection{}\label{assec:sep-quant}
%-----------
We now narrow our focus to $\mathfrak{a}=\g[z]$. Let us begin by
recalling how one passes from $\Yhg$ to the standard Lie bialgebra 
structure on $\g[z]$, given by \eqref{Yhg-LB} below.
Since  $\Yhg$ is a filtered Hopf deformation of the 
cocommutative Hopf algebra $U(\g[z])$ (see Section \ref{ssec:filt-tau}), 
the linear map
\begin{gather*}
\Delta_{\mathrm{Alt}}:=
\Delta-\Delta^{\scriptscriptstyle{\operatorname{op}}}: 
\Yhg\to \Yhg\otimes \Yhg
\end{gather*}
is a filtered linear map of degree $-1$. That is, it satisfies 
\[
\Delta_{\mathrm{Alt}}(\filt{k}{\Yhg})
\subset \filt{k-1}{\Yhg^{\otimes 2}}\quad \forall \; k\in \N
\]
where $\filt{-n}{\Yhg^{\otimes 2}}=\{0\}$ for $n>0$. We may 
therefore view it  as a filtered map 
\[
(\Yhg,\filt{\bullet}{\Yhg})
\to (\Yhg^{\otimes 2}, \mathrm{F}_\bullet(\Yhg^{\otimes 2}))
\]
where $\mathrm{F}_\bullet(\Yhg^{\otimes 2})$ is the vector space 
filtration on $\Yhg^{\otimes 2}$ defined by 
$
\mathrm{F}_k(\Yhg^{\otimes 2}):= \filt{k-1}{\Yhg^{\otimes 2}}$ 
for all 
$k\in \N$. 
By definition, the semiclassical limit of $\Delta$ is the 
associated graded map
\begin{equation} \label{Yhg-SCL}
\delta:=\gr\!\left(\frac{\Delta_{\mathrm{Alt}}}{\hbar}\right):
U(\g[z])\to \gr_\mathrm{F}(\Yhg\otimes \Yhg)\cong U(\g[z])\otimes U(\g[w])
\end{equation}
It is a Poisson cobracket which endows $U(\g[z])$ with the 
structure of a co-Poisson Hopf algebra. Moreover, this co-Poisson 
structure induces the standard Lie bialgebra structure on $\g[z]$,
 with cocommutator
  \begin{equation*}
  \delta|_{\g[z]}:\g[z]\to \g[z]\otimes \g[w]\cong (\g\otimes\g)[z,w]
  \end{equation*}
  given on $f(z)\in \g[z]$ by the formula
\begin{equation}\label{Yhg-LB}
\delta(f(z)) = \left[f(z)\otimes 1 + 1\otimes f(w),
\frac{\Omega_{\g}}{z-w} \right] \in (\g\otimes\g)[z,w]
\end{equation}
We may summarize the above discussion concisely by saying that 
$\Yhg$ is a \textit{filtered quantization} of the Lie bialgebra 
$\g[z]$, equipped with the above cocommutator.\\

The last ingredient we shall need is the following lemma,
 which is immediately obtained from Corollary 8.9 of 
 \cite{sachin-valerio-1} with the help of Lemma \ref{L:Hopfi}.
\begin{lem}\label{L:Hopf-g[z]} Let $\epsilon_U:U(\g[z])\to \C$ be the counit. 
\begin{enumerate}[font=\upshape]
\item If $\mathfrak{l}\subset \g[z]$ is a Lie bialgebra ideal, 
      then $\mathfrak{l}=\{0\}$ or $\mathfrak{l}=\g[z]$. 
\item If $\mathrm{J}\subset U(\g[z])$ is a co-Poisson bialgebra
       ideal, then $\mathrm{J}=\{0\}$ or $\mathrm{J}=\Ker(\epsilon_U)$. 
\end{enumerate}
\end{lem}

\subsection{Proof that $\J$ vanishes}\label{assec:sep-J}
%-----------

Consider now the filtration $\filt{\bullet}{\J}$ on $\J$  
given by $\filt{k}{\J}=\filt{k}{\Yhg}\cap \J$ 
for all $k\in \N$, and the associated graded ideal 
\begin{equation*}
\gr(\J)
=\bigoplus_{k\geq 0} \filt{k}{\J}/\filt{k-1}{\J}
\subset \gr(\Yhg)=U(\g[z])
\end{equation*}
If $x\in \J$ is nonzero and $k\in \N$ is minimal such 
that $x\in \filt{k}{\J}$, then the image of $x$ in 
$\gr_k(\J)$ is a nonzero element. 
Hence, our task is reduced to proving that $\gr(\J)=\{0\}$.  

Next, note that the condition \eqref{V:C1} imposed on $\mathcal{V}$ 
guarantees that $\J$ is a bialgebra ideal in $\Yhg$. 
As $\Delta$ and $\epsilon$ are filtered, it follows that 
$\gr(\J)$ is necessarily a co-Poisson bialgebra ideal of 
$U(\g[z])$. Using Lemma \ref{L:Hopf-g[z]}, we deduce that  
$\gr(\J)=\{0\}$ or $\gr(\J)=\Ker(\epsilon_U)$. 
If $\gr(\J)=\Ker(\epsilon_U)$, then $\g\subset \gr(\J)$,
 and thus 
\[\g\subset \gr_0(\J)=\filt{0}{\J}\subset \J\]
This contradicts the assumption \eqref{V:C2} on $\mathcal{V}$, 
which guarantees that  $\J$ intersects $\g$ trivially. 
Hence, we may conclude that $\gr(\J)$, and thus $\J$ 
itself, vanishes. 

\section{Uniqueness of the universal $R$--matrix}\label{asec:UR-unique}
%=================================================

The aim of this section is to give a proof of the
uniqueness part of Drinfeld's theorem (Theorem \ref{thm:intro1}).
Namely, we assume that we are given two formal series
\[
\RR^{(1)}(s), \RR^{(2)}(s)
\in 1 + s^{-1}\hext{(\Yhg\otimes\Yhg)}{s^{-1}}
\]
satisfying the hypotheses of Theorem \ref{thm:intro1}. That is,
for $i=1,2$,
\begin{align*}
\Delta\otimes \id(\RR^{(i)}(s)) &= \RR^{(i)}_{13}(s)\RR^{(i)}_{23}(s)\\
\id \otimes\Delta(\RR^{(i)}(s)) &= \RR^{(i)}_{13}(s)\RR^{(i)}_{12}(s)
\end{align*}
and, for every $a\in\Yhg$: 
\[
\tau_s\otimes\id\circ \Delta\sop(a)
= 
\RR^{(i)}(s)\cdot
\tau_s\otimes\id\circ \Delta(a)
\cdot
\RR^{(i)}(s)^{-1}.
\]

%In this section $\Delta = \kmDelta$. 
Our main tool will be the
following lemma.

\subsection{}
%-----------

\begin{lem}\label{lem:primitiveYKM}
The Lie algebra of primitive elements
 \begin{equation*}
 \Primitive{\Yhg}{\Delta}=\{y\in \Yhg\,:\,\Delta(y)=y\otimes 1+1\otimes y\}
 \end{equation*}
 is equal to $\g$.
 \end{lem}

\begin{pf}
  We shall again exploit the fact, recalled in Section \ref{assec:sep-quant}, 
  that the Yangian $\Yhg$ provides a filtered quantization of the Lie bialgebra 
  structure $(\g[z],\delta)$ on $\g[z]$ given by \eqref{Yhg-LB}. 
  By \eqref{Yhg-SCL}, this means that, for each $x\in \g$ and $k\geq 0$, we have
 \begin{equation}\label{Delta-delta}
  \hbar\cdot\delta(x.z^k)=\Delta(y)-\Delta^{\scriptscriptstyle{\operatorname{op}}}(y) \mod \mathcal{F}_{k-2}(\Yhg\otimes \Yhg )
 \end{equation}
 for any $y\in \mathcal{F}_{k}(\Yhg)$ whose image $\bar y\in \gr_{k}(\Yhg)\subset U(\g[z])$ coincides with $x.z^k$. \\

Now let $y\in\Yhg$ be an arbitrary nonzero primitive element. Assume that
$k\geq 0$ is such that $y\in\filt{k}{\Yhg}\setminus
\filt{k-1}{\Yhg}$. As $\Delta$ is filtered with $\gr(\Delta)$ recovering the standard coproduct on $U(\g[z])$ (see Section \ref{ssec:filt-tau}), we can conclude that the image $\bar y$ of $y$ in $\gr_k(\Yhg)\subset U(\g[z])$ is a nonzero, primitive degree $k$ element, and thus
of the form $\bar y = x.z^k$ for some $x\in\g$.\\

Using that $\Delta(y) = \Delta^{\scriptscriptstyle{\operatorname{op}}}(y)$, we deduce from \eqref{Delta-delta} that
$\delta(x.z^k)=0$.
On the other hand, by definition of $\delta$, we have
\[
0 = \delta(x.z^k) = \frac{z^k-w^k}{z-w} [x\otimes 1, \Omega_{\g}]
\]
Hence, $k=0$ and $y\in\g\subset\Yhg$.
\end{pf}

\subsection{}
%--------------

Now let $n\geq 1$ and $X\in\Yhg\otimes\Yhg$ be such that
$\RR^{(1)}(s)-\RR^{(2)}(s) = s^{-n}X + O(s^{-n-1})$.
We will prove that $X=0$.\\

Comparing the coefficients of $s^{-n}$ on both sides
of the cabling identities, we obtain the following
\begin{align*}
\Delta\otimes \id (X) &= X_{13} + X_{23} \\
\id \otimes \Delta (X) &= X_{13} + X_{12}
\end{align*}
In other words, $X\in\Primitive{\Yhg}{\Delta}^{\otimes 2}$,
that is, $X\in \g\otimes\g \subset \Yhg\otimes\Yhg$, by the lemma above.\\

By the intertwining equation for $a\in\g$, we conclude
that $X\in (\g\otimes\g)^{\g}$.
Hence $X$ is a scalar
multiple of the Casimir tensor: $X = c\Omega_{\g}$,
for some $c\in\C$.\\

Let us now consider the intertwining equation
for $a = \Top{h}$
\[
\ad(\primitive{\Top{h}} + sh\otimes 1)\cdot
\RR^{(i)}(s)
=\hbar\fopop(h)\RR^{(i)}(s)
+\hbar\RR^{(i)} \fop(h)
\]
Take the difference of the two equations, for $i=1,2$,
and compare the coefficient of $s^{-n+1}$, to get
$c[h\otimes 1,\Omega_{\g}]=0$, for every $h\in\h$.
But that means $c=0$ and hence $X=0$, which is 
what we wanted to show.

\bibliographystyle{amsplain}
%\bibliography{LambOfGod,Hopf,Yangians}

 \newcommand{\noop}[1]{}
\providecommand{\bysame}{\leavevmode\hbox to3em{\hrulefill}\thinspace}
\providecommand{\MR}{\relax\ifhmode\unskip\space\fi MR }
\providecommand{\MRhref}[2]{%
  \href{http://www.ams.org/mathscinet-getitem?mr=#1}{#2}
}
\providecommand{\href}[2]{#2}

\end{document}